\documentclass[12pt]{article}
\usepackage{amsmath, amsthm, amssymb, stmaryrd, fullpage}
\usepackage[all]{xy}

\numberwithin{equation}{subsection}
\newtheorem{theorem}[equation]{Theorem}
\newtheorem{lemma}[equation]{Lemma}

\newtheorem{cor}[equation]{Corollary}
\newtheorem{prop}[equation]{Proposition}

\theoremstyle{definition}
\newtheorem{defn}[equation]{Definition}
\newtheorem{example}[equation]{Example}
\newtheorem{remark}[equation]{Remark}

\newtheorem{convention}[equation]{Convention}

\def\CC{\mathbb{C}}
\def\FF{\mathbb{F}}

\def\QQ{\mathbb{Q}}
\def\RR{\mathbb{R}}
\def\ZZ{\mathbb{Z}}

\def\calO{\mathcal{O}}

\def\calR{\mathcal{R}}

\newcommand{\gothm}{\mathfrak{m}}

\def\alg{\mathrm{alg}}
\def\an{\mathrm{an}}
\def\con{\mathrm{con}}
\def\sep{\mathrm{sep}}
\def\perf{\mathrm{perf}}

\def\beq{\begin{equation}}
\def\eeq{\end{equation}}

\def\dual{\vee}

\def\fp{\pi^{-1}}

\def\GK{\Gamma^K}
\def\Gperf{\Gamma^{\perf}}

\def\Galg{\Gamma^{\alg}}

\def\be{\mathbf{e}}

\def\bv{\mathbf{v}}
\def\bw{\mathbf{w}}

\def\Galgcon{\Gamma^{\alg}_{\con}}
\def\Galgancon{\Gamma^{\alg}_{\an,\con}}
\def\Gperfcon{\Gamma^{\perf}_{\con}}
\def\Gancon{\Gamma_{\an,\con}}
\def\GKcon{\Gamma^K_{\con}}
\def\GKancon{\Gamma^K_{\an,\con}}
\def\GLcon{\Gamma^L_{\con}}
\def\GLancon{\Gamma^L_{\an,\con}}
\def\Gcon{\Gamma_{\con}}

\DeclareMathOperator{\End}{End}
\DeclareMathOperator{\Ext}{Ext}

\DeclareMathOperator{\Frac}{Frac}
\DeclareMathOperator{\Gal}{Gal}

\DeclareMathOperator{\GL}{GL}
\DeclareMathOperator{\Hom}{Hom}

\DeclareMathOperator{\image}{im}

\DeclareMathOperator{\naive}{naive}

\DeclareMathOperator{\rank}{rank}

\DeclareMathOperator{\Spec}{Spec}

\newcounter{fixmectr}

\begin{document}

\title{Slope filtrations revisited}
\author{Kiran S. Kedlaya \\ Department of Mathematics \\ Massachusetts
Institute of Technology \\ 77 Massachusetts Avenue \\
Cambridge, MA 02139 \\
\texttt{kedlaya@math.mit.edu}}
\date{October 25, 2005}

\maketitle

\begin{abstract}
We give a ``second generation'' exposition of the
slope filtration theorem
for modules with Frobenius action over the Robba ring, providing a number of 
simplifications in the arguments. 
Some of these are inspired by parallel work of Hartl and Pink, which points
out some analogies with the formalism of stable vector bundles.
\end{abstract}

\tableofcontents

\section{Introduction}

This paper revisits the slope filtration theorem given
by the author in \cite{me-local}. Its main purpose is expository: it
provides a simplified and clarified presentation of the theory of slope
filtrations over rings of Robba type. In the process, we generalize
the theorem in a fashion useful for certain applications, such as the
semistable reduction problem for overconvergent $F$-isocrystals
\cite{me-part1}.

In the remainder of this introduction,
we briefly describe the theorem and some applications, then say a bit
more about the nature and structure of this particular paper.

\subsection{The slope filtration theorem}

The Dieudonn\'e-Manin classification 
\cite{katz}, \cite{manin} describes the category of
finite free modules equipped with a Frobenius action,
over a complete discrete valuation ring with algebraically closed
residue field, loosely analogous to the eigenspace decomposition of
a vector space over an algebraically
closed field equipped with a linear transformation. When the residue field is unrestricted, the classification
no longer applies, but one does retrieve a canonical filtration whose
successive quotients are all isotypical of different types if one
applies the Dieudonn\'e-Manin classification after enlarging the residue
field.

The results of \cite{me-local} give an analogous pair of assertions for
finite free modules equipped with a Frobenius action over the
Robba ring over a complete discretely valued field of mixed characteristic.
(The Robba ring consists of those formal Laurent series over the given
coefficient field converging on some open annulus of outer radius 1.)
Namely, over a suitable ``algebraic closure'' of the Robba ring,
every such module admits a decomposition into the same sort of standard
pieces as in Dieudonn\'e-Manin (\cite[Theorem~4.16]{me-local}
and Theorem~\ref{T:Dieudonne-Manin} herein),
and the analogous canonical slope filtration
descends back down to the original module (\cite[Theorem~6.10]{me-local}
and Theorem~\ref{T:descend special} herein).

\subsection{Applications}

The original application of the slope filtration theorem was to
the $p$-adic local monodromy theorem 
on quasi-unipotence of $p$-adic differential equations
with Frobenius structure over the Robba ring. 
(The possibility of, and need for, such a theorem first arose in the work of Crew
\cite{crew1}, \cite{crew2}
on the rigid cohomology of curves with coefficients, 
and so the theorem 
is commonly referred to as ``Crew's conjecture''.)
Specifically, the slope filtration
theorem reduces the $p$LMT to its unit-root case, established
previously by Tsuzuki \cite{tsuzuki-unitroot}. We note, for now
in passing, that Crew's conjecture
has also been proved by Andr\'e \cite{andre} and by Mebkhout 
\cite{mebkhout}, using the Christol-Mebkhout index theory for
$p$-adic differential equations; for more on the relative merits
of these proofs, see Remark~\ref{R:compare pLMT}.

In turn, the $p$-adic local monodromy theorem is already known to have
several applications. Many of these are in the study of rigid
$p$-adic cohomology of varieties over fields of characteristic $p$:
these include
a full faithfulness theorem for restriction between the categories
of overconvergent and convergent $F$-isocrystals \cite{me-full}, a
finiteness theorem with coefficients \cite{me-finite}, and an analogue
of Deligne's ``Weil II'' theorem \cite{me-weil}. The $p$LMT also gives
rise to a proof of Fontaine's conjecture that every de Rham representation (of
the absolute Galois group of a mixed characteristic local field)
is potentially semistable, via a construction of Berger 
\cite{berger-cst} linking the theory of $(\phi,\Gamma)$-modules to 
the theory of $p$-adic differential equations.

Subsequently, other applications of the slope filtration theorem have
come to light. Berger \cite{berger-weak} has used it to give a new proof
of the theorem of Colmez-Fontaine that weakly admissible 
$(\phi, \Gamma)$-modules are admissible. A variant of Berger's proof
has been given by Kisin \cite{kisin}, who goes on to give
a classification of crystalline representations with nonpositive
Hodge-Tate weights in terms of certain Frobenius modules; as corollaries,
he obtains classification results for $p$-divisible groups conjectured
by Breuil and Fontaine.
Colmez \cite{colmez} has used the slope filtration theorem to construct
a category of ``trianguline representations'' involved in a proposed
$p$-adic Langlands correspondence.
Andr\'e and di Vizio 
\cite{andre-divizio} have used the slope filtration theorem to prove
an analogue of Crew's conjecture for $q$-difference equations, by establishing
an analogue of Tsuzuki's theorem for such equations. (The replacement of
differential equations by $q$-difference equations does not affect the
Frobenius structure, so the slope filtration theorem applies unchanged.)
We expect to see additional applications in the future.

\subsection{Purpose of the paper}

The purpose of this paper is to give a ``second generation'' exposition
of the proof of the slope filtration theorem, using ideas we have learned
about since \cite{me-local} was written. These ideas include a close analogy
between the theory of slopes of Frobenius modules and the formalism of 
semistable vector bundles; this analogy is visible in the work of 
Hartl and Pink \cite{hartl-pink}, which strongly 
resembles our Dieudonn\'e-Manin classification but takes place in equal
characteristic $p>0$. It is also visible in the theory of filtered
$(\phi,N)$-modules, used to study $p$-adic Galois representations; indeed,
this theory is directly related to slope filtrations via the work of 
Berger \cite{berger-weak} and Kisin \cite{kisin}.

In addition to clarifying the exposition, we have phrased the results at a
level of generality that may be useful for additional applications. In
particular, the results apply to Frobenius modules over what might be
called ``fake annuli'', which occur in the context of semistable reduction
for overconvergent $F$-isocrystals (a higher-dimensional analogue of
Crew's conjecture). See \cite{me-fake} for an analogue of the $p$-adic
local monodromy theorem in this setting.

\subsection{Structure of the paper}

We conclude this introduction with a summary of the various
chapters of the paper.

In Chapter~\ref{sec:basic rings}, we construct a number of rings
similar to (but more general than)
those occurring in \cite[Chapters~2 and~3]{me-local},
and prove that a
certain class of these are B\'ezout rings (in which every finitely
generated ideal is principal). 

In Chapter~\ref{sec:sigma-mod}, we introduce $\sigma$-modules and
some basic terminology for dealing with them.
Our presentation is informed by some
strongly analogous work (in equal characteristic $p$)
of Hartl and Pink.

In Chapter~\ref{sec:dm}, we give a uniform presentation of the standard
Dieudonn\'e-Manin decomposition theorem and of the variant form proved
in \cite[Chapter~4]{me-local}, again using the Hartl-Pink framework.

In Chapter~\ref{sec:generic special}, we recall some results mostly from
\cite[Chapter~5]{me-local} on $\sigma$-modules over the bounded subrings
of so-called
analytic rings. In particular, we compare the ``generic'' and ``special''
polygons and slope filtrations.

In Chapter~\ref{sec:descents}, we give a streamlined form of the arguments
of \cite[Chapter~6]{me-local}, which deduce from the Dieudonn\'e-Manin-style
classification the slope filtration theorem
for $\sigma$-modules over arbitrary analytic rings.

In Chapter~\ref{sec:complements}, we make some related observations.
In particular, we explain how the slope filtration theorem, together
with Tsuzuki's theorem on unit-root
$\sigma$-modules with connection, implies Crew's conjecture.
We also explain the relevance of the terms ``generic'' and ``special''
to the discussion of Chapter~\ref{sec:generic special}.

\subsection*{Acknowledgments}

Thanks to Olivier Brinon for providing detailed
feedback on \cite[Chapter~3]{me-local},
to Francesco Baldassarri and Frans Oort for some suggestions
on the exposition of this paper, to Mark Kisin for suggesting
the statement of Theorem~\ref{T:isoclinic lattice},
to Bruno Chiarellotto for pointing out an error in a prior
version of Proposition~\ref{P:filtrations}, and to Eike Lau for
pointing out some additional errors.
Thanks also to Baldassarri and to Pierre Berthelot for organizing
two workshops on $F$-isocrystals and rigid cohomology in 
December 2004 and June 2005, which helped spark the aforementioned discussions.
The author was supported by NSF grant number DMS-0400727.

\section{The basic rings}
\label{sec:basic rings}

In this chapter, we recall and generalize the ring-theoretic setup of 
\cite[Chapter~3]{me-local}.

\begin{convention}
Throughout this chapter, fix a prime number $p$
and a power $q = p^a$ of $p$.
Let $K$ be a field of
characteristic $p$, equipped with a valuation $v_K$;
we will allow $v_K$ to be trivial unless otherwise specified.
Let $K_0$ denote a subfield of $K$ on which $v_K$ is trivial.
We will frequently do matrix calculations; in so doing, we apply
a valuation to a matrix by taking its minimum over entries,
and write $I_n$ for the $n \times n$ identity matrix over any ring.
See Conventions~\ref{conv:more notation1}
and~\ref{conv:more notation2} for some further notations.
\end{convention}

\subsection{Witt rings}

\begin{convention} \label{conv:perfect}
Throughout this section only, assume that $K$ and $K_0$ are perfect.
\end{convention}

\begin{defn}
Let $W(K)$
denote the ring of $p$-typical Witt vectors over $K$.
Then $W$ gives a covariant
functor from perfect fields of characteristic $p$ to
complete discrete valuation rings of characteristic 0,
with maximal ideal $p$ and perfect residue field; this functor is in
fact an equivalence of categories, being a quasi-inverse of the
residue field functor.
In particular, the absolute ($p$-power) Frobenius lifts uniquely to
an automorphism $\sigma_0$ of $W(K)$; write $\sigma$ for the $\log_p(q)$-th
power of $\sigma_0$.
Use a horizontal overbar to denote the reduction map from $W(K)$ to $K$.
In this notation, we have 
$\overline{u^{\sigma}} = \overline{u}^q$ for all $u \in W(K)$.
\end{defn}

We will also want to allow some ramified extensions of Witt rings.
\begin{defn}
Let $\calO$ be a finite totally ramified extension of $W(K_0)$,
equipped with an extension of $\sigma$; let $\pi$ denote a uniformizer
of $\calO$.
Write $W(K, \calO)$ for $W(K) \otimes_{W(K_0)} \calO$, and extend the notations
$\sigma, \overline{x}$ to $W(K, \calO)$ in the natural fashion.
\end{defn}

\begin{defn}
For $\overline{z} \in K$, let $[\overline{z}] \in W(K)$ denote the Teichm\"uller lift of $K$;
it can be constructed as $\lim_{n \to \infty} y_n^{p^n}$ for
any sequence $\{y_n\}_{n=0}^\infty$ with 
$\overline{y}_n = \overline{z}^{1/p^n}$. (The point is that this limit is well-defined:
if $\{y'_n\}_{n=0}^\infty$ is another such sequence,
we have $y_n^{p^n} \equiv (y'_n)^{p^n} \pmod{p^n}$.)
Then $[\overline{z}]^\sigma = [\overline{z}]^q$, and if $\overline{z}' \in K$, then $[\overline{z} \overline{z}'] = [\overline{z}][\overline{z}']$.
Note that
each $x \in W(K, \calO)$
can be written uniquely as $\sum_{i=0}^\infty [\overline{z_i}] \pi^i$
for some $\overline{z_0}, \overline{z_1}, \dots \in K$; similarly,
each $x \in W(K, \calO)[\fp]$
can be written uniquely as $\sum_{i \in \ZZ} [\overline{z_i}] 
\pi^i$ for some $z_i \in K$ with $\overline{z_i} = 0$ for $i$ sufficiently small.
\end{defn}

\begin{defn} \label{D:partial}
Recall that $K$ was assumed to be
equipped with a valuation $v_K$.
Given $n \in \ZZ$, we define the ``partial 
valuation'' $v_n$ on $W(K, \calO)[\fp]$ by
\begin{equation} \label{eq:vn}
v_n\left(\sum_i [\overline{z_i}] \pi^i\right) = \min_{i \leq n} \{ v_K(\overline{z_i})\};
\end{equation}
it satisfies the properties
\begin{align*}
v_n(x+y) &\geq \min\{v_n(x), v_n(y)\} \qquad (x,y \in W(K, \calO)[\fp], 
n \in \ZZ) \\
v_n(xy) &\geq \min_{m \in \ZZ} \{v_m(x) + v_{n-m}(y)\} \qquad (x,y \in 
W(K, \calO)[\fp], n \in \ZZ) \\
v_n(x^{\sigma}) &= q v_n(x) \qquad (x \in W(K, \calO)[\fp], n \in \ZZ) \\
v_n([\overline{z}]) &= v_K(\overline{z}) \qquad (\overline{z} \in K, n \geq 0).
\end{align*}
In each of the first two inequalities, one has equality if
the minimum is achieved exactly once.
For $r > 0$, $n \in \ZZ$, and $x \in W(K, \calO)[\fp]$, put
\[
v_{n,r}(x) = r v_n(x) + n;
\]
for $r = 0$, put $v_{n,r}(x) = n$ if $v_n(x) < \infty$ and $v_{n,r}(x) = \infty$
if $v_n(x) = \infty$.
For $r \geq 0$, let $W_r(K,\calO)$ be the subring of $W(K, \calO)$ 
consisting of those $x$
for which $v_{n,r}(x) \to \infty$ as $n \to \infty$; then
$\sigma$ sends $W_{qr}(K, \calO)$ onto $W_r(K, \calO)$.
(Note that there is no restriction when $r=0$.)
\end{defn}

\begin{lemma} \label{L:additivity comp}
Given $x,y \in W_r(K, \calO)[\fp]$ nonzero, let $i$ and $j$ be the smallest and
largest integers $n$ achieving $\min_n \{v_{n,r}(x)\}$, and let
$k$ and $l$ be the smallest and largest integers $n$ achieving 
$\min_n\{v_{n,r}(y)\}$. Then $i+k$ and $j+l$ are the smallest and largest
integers $n$ achieving
$\min_n\{v_{n,r}(xy)\}$, and this minimum equals
$\min_n\{v_{n,r}(x)\} + \min_n\{v_{n,r}(y)\}$.
\end{lemma}
\begin{proof}
We have
\[
v_{m,r}(xy) \geq \min_n \{ v_{n,r}(x) + v_{m-n,r}(y)\},
\]
with equality if the minimum on the right is achieved only once.
This means that:
\begin{itemize}
\item
for all $m$, the minimum is at least $v_{i,r}(x) + v_{k,r}(y)$;
\item
for $m = i+k$ and $m = j+l$, the value $v_{i,r}(x) + v_{k,r}(y)$
 is achieved exactly once
(respectively by $n=i$ and $n=j$);
\item
for $m < i+k$ or $m > j+l$, the value $v_{i,r}(x) + v_{k,r}(y)$
is never achieved.
\end{itemize}
This implies the desired results.
\end{proof}

\begin{defn} \label{D:wr}
Define the map $w_r: W_r(K, \calO)[\fp] \to \RR \cup \{\infty\}$ by
\begin{equation} \label{eq:wr}
w_r(x) = \min_n \{ v_{n,r}(x) \};
\end{equation}
also write $w$ for $w_0$.
By Lemma~\ref{L:additivity comp},
$w_r$ is a valuation on $W_r(K, \calO)[\fp]$;
moreover, $w_r(x) = w_{r/q}(x^{\sigma})$.
Put 
\[
W_{\con}(K, \calO) = \cup_{r>0} W_r(K, \calO);
\]
note that $W_{\con}(K, \calO)$ is a discrete valuation ring with residue 
field $K$
and maximal ideal generated by $\pi$, but is not complete if $v_K$ is 
nontrivial.
\end{defn}
\begin{remark}
Note that $u$ is a unit in $W_r(K, \calO)$ if and only if
$v_{n,r}(u) > v_{0,r}(u)$ for $n > 0$. We will generalize this observation
later in Lemma~\ref{L:additivity}.
\end{remark}

\begin{remark} \label{R:normalize}
Note that $w$ is a $p$-adic valuation on $W(K,\calO)$
normalized so that $w(\pi) = 1$.
This indicates two discrepancies from choices made in \cite{me-local}.
First, we have normalized $w(\pi) = 1$ instead of $w(p) = 1$ for internal
convenience; the normalization will not affect any of the final results.
Second, we use $w$ for the $p$-adic valuation instead of $v_p$ (or
simply $v$) because we are using $v$'s for
valuations in the ``horizontal'' direction, such as the valuation on $K$, 
and the partial valuations of Definition~\ref{D:partial}.
By contrast, decorated $w$'s denote ``nonhorizontal'' valuations,
as in Definition~\ref{D:wr}.
\end{remark}

\begin{lemma} \label{L:henselian}
The (noncomplete) discrete valuation ring $W_{\con}(K, \calO)$ is henselian.
\end{lemma}
\begin{proof}
It suffices to verify that if
$P(x)$ is a polynomial over $W_{\con}(K, \calO)$ and 
$y \in W_{\con}(K, \calO)$ satisfies $P(y) \equiv 0 \pmod{\pi}$
and $P'(y) \not\equiv 0 \pmod{\pi}$, then there exists 
$z \in W_{\con}(K, \calO)$ with $z \equiv y \pmod{\pi}$ and
$P(z) = 0$. To see this, pick $r>0$ such that
$w_r(P(y)/P'(y)) > 0$; then the usual Newton iteration gives a series
converging under $w$ to a root $z$ of $P$ in $W(K,\calO)$
with $z \equiv y \pmod{\pi}$. However,
the iteration also converges under $w_r$, so we must have
$z \in W_r(K, \calO)$. (Compare \cite[Lemma~3.9]{me-local}.)
\end{proof}

\subsection{Cohen rings}
\label{subsec:cohen}

Remember that Convention~\ref{conv:perfect} is no longer in force,
i.e., $K_0$ and $K$ no longer need be perfect.

\begin{defn}
Let $C_K$ denote a \emph{Cohen ring} of $K$, i.e., a complete
discrete valuation ring with maximal ideal $pC_K$ and 
residue field $K$. Such a ring necessarily exists and is unique up
to noncanonical isomorphism 
\cite[Proposition~0.19.8.5]{ega4}. Moreover, any map $K \to K'$
can be lifted, again noncanonically, to a map $C_K \to C_{K'}$.
\end{defn}

\begin{convention} \label{conv:more notation1}
For the remainder of the chapter, assume chosen and fixed a map
(necessarily injective) $C_{K_0} \to C_K$.
Let $\calO$ be a finite totally ramified extension of $C_{K_0}$,
and let $\pi$ denote a uniformizer of $\calO$.
Write $\Gamma^K$ for $C_K \otimes_{C_{K_0}} \calO$; we write $\Gamma$ for
short if $K$ is to be understood, as it will usually be in this chapter.
\end{convention}

\begin{defn}
By a \emph{Frobenius lift} on $\Gamma$, we mean any endomorphism
$\sigma: \Gamma \to \Gamma$ lifting the absolute $q$-power Frobenius on $K$.
Given $\sigma$, we may form the completion of the direct limit
\begin{equation} \label{eq:comp dir lim}
\Gamma^K \stackrel{\sigma}{\to} \Gamma^K \stackrel{\sigma}{\to} \cdots;
\end{equation}
for $K = K_0$, this ring is a finite totally ramified extension of
$\Gamma^{K_0} = \calO$, which we denote by $\calO^{\perf}$. In general, 
if $\sigma$
is a Frobenius lift on $\Gamma^K$ which maps $\calO$ into
itself, we may
identify the completed direct limit of
\eqref{eq:comp dir lim} with $W(K^{\perf}, \calO^{\perf})$;
we may thus use the induced embedding $\Gamma^K \hookrightarrow W(K^{\perf},
\calO^{\perf})$ to define $v_n, v_{n,r}, w_r, w$ on $\Gamma$.
\end{defn}
\begin{remark}
In \cite{me-local}, a Frobenius lift is assumed to be a power of
a $p$-power Frobenius lift, but all calculations therein work in this
slightly less restrictive setting.
\end{remark}

\begin{convention} \label{conv:more notation2}
For the remainder of the chapter, assume chosen and fixed a Frobenius lift
$\sigma$ on $\Gamma$ which carries $\calO$ into itself.
\end{convention}

\begin{defn}
Define the \emph{levelwise topology} on $\Gamma$ by declaring that
a sequence $\{x_l\}_{l=0}^\infty$ converges to zero if and only if for each $n$,
$v_n(x_l) \to \infty$ as $l \to \infty$. This topology is coarser
than the usual $\pi$-adic topology.
\end{defn}

\begin{defn}
For $L/K$ finite separable, we may view
$\Gamma^L$ as a finite unramified extension of $\Gamma^K$,
and $\sigma$ extends uniquely to $\Gamma^L$; if $L/K$ is Galois,
then $\Gal(L/K)$ acts on $\Gamma^L$ fixing $\Gamma^K$.
More generally, we say $L/K$ is \emph{pseudo-finite separable} if $L
= M^{1/q^n}$ for some $M/K$ finite separable and some nonnegative
integer $n$; in this case, we define $\Gamma^L$ to be a copy of
$\Gamma^M$ viewed as a $\Gamma^M$-algebra via $\sigma^n$.
In particular, we have a unique extension of $v_K$ to $L$,
under which $L$ is complete,
and we have a distinguished extension of $\sigma$ to $\Gamma^L$
(but only because we built the choice of $\sigma$ into the definition
of $\Gamma^L$).
\end{defn}

\begin{remark}
One can establish a rather strong functoriality for the formation
of the $\Gamma^L$, as in \cite[Section~2.2]{me-local}.
One of the simplifications introduced here is to avoid having to elaborate
upon this.
\end{remark}

\begin{defn}
For $r>0$,
put $\Gamma_r = \Gamma \cap W_r(K^{\perf}, \calO)$.
We say $\Gamma$ \emph{has enough $r$-units} if every nonzero element of
$K$ can be lifted to a unit of $\Gamma_r$. We say $\Gamma$ \emph{has enough
units} if $\Gamma$ has enough $r$-units for some $r>0$. 
\end{defn}

\begin{remark} \label{R:enough units}
\begin{enumerate}
\item[(a)] If $K$ 
is perfect, then $\Gamma$ has enough $r$-units for any $r>0$, because
a nonzero Teichm\"uller element is a unit in every $\Gamma_r$.
\item[(b)] If $\Gamma^K$ has enough $r$-units, then $\Gamma^{K^{1/q}}$ has
enough $qr$-units, and vice versa.
\end{enumerate}
\end{remark}

\begin{lemma} \label{L:has enough units}
Suppose that $\Gamma^K$ has enough units, and let $L$ be a pseudo-finite
separable extension of $K$. Then $\Gamma^L$ has enough units.
\end{lemma}
\begin{proof}
It is enough to check the case when $L$ is actually finite
separable. Put $d = [L:K]$.
Apply the primitive element theorem to produce $\overline{x} \in L$
which generates $L$ over $K$, and apply Lemma~\ref{L:henselian} to produce
$x \in \Gamma^L_{\con}$ lifting $\overline{x}$.

Recall that any two Banach norms on a finite 
dimensional vector space over a complete
normed field are equivalent \cite[Proposition~4.13]{schneider}. In particular,
if we let $v_L$ denote the unique extension of $v_K$ to $L$, then
there exists a constant $a>0$ such that whenever
$\overline{y} \in L$ and $\overline{c_0}, \dots, \overline{c_{d-1}} \in K$
satisfy $\overline{y} = \sum_{i=0}^{d-1} \overline{c_i} \overline{x}^i$, 
we have $v_L(\overline{y}) \leq \min_i\{v_K(\overline{c_i} \overline{x}^i)\} + a$.

Pick $r>0$ such that $\Gamma^K$ has enough $r$-units and
$x$ is a unit in $\Gamma^L_r$, and choose
$s > 0$ such that $1 - s/r > sa$.
Given $\overline{y} \in L$, lift each $\overline{c_i}$ to either zero
or a unit in $\Gamma_r$, and set $y = \sum_{i=0}^{d-1} c_i x^i$. Then
for all $n \geq 0$,
\begin{align*}
v_{n,r}(y) &\geq \min_i \{v_{n,r}(c_i x^i)\} \\
&\geq \min_i \{ r v_K(\overline{c_i}) + r i v_L(\overline{x})\} \\
&\geq r v_L(\overline{y}) - ra.
\end{align*}
In particular, $v_{n,s}(y) > v_{0,s}(y)$ for $n > 0$, so $y$ is a unit in
$\Gamma^L_s$. Since $s$ does not depend on $y$, we conclude that
$\Gamma^L$ has enough $s$-units, as desired.
\end{proof}

\begin{defn}
Suppose that $\Gamma$ has enough units. Define
$\Gcon = \cup_{r>0} \Gamma_r = \Gamma \cap W_{\con}(K, \calO)$;
then $\Gcon$ is again a discrete valuation ring with maximal ideal generated
by $\pi$. Although $\Gcon$ is not complete, it is
henselian thanks to Lemma~\ref{L:henselian}.
For $L/K$ pseudo-finite separable, we may view $\GLcon$ as an
extension of $\GKcon$, which is finite unramified if $L/K$ is finite separable.
\end{defn}

\begin{remark}
Remember that $v_K$ is allowed to be trivial, in which case the distinction
between $\Gamma$ and $\Gcon$ collapses.
\end{remark}

\begin{prop}
Let $L$ be a finite separable extension of $K$. Then for any $x \in 
\GLcon$ such that
$\overline{x}$ generates $L$ over $K$, we have
$\GLcon \cong \GKcon[x]/(P(x))$, where $P(x)$ denotes the minimal 
polynomial of $x$.
\end{prop}
\begin{proof}
Straightforward.
\end{proof}

\begin{convention} \label{conv:no completion}
For $L$ the \emph{completed} perfect closure or 
algebraic closure of $K$, we replace
the superscript $L$ by ``$\perf$'' or ``$\alg$'', respectively,
writing $\Gperf$ or $\Galg$ for $\Gamma^L$ and so 
forth. (Recall that these are obtained by embedding
$\Gamma^K$ into $W(K^{\perf}, \calO)$ via $\sigma$,
and then embedding the latter into $W(L, \calO)$ via Witt
vector functoriality.)
Beware that this convention disagrees with a convention of
\cite{me-local}, in which
$\Galg = W(K^{\alg}, \calO)$, without the completion; we will
comment further on this discrepancy in Remark~\ref{R:no completion}.
\end{convention}

The next assertions are essentially \cite[Proposition~8.1]{dejong},
only cast a bit more generally; compare also
\cite[Proposition~4.1]{me-full}.

\begin{defn}
By a \emph{valuation $p$-basis} of $K$, we mean
a subset $S \subset K$ such that the
set $U$ of monomials in $S$ of degree $<p$ in each factor (and degree
0 in almost all factors) is a valuation basis of $K$ over $K^p$. That is,
each $\overline{x} \in K$ has a unique expression of the form $\sum_{\overline{u} \in U}
\overline{c_u} \overline{u}$, with each $\overline{c_u} \in K^p$ and almost all zero, and one has
\[
v_K(\overline{x}) = \min_{\overline{u} \in U} \{ v_K(\overline{c_u} \overline{u})\}.
\]
\end{defn}

\begin{example}
For example, 
$K = k((t))$ admits a valuation $p$-basis consisting of $t$ plus
a $p$-basis of $k$ over $k^p$.
In a similar vein, if $[v(K^*):v((K^p)^*)] = [K:K^p] < \infty$, then
one can choose a valuation $p$-basis for $K$ by selecting
elements of $K^*$ whose images under $v$ generate $v(K^*)/v((K^p)^*)$.
(See also the criterion of \cite[Chapter~9]{kuhlmann}.)
\end{example}

\begin{lemma} \label{L:projection}
Suppose that $\Gamma$ has enough units and that $K$ admits a 
valuation $p$-basis $S$.
Then there exists a $\Gamma$-linear map $f: \Gperf \to \Gamma$
sectioning the inclusion $\Gamma \to \Gperf$, which maps
$\Gperfcon$ to $\Gcon$.
\end{lemma}
\begin{proof}
Choose $r>0$ such that $\Gamma$ has enough $r$-units,
and, for each $\overline{s} \in S$, choose a unit $s$ of $\Gamma_r$ lifting $\overline{s}$.
Put $U_0 = \{1\}$. For $n$ a positive integer, let $U_n$ be the set
of products
\[
\prod_{\overline{s} \in S} (s^{e_s})^{\sigma^{-n}}
\]
in which each $e_s \in \{0, \dots, q^n-1\}$, all but finitely many
$e_s$ are zero (so the product makes sense), and the $e_s$ are not
all divisible by $q$. 
Put $V_n = U_0 \cup \cdots \cup U_n$;
then the reductions of $V_n$ form a basis of $K^{q^{-n}}$ over $K$.
We can thus write each element of $\Gamma^{\sigma^{-n}}$ uniquely as
a sum $\sum_{u \in V_n} x_u u$, with each $x_u \in \Gamma$ and for any
integer $m>0$, only finitely many of the $x_u$ nonzero modulo $\pi^m$.
Define the map $f_n: \Gamma^{\sigma^{-n}} \to \Gamma$ sending $x 
= \sum_{u \in V_n} x_u u$ to $x_1$.

Note that each element of each $U_n$ is a unit in $(\Gamma^{\sigma^{-n}})_r$.
Since $S$ is a valuation $p$-basis, it follows (by induction on $m$)
that if we write
$x = \sum_{u \in V_n} x_u u$, then
\[
\min_{j \leq m} \{v_{j,r}(x)\}
= \min_{u \in V_n} \min_{j \leq m} \{v_{j,r}(x_u u)\}.
\]
Hence for any $r' \in (0,r]$, $f_n$ sends $(\Gamma^{\sigma^{-n}})_{r'}$
to $\Gamma_{r'}$. That means in particular that the $f_n$ fit together to
give a function $f$ that extends by continuity to all of $\Gperf$, sections
the map $\Gamma \to \Gperf$, and carries $\Gperfcon$ to $\Gcon$.
\end{proof}

\begin{remark}
It is not clear to us whether it should be
possible to loosen the restriction that $K$ must have a valuation
$p$-basis, e.g., by imitating the proof strategy of
Lemma~\ref{L:has enough units}.
\end{remark}

\begin{prop} \label{P:injective}
Suppose that $\Gamma$ has enough units and that $K$ admits a 
valuation $p$-basis.
Let $\mu: \Gamma \otimes_{\Gcon} \Galgcon \to \Galg$ denote the
multiplication map, so that $\mu(x \otimes y) = xy$.
\begin{enumerate}
\item[(a)]
If $x_1, \dots, x_n \in \Gamma$ are linearly independent
over $\Gcon$, and $\mu(\sum_{i=1}^n x_i \otimes y_i) = 0$, then
$y_i = 0$ for $i=1, \dots, n$.
\item[(b)] 
If $x_1, \dots, x_n \in \Gamma$ are linearly independent
over $\Gcon$, and $\mu(\sum_{i=1}^n x_i \otimes y_i) \in \Gamma$, then
$y_i \in \Gcon$ for $i=1, \dots, n$.
\item[(c)]
The map $\mu$ is injective.
\end{enumerate}
\end{prop}
\begin{proof}
\begin{enumerate}
\item[(a)]
Suppose the contrary; choose a counterexample with $n$ minimal.
We may assume without loss of generality that $w(y_1) = \min_i \{w(y_i)\}$;
we may then divide through by $y_1$ to reduce to the case $y_1 = 1$,
where we will work hereafter.

Any $g \in \Gal(K^{\alg}/K^{\perf})$ extends uniquely to an
automorphism of $\Galg$ over $\Gperf$, and to an automorphism
of $\Galgcon$ over $\Gperfcon$. Then
\[
0 = \sum_{i=1}^n x_i y_i = \sum_{i=1}^n x_i y_i^g
= \sum_{i=2}^n x_i (y_i^g - y_i);
\]
by the minimality of $n$, we have $y_i^g = y_i$ for $i=2, \dots, n$.
Since this is true for any $g$, we have $y_i \in \Gperfcon$ for each $i$.

Let $f$ be the map from Lemma~\ref{L:projection};
then 
\[ 
0 = \sum x_i y_i = f\left(\sum x_i y_i\right) = \sum x_i f(y_i)
= \sum x_i (y_i - f(y_i)),
\]
so again $y_i = f(y_i)$ for $i=2, \dots, n$.
Hence $x_1 = -\sum_{i=2}^n x_i y_i$, contradicting the
linear independence of the $x_i$ over $\Gcon$.
\item[(b)]
For $g$ as in (a), we have $0 = \sum x_i (y_i^g - y_i)$; by
(a), we have $y_i^g = y_i$ for all $i$ and $g$, so $y_i \in \Gperfcon$.
Now $0 = \sum x_i (y_i - f(y_i))$, so $y_i = f(y_i) \in \Gcon$.
\item[(c)]
Suppose on the contrary that $\sum_{i=1}^n x_i \otimes y_i \neq 0$
but $\sum_{i=1}^n x_i y_i = 0$; choose such a counterexample with
$n$ minimal. By (a), the $x_i$ must be linearly dependent over $\Gcon$; without
loss of generality, suppose we can write
$x_1 = \sum_{i=2}^n c_i x_i$ with $c_i \in \Gcon$.
Then $\sum_{i=1}^n x_i \otimes y_i = 
\sum_{i=2}^n x_i \otimes (y_i + c_i)$ is a counterexample with only
$n-1$ terms, contradicting the minimality of $n$.
\end{enumerate}
\end{proof}

\subsection{Relation to the Robba ring}
\label{subsec:robba}

We now recall how the constructions in the previous section relate to
the usual Robba ring.

\begin{convention}
Throughout this section, assume that 
$K = k((t))$ and $K_0 = k$; we may then describe
$\Gamma$ as the ring of formal Laurent series $\sum_{i \in \ZZ} c_i u^i$
with each $c_i \in \calO$, and $w(c_i) \to \infty$ as $i \to -\infty$.
Suppose further that the Frobenius lift is given by
$\sum c_i u^i \mapsto \sum c_i^\sigma (u^\sigma)^i$, where
$u^\sigma = \sum a_i u^i$ with $\liminf_{i \to -\infty} w(a_i)/(-i) > 0$.
\end{convention}

\begin{defn}
Define the na\"\i ve partial valuations $v_n^{\naive}$ on $\Gamma$ by
the formula
\[
v^{\naive}_n \left(\sum c_i u^i \right) 
= \min\{i: w(c_i) \leq n\}.
\]
These functions satisfy some identities analogous to those in
Definition~\ref{D:partial}:
\begin{align*}
v^{\naive}_n(x+y) &\geq \min\{v^{\naive}_n(x), v^{\naive}_n(y)\} 
\qquad (x,y \in \Gamma[\fp], 
n \in \ZZ) \\
v^{\naive}_n(xy) &\geq \min_{m \leq n} \{v^{\naive}_m(x) 
+ v^{\naive}_{n-m}(y)\} \qquad (x,y \in 
\Gamma[\fp], n \in \ZZ).
\end{align*}
Again, equality holds in each case if the minimum on the right side
is achieved exactly once.
Put 
\[
v^{\naive}_{n,r}(x) = r v^{\naive}_n(x) + n.
\]
For $r>0$, let $\Gamma_r^{\naive}$ be the set of $x \in \Gamma$ such that
$v^{\naive}_{n,r}(x) \to \infty$ as $n \to \infty$.
Define the map $w^{\naive}_r$ on $\Gamma_r^{\naive}$ by
\[
w^{\naive}_r(x) = \min_n \{ v^{\naive}_{n,r}(x)\};
\]
then $w^{\naive}_r$ is a valuation on $\Gamma_r^{\naive}[\fp]$ by the same
argument as in Lemma~\ref{L:additivity comp}.
Put 
\[
\Gcon^{\naive} = \cup_{r>0} \Gamma_r^{\naive}.
\]
By the hypothesis on the Frobenius lift,
we can choose $r>0$ such that $u^\sigma/u^q$ is a unit in 
$\Gamma_r^{\naive}$.
\end{defn}

\begin{lemma} \label{L:naive iterate}
For $r>0$ such that $u^\sigma/u^q$ is a unit in $\Gamma_r^{\naive}$,
and $s \in (0,qr]$, we have
\begin{equation} \label{eq:compare naive1}
\min_{j \leq n} \{v_{j,s}^{\naive}(x) \} = 
\min_{j \leq n} \{v_{j,s/q}^{\naive}(x^\sigma) \}
\end{equation}
for each $n \geq 0$ 
and each $x \in \Gamma$.
\end{lemma}
\begin{proof}
The hypothesis ensures that \eqref{eq:compare naive1} holds for
$x = u^i$ for any $i \in \ZZ$ and any $n$.
For general $x$,
write $x = \sum_i c_i u^i$; then on one hand,
\begin{align*}
\min_{j \leq n} \{ v_{j,s/q}^{\naive}(x^\sigma)
\} &\geq \min_{i \in \ZZ} \{ \min_{j \leq n}
\{ v_{j,s/q}^{\naive}(c_i^\sigma (u^\sigma)^i) \} \}\\
&= \min_{i \in \ZZ} \{\min_{j \leq n}
\{ v_{j,s}^{\naive}(c_i u^i) \} \}\\
&= \min_{j \leq n} \{v_{j,s}^{\naive}(x)\}.
\end{align*}
On the other hand,
if we take the smallest $j$ achieving the minimum on the left side of
\eqref{eq:compare naive1}, then the minimum of
$v_{j,s}^{\naive}(c_i u^i)$ is achieved by a \emph{unique} integer $i$.
Hence the one inequality in the previous sequence is actually an equality.
\end{proof}

\begin{lemma}
For $r>0$ such that $u^\sigma/u^q$ is a unit in 
$\Gamma_r^{\naive}$, and $s \in (0, qr]$, we have
\begin{equation} \label{eq:compare naive}
\min_{j \leq n} \{v_{j,s}(x)\} = \min_{j \leq n} \{ v_{j,s}^{\naive}(x)
\}
\end{equation}
for each $n \geq 0$ 
and each $x \in \Gamma$. In particular, $\Gamma^{\naive}_s = 
\Gamma_s$, and $w_s(x) = w_s^{\naive}(x)$ for all $x \in \Gamma_s$.
\end{lemma}
\begin{proof}
Write $x = \sum_{i=0}^\infty [\overline{x_i}] \pi^i$ with each
$\overline{x_i} \in K^{\perf}$. Choose an integer $l$
such that $\overline{x_i}^{q^l} \in K$ for $i=0,\dots,n$,
and write $\overline{x_i}^{q^l} = \sum_{h \in \ZZ} \overline{c_{hi}} t^h$
with $\overline{c_{hi}}
 \in k$. Choose $c_{hi} \in \calO$ lifting $\overline{c_{hi}}$,
with $c_{hi} = 0$ whenever $\overline{c_{hi}} = 0$, and put
$y_i = \sum_h c_{hi} u^h$.

Pick an integer $m>n$, and define
\[
x' = \sum_{i=0}^n y_i^{q^m} (\pi^i)^{\sigma^{l+m}};
\]
then $w(x' - x^{\sigma^{l+m}}) > n$. Hence for $j \leq n$,
$v_j(x') = v_j(x^{\sigma^{l+m}}) = q^{l+m} v_j(x)$
and $v_j^{\naive}(x') = v_j^{\naive}(x^{\sigma^{l+m}})$.

From the way we chose the $y_i$, we have
\[
v_j^{\naive}(y_i^{q^m} (\pi^i)^{\sigma^{l+m}}) = q^{l+m} v_0(\overline{x_i}) \qquad
(j \geq i).
\]
It follows that $v_j^{\naive}(x') = q^{l+m} v_j(x)$ for $j \leq n$;
that is, we have $v_j^{\naive}(x^{\sigma^{l+m}}) = q^{l+m} v_j(x)$ for 
$j \leq n$. In particular, we have
\[
\min_{j \leq n} \{v_{j,s}(x)\} = \min_{j \leq n} \{ 
v_{j,s/q^{l+m}}^{\naive}(x^{\sigma^{l+m}})\}.
\]
By Lemma~\ref{L:naive iterate}, this yields the desired result.
(Compare \cite[Lemmas~3.6 and~3.7]{me-local}.)
\end{proof}
\begin{cor}
For $r>0$ such that $u^\sigma/u^q$ is a unit in $\Gamma_r^{\naive}$,
$\Gamma$ has enough $qr$-units, and $\Gcon = \Gamma^{\naive}_{\con}$.
\end{cor}

\begin{remark} \label{R:annulus}
The ring $\Gamma^{\naive}_r[\fp]$ is the ring of bounded rigid analytic
functions on the annulus $|\pi|^r \leq |u| < 1$, and the valuation $w^{\naive}_s$ is the supremum norm on
the circle $|u| = |\pi|^s$. This geometric interpretation motivates the
subsequent constructions, and so is worth keeping in mind; indeed, much of
the treatment of analytic rings in the rest of this chapter is modeled on
the treatment of rings of functions on annuli given by
Lazard \cite{lazard}, and our results generalize some of the results
in \cite{lazard} (given Remark~\ref{R:standard} below).
\end{remark}

\begin{remark} \label{R:standard}
In the context of this section, the ring $\Gancon$ is what is
usually called the \emph{Robba ring} over $K$.
The point of view of \cite{me-local}, maintained here, is that
the Robba ring should always be viewed as coming with the 
``equipment'' of a Frobenius lift $\sigma$; this seems to be
the most convenient angle from which to approach
$\sigma$-modules. However, when discussing a
statement about $\Gancon$ that depends only on its underlying topological
ring (e.g., the B\'ezout property, as in Theorem~\ref{T:bezout}),
one is free to use any Frobenius, and so it is sometimes convenient
to use a ``standard'' Frobenius lift, under which $u^{\sigma} = u^q$
and $v^{\naive}_n = v_n$ for all $n$.
In general, however, one cannot get away with only standard Frobenius lifts
because the property of standardness is not preserved by passing to
$\Gamma^L_{\an,\con}$ for $L$ a finite separable extension of $k((t))$.
\end{remark}

\begin{remark} \label{R:no spherical0}
It would be desirable to be able to have it possible for
$\Gamma^{\naive}_r$ to be the ring
of rigid analytic functions on an annulus over a $p$-adic field whose
valuation is not discrete (e.g., the completed algebraic closure
$\CC_p$ of $\QQ_p$), since the results of Lazard we are analogizing
hold in that context. However, this seems rather difficult to accommodate
in the formalism developed above; 
for instance, the $v_n$ cannot be described in terms
of Teichm\"uller elements, so an axiomatic characterization is probably
needed.
There are additional roadblocks later in the story; we will flag some of
these as we go along.
\end{remark}

\begin{remark}
One can carry out an analogous comparison between na\"\i ve and true
partial valuations when $K$ is the completion of $k(x_1, \dots, x_n)$
for a ``monomial'' valuation, in which $v(x_1), \dots, v(x_n)$ are
linearly independent over $\QQ$; this gives additional examples in which
the hypothesis ``$\Gamma$ has enough units'' can be checked, and hence
additional examples in which the framework of this paper applies.
See \cite{me-fake} for details.
\end{remark}

\subsection{Analytic rings}

We now proceed roughly as in \cite[Section~3.3]{me-local};
however, we will postpone certain ``reality checks'' on the definitions
until the next section.
\begin{convention}
Throughout this section, and for the rest of the chapter,
assume that the field $K$ is 
complete with respect to the valuation $v_K$, and that $\Gamma^K$
has enough $r_0$-units for some fixed $r_0>0$. Note that the assumption that
$K$ is complete ensures that $\Gamma_r$ is complete under $w_r$ for any 
$r \in [0,r_0)$. 
\end{convention}

\begin{defn}
Let $I$ be a subinterval of $[0,r_0)$ bounded away from $r_0$, i.e.,
$I \subseteq [0,r]$ for some $r<r_0$.
Let $\Gamma_{I}$ be the Fr\'echet completion of $\Gamma_{r_0}[\fp]$
for the valuations $w_s$ for $s \in I$; note that the functions
$v_n, v_{n,s}, w_s$ extend to $\Gamma_{I}$ extend by continuity,
and that $\sigma$ extends to a map 
$\sigma: \Gamma_{I} \to \Gamma_{q^{-1} I}$.
For $I \subseteq J$ subintervals of $[0,r_0)$ bounded away from 0,
we have a natural map $\Gamma_J \to \Gamma_I$; this map is injective
with dense image.
For $I = [0,s]$, note that $\Gamma_I = \Gamma_s[\fp]$.
For $I = (0,s]$, we write $\Gamma_{\an,s}$ for $\Gamma_I$.
\end{defn}

\begin{remark}
In the context of Section~\ref{subsec:robba}, $\Gamma_I$ is the ring of
rigid analytic functions on the subspace of the open unit disc
defined by the condition $\log_{|\pi|} |u| \in I$;
compare Remark~\ref{R:annulus}.
\end{remark}

\begin{defn}
For $I$ a subinterval of $[0,r_0)$ bounded away from $r_0$,
and for $x \in \Gamma_{I}$ nonzero, define the \emph{Newton polygon}
of $x$ to be the lower convex hull of the set of points
$(v_n(x), n)$, minus any segment the negative of whose slope is
not in $I$. Define the \emph{slopes} of
$x$ to be the negations of the slopes of the Newton polygon of $x$.
Define the \emph{multiplicity} of $s \in (0,r]$ as a slope of $x$ to
be the difference in vertical coordinates between the endpoints of 
the segment of the Newton polygon of $x$ of slope $-s$, or 0 if no
such segment exists. 
If $x$ has only finitely many slopes, define 
the \emph{total multiplicity} of $x$ to be the sum of the multiplicities
of all slopes of $x$.
If $x$ has only one slope, we say $x$ is \emph{pure} of that slope.
\end{defn}

\begin{remark}
The analogous definition of total multiplicity for $\Gamma^{\naive}_r$ counts the total number
of zeroes (with multiplicities) that a function has in the annulus
$|\pi|^r \leq |u| < 1$.
\end{remark}

\begin{remark}
Note that the multiplicity of any given slope is always finite.
More generally, for any closed subinterval $I = [r',r]$ of $[0,r_0)$,
the total multiplicity of any $x \in \Gamma_I$ is finite.
Explicitly, the total multiplicity equals $i-j$, where
$i$ is the largest $n$ achieving $\min_n \{v_{n,r}(x)\}$ and
$j$ is the smallest $n$ achieving
$\min_n\{v_{n,r'}(x)\}$.
In particular, if $x \in \Gamma_{\an,r}$, the slopes of $x$ form a
sequence decreasing to zero.
\end{remark}

\begin{lemma} \label{L:additivity}
For $x,y \in \Gamma_{I}$ nonzero, the multiplicity of each
$s \in I$ as a slope of $xy$ is the sum of the multiplicities
of $s$ as a slope of $x$ and of $y$. In particular, $\Gamma_I$
is an integral domain.
\end{lemma}
\begin{proof}
For $x,y \in \Gamma_r[\fp]$, this follows at once from
Lemma~\ref{L:additivity comp}. In the general case, note that the
conclusion of Lemma~\ref{L:additivity comp} still holds, by
approximating $x$ and $y$ suitably well by elements of $\Gamma_r[\fp]$.
\end{proof}

\begin{defn}
Let $\GKancon$ be the union of the $\GK_{\an,r}$ over all $r \in (0,r_0)$;
this ring is an integral domain by Lemma~\ref{L:additivity}.
Remember that we are allowing $v_K$ to be trivial, in which case
$\Gancon = \Gcon[\fp] = \Gamma[\fp]$.
\end{defn}

\begin{example}
In the context of Section~\ref{subsec:robba},
the ring
$\Gancon$ consists of formal Laurent series $\sum_{n \in \ZZ} c_n u^n$
with each $c_n \in \calO[\fp]$, $\liminf_{n \to -\infty} w(c_n)/(-n) > 0$,
and $\liminf_{n \to \infty} w(c_n)/n \geq 0$. The latter is none other than
the Robba ring over $\calO[\fp]$.
\end{example}

We make a few observations about finite extensions of $\Gancon$.

\begin{prop} \label{P:tensor}
Let $L$ be a finite separable extension of $K$. Then the multiplication map
\[
\mu: \GKancon \otimes_{\GKcon} \GLcon \to \GLancon
\]
is an isomorphism. More precisely,
for any $x \in \GLcon$ such that $\overline{x}$ generates $L$ over $K$,
we have $\GLancon \cong \GKancon[x]/(P(x))$.
\end{prop}
\begin{proof}
For $s>0$ sufficiently small, we have $\Gamma^L_s \cong \Gamma^K_s[x]/(P(x))$
by Lemma~\ref{L:henselian}, from which the claim follows.
\end{proof}
\begin{cor} \label{C:invariant}
Let $L$ be a finite Galois extension of $K$. Then the fixed subring of
$\Frac \GLancon$ under the action of $G = \Gal(L/K)$ is equal to $\Frac 
\GKancon$.
\end{cor}
\begin{proof}
By Proposition~\ref{P:tensor}, the fixed subring of $\GLancon$ under the 
action of $G$ is equal to $\GKancon$.
Given $x/y \in \Frac \GLancon$ fixed under $G$, 
put $x' = \prod_{g \in G} x^g$; since $x'$ is $G$-invariant,
we have $x' \in \GKancon$.
Put $y' = x'y/x \in \GLancon$; then $x'/y' =x/y$, and both $x'$ and 
$x'/y'$ are $G$-invariant, so $y'$ is as well. 
Thus $x/y \in \Frac \GKancon$, as desired.
\end{proof}

\begin{lemma} \label{L:dense}
Let $I$ be a subinterval of $(0,r_0)$ bounded away from $r_0$.
Then the union $\cup \Gamma^L_I$,
taken over all pseudo-finite separable extensions $L$ of $K$,
is dense in $\Galg_I$.
\end{lemma}
\begin{proof}
Let $M$ be the algebraic closure (not completed) of $K$.
Then $\cup \Gamma^L$ is clearly dense in $\Gamma^M$ for the $p$-adic 
topology. By Remark~\ref{R:enough units} and
Lemma~\ref{L:has enough units}, the set of pseudo-finite
separable extensions $L$ such that $\Gamma^L$ has enough $r_0$-units
is cofinal. Hence the set $U$ of $x \in \cup \Gamma^L_{r_0}$ with $w_{r_0}(x) \geq 0$
is dense in the set $V$ of $x \in \Gamma^M_{r_0}$ with $w_{r_0}(x) \geq 0$ for the
$p$-adic topology. On these sets, the topology induced on $U$ or $V$
by any one $w_s$ with $s \in (0,r_0)$ is coarser than the $p$-adic topology.
Thus $U$ is also dense in $V$ for the Fr\'echet topology induced by the
$w_s$ for $s \in I$. It follows that $\cup \Gamma^L_I$ is dense in
$\Gamma^M_I$; however, the condition that $0 \notin I$ ensures that
$\Gamma^M_I = \Galg_I$, so we have the desired result.
\end{proof}

\begin{remark} \label{R:no completion}
Recall that in \cite{me-local} (contrary to our present
Convention~\ref{conv:no completion}),
the residue field of $\Gamma^{\alg}$ is the algebraic
closure of $K$, rather than the completion thereof. However,
the definition of $\Galgancon$ comes out the same, and our
convention here makes a few statements a bit easier to make.
For instance, in the notation of \cite{me-local}, an element $x$ of 
$\Galgancon$ can satisfy $v_n(x) = \infty$ for all $n<0$
without belonging to $\Galgcon$.
(Thanks to Francesco Baldassarri for suggesting this change.)
\end{remark}

\subsection{Reality checks}

Before proceeding further, we must make some tedious but necessary
``reality checks'' concerning the analytic rings. 
This is most easily done for $K$ perfect,
where elements of $\Gamma_I$ have canonical decompositions
(related to the
``strong semiunit decompositions'' of \cite[Proposition~3.14]{me-local}.)
\begin{defn} \label{D:teich}
For $K$ perfect, define the functions $f_n: \Gamma[\fp] \to K$ for $n \in \ZZ$
by the formula $x = \sum_{n \in \ZZ} [f_n(x)]\pi^n$, where the brackets
again denote Teichm\"uller lifts. Then
\[
v_n(x) = \min_{m \leq n} \{ v_K(f_m(x)) \} \leq v_K(f_n(x)),
\] 
which implies that $f_n$ extends uniquely to
a continuous function $f_n: \Gamma_I \to K$ for any subinterval $I \subseteq
[0, \infty)$, and that the sum
$\sum_{n \in \ZZ} [f_n(x)] \pi^n$ converges to $x$ in $\Gamma_I$.
We call this sum the \emph{Teichm\"uller presentation} of $x$.
Let $x_+, x_-, x_0$ be the sums of
$[f_n(x)] \pi^n$ over those $n$ for which $v_K(f_n(x))$ is positive, negative,
or zero; we call the presentation $x = x_+ + x_- + x_0$ the
\emph{plus-minus-zero presentation} of $x$.
\end{defn}

From the existence of Teich\"muller presentations, it is obvious that
for instance, if $x \in \Gamma_{\an,r}$ satisfies $v_n(x) = \infty$ for all
$n<0$, then $x \in \Gamma_r$. In order to make such statements evident
in case $K$ is not perfect, we need an approximation of the same technique.
\begin{defn}
Define a \emph{semiunit} to be an element of $\Gamma_{r_0}$ which is either
zero or a unit.
For $I \subseteq [0,r]$ and $x \in \Gamma_I$, a \emph{semiunit
presentation} of $x$ (over $\Gamma_I$)
is a convergent sum $x = \sum_{i \in \ZZ} u_i \pi^i$,
in which each $u_i$ is a semiunit.
\end{defn}

\begin{lemma} \label{L:sum}
Suppose that $u_0, u_1, \dots$ are semiunits.
\begin{enumerate}
\item[(a)]
For each $i \in \ZZ$ and $r \in (0,r_0)$,
\[
w_r(u_i \pi^i) \geq \min_{n \leq i} \left\{v_{n,r}\left(\sum_{j=0}^i 
u_i \pi^i\right)\right\}.
\]
\item[(b)]
Suppose that $\sum_{i=0}^\infty u_i \pi^i$ converges $\pi$-adically
to some $x$ such that for some $r \in (0,r_0)$,
$v_{n,r}(x) \to \infty$ as $n \to \infty$.
Then $w_r(u_i \pi^i) \to \infty$ as $i \to \infty$, so that
$\sum_i u_i \pi^i$ is a semiunit presentation of $x$ over $\Gamma_r$.
\end{enumerate}
\end{lemma}
\begin{proof}
\begin{enumerate}
\item[(a)]
The inequality is evident for $i=0$; we prove the general claim by induction
on $i$. If $w_r(u_i \pi^i) \geq w_r(u_j \pi^j)$ for some $j<i$, then
the induction hypothesis yields the claim. Otherwise,
$w_r(u_i \pi^i) < w_r(\sum_{j<i} u_j \pi^j)$, so
$v_{n,r}(\sum_{j=0}^i u_j \pi^j) 
= v_{n,r}(u_i \pi^i)$, again yielding the claim.
\item[(b)]
Choose $r' \in (r,r_0)$; we can then apply (a) to deduce that
\begin{align*}
w_{r'}(u_i \pi^i) &\geq \min_{n \leq i} \{v_{n,r'}(x)\} \\
&= \min_{n \leq i} \{(r'/r) v_{n,r}(x) + (1-r'/r)n\}.
\end{align*}
It follows that
\begin{align*}
w_r(u_i \pi^i) &\geq
\min_{n \leq i} \{v_{n,r}(x) + (r/r'-1)n\} + (1-r/r')i \\
&= \min_{n \leq i} \{v_{n,r}(x) + (1-r/r')(i-n)\}.
\end{align*}
Since $v_{n,r}(x) \to \infty$ as $n \to \infty$, the right
side tends to $\infty$ as $n \to \infty$.
\end{enumerate}
\end{proof}

\begin{lemma} \label{L:semiunit2}
Given a subinterval $I$ of $[0,r_0)$ bounded away from $r_0$,
and $r \in I$,
suppose that $x \in \Gamma_{[r,r]}$ has the property that for any $s \in I$,
$v_{n,s}(x) \to \infty$ as $n \to \pm \infty$.
Suppose also that $\sum_i u_i \pi^i$ is a semiunit presentation of
$x$ over $\Gamma_{[r,r]}$. Then
$\sum_i u_i \pi^i$ converges in $\Gamma_I$; in particular,
$x \in \Gamma_I$.
\end{lemma}
\begin{proof}
By applying Lemma~\ref{L:sum}(a) to $\sum_{i=-N}^N u_i \pi^i$
and using continuity,
we deduce that $w_{r}(u_i \pi^i) \geq \min_{n \leq i} \{v_{n,r}(x)\}$.
For $s \in I$ with $s \geq r$, we have $w_s(u_i \pi^i) \geq
(s/r) w_r(u_i \pi^i) + (s/r-1)(-i)$, so $w_s(u_i \pi^i) \to \infty$
as $i \to -\infty$. On the other hand, for $s \in I$ with $s < r$,
we have
\begin{align*}
w_s(u_i \pi^i) &= (s/r) w_r(u_i \pi^i) + (1-s/r)i \\
&\geq (s/r) \min_{n \leq i} \{v_{n,r}(x)\} + (1-s/r)i \\
&= (s/r) \min_{n \leq i} \{(r/s) v_{n,s}(x) + (1-r/s)n\} + (1-s/r)i \\
&= \min_{n \leq i} \{v_{n,s}(x) + (s/r - 1)(n-i) \} \\
&\geq \min_{n \leq i} \{v_{n,s}(x)\};
\end{align*}
by the hypothesis that $v_{n,s}(x) \to \infty$ as $n \to \pm \infty$,
we have $w_s(u_i \pi^i) \to \infty$ as $i \to -\infty$ also in this case.

We conclude that $\sum_{i<0} u_i \pi^i$ converges in $\Gamma_I$;
put $y = x - \sum_{i<0} u_i \pi^i$.
Then $\sum_{i=0}^\infty u_i \pi^i$ converges to $y$ under $w_r$,
hence also $\pi$-adically. By
Lemma~\ref{L:sum}(b), $\sum_{i=0}^\infty u_i \pi^i$ converges
in $\Gamma_r$, so we have $x \in \Gamma_I$, as desired.
\end{proof}

One then has the following variant of \cite[Proposition~3.14]{me-local}.
\begin{prop} \label{P:semiunit pres}
For $I$ a subinterval of $[0,r_0)$ bounded away from $r_0$,
every $x \in \Gamma_{I}$ admits a semiunit presentation.
\end{prop}
\begin{proof}
We first verify that for $r \in (0,r_0)$, every element of
$\Gamma_r$ admits a semiunit presentation. Given $x \in \Gamma_r$,
we can construct a sum 
$\sum_{i} u_i \pi^i$ converging $\pi$-adically to $x$,
in which each $u_i$ is a semiunit.
By Lemma~\ref{L:sum}(b), this sum actually converges under $w_s$
for each $s \in [0,r]$, hence yields a semiunit presentation.

We now proceed to the general case;
by Lemma~\ref{L:semiunit2}, it is enough to treat the case
$I = [r,r]$.
Choose a sum $\sum_{i=0}^\infty x_i$ converging to $x$
in $\Gamma_{[r,r]}$, with each $x_i\in \Gamma_r[\fp]$.
We define elements $y_{il} \in \Gamma_r[\fp]$ for $i \in \ZZ$
and $l \geq 0$, such that for each $l$, there are only finitely many
$i$ with $y_{il} \neq 0$, as follows.
By the vanishing condition on the $y_{il}$,
$x_0 + \cdots + x_l - \sum_{j<l} \sum_i y_{ij} \pi^i$ belongs to
$\Gamma_r[\fp]$ and so admits a semiunit presentation
$\sum_{i} u_i \pi^i$ by virtue of the previous paragraph.
For each $i$ with $w_r(u_i \pi^i) < w_r(x_{l+1})$ (of which there are only
finitely many), put $y_{il} = u_i$;
for all other $i$, put $y_{il} = 0$.
Then
\[
w_r \left( x_0 + \cdots + x_l - \sum_{j \leq l} \sum_i y_{ij} \pi^i
\right) \geq w_r(x_{l+1}).
\]
In particular, the doubly infinite sum $\sum_{i,l} y_{il} \pi^i$
converges to $x$ under $w_r$.
If we set $z_i = \sum_l y_{il}$, then the sum $\sum_{i} z_i \pi^i$
converges to $x$ under $w_r$.

Note that whenever $y_{il} \neq 0$, $w_r(x_l) \leq w_r(y_{il} \pi^i)$
by Lemma~\ref{L:sum}, whereas
$w_r(y_{il} \pi^i) < w_r(x_{l+1})$ by construction.
Thus for any fixed $i$, the values of 
$w_r(y_{il} \pi^i)$, taken over all $l$ such that $y_{il} \neq 0$,
form a strictly increasing sequence. Since each such $y_{il}$ is
a unit in $\Gamma_{r_0}$, we have $w_{r_0}(y_{il} \pi^i) = (r_0/r) w_r(y_{il}
\pi^i) + (1-r_0/r)i$; hence the values of $w_{r_0}(y_{il} \pi^i)$ also form
an increasing sequence. Consequently, the sum $\sum_l y_{il}$
converges in $\Gamma_{r_0}$ (not just under $w_r$) and its limit $z_i$
is a semiunit. Thus $\sum_i z_i \pi^i$ is a semiunit presentation
of $x$ over $\Gamma_{[r,r]}$, as desired.
\end{proof}

\begin{cor} \label{C:subring}
For $r \in (0,r_0)$ and
$x \in \Gamma_{[r,r]}$, we have $x \in \Gamma_r$ if and only if 
$v_n(x) = \infty$ for all $n<0$.
\end{cor}
\begin{proof}
If $x \in \Gamma_r$, then $v_n(x) = \infty$ for all $n<0$.
Conversely, suppose that $v_n(x) = \infty$ for all $n<0$.
Apply Proposition~\ref{P:semiunit pres}
to produce a semiunit presentation $x = \sum_i u_i \pi^i$.
Suppose there exists
$j<0$ such that $u_j \neq 0$; pick such a $j$ minimizing
$w_r(u_j \pi^j)$. Then $v_{j,n}(x) = w_r(u_j \pi^j) \neq \infty$,
contrary to assumption. Hence $u_j = 0$ for $j < 0$, and
so $x = \sum_{i=0}^\infty u_i \pi^i \in \Gamma_r$.
\end{proof}

\begin{cor} \label{C:subring2}
Let $I \subseteq J$ be subintervals of $[0,r_0)$ bounded away
from $r_0$. Suppose $x \in \Gamma_{I}$ has the property that
for each $s \in J$, $v_{n,s}(x) \to \infty$ as $n \to
\pm \infty$. Then $x \in \Gamma_J$.
\end{cor}
\begin{proof}
Produce a semiunit presentation of $x$ over $\Gamma_J$
using Proposition~\ref{P:semiunit pres}, then apply
Lemma~\ref{L:semiunit2}.
\end{proof}
The numerical criterion provided by Corollary~\ref{C:subring2}
in turn implies a number of results that are evident in the case
of $K$ perfect (using Teichm\"uller presentations).
\begin{cor} \label{C:intersect}
For $K \subseteq K'$ an extension of complete fields such that
$\Gamma^K$ and $\Gamma^{K'}$ have enough $r_0$-units, and
$I \subseteq J \subseteq [0,r_0)$ bounded away from $r_0$,
we have 
\[
\Gamma^K_I \cap \Gamma^{K'}_J = \Gamma^K_{J}.
\]
\end{cor}
\begin{cor} \label{C:splitting}
Let $I = [a,b]$ and $J = [c,d]$ be 
subintervals of $[0,r_0)$ bounded away from $r_0$
with $a \leq c \leq b \leq d$. Then the intersection of
$\Gamma_I$ and $\Gamma_J$ within $\Gamma_{I \cap J}$ is equal to
$\Gamma_{I \cup J}$. Moreover, 
any $x \in \Gamma_{I \cap J}$ with $w_s(x) >0$ for $s \in I \cap J$
can be written as $y + z$ with $y \in \Gamma_I$, $z \in \Gamma_J$, and
\begin{align*}
w_s(y) &\geq (s/c) w_c(x) \qquad (s \in [a,c]) \\
w_s(z) &\geq (s/b) w_b(x) \qquad (s \in [b,d]) \\
\min\{w_s(y), w_s(z)\} &\geq w_s(x) \qquad (s \in [c,b]).
\end{align*}
\end{cor}
\begin{proof}
The first assertion follows from Corollary~\ref{C:subring2}.
For the second assertion, apply Proposition~\ref{P:semiunit pres}
to obtain a semiunit presentation $x = \sum u_i \pi^i$.
Put $y = \sum_{i\leq 0} u_i \pi^i$ and $z = \sum_{i > 0} u_i \pi^i$;
these satisfy the claimed inequalities.
\end{proof}

\begin{remark}
The notion of a semiunit presentation is similar to that of a
``semiunit decomposition'' as in \cite{me-local}, but somewhat
less intricate.
In any case, we will have only limited direct use for semiunit 
presentations;
we will mostly exploit them indirectly, via their role in proving
Lemma~\ref{L:approx} below.
\end{remark}

\begin{lemma} \label{L:approx}
Let $I$ be a closed subinterval of $[0,r]$ for some $r \in (0,r_0)$,
and suppose $x \in \Gamma_{I}$. Then there
exists $y \in \Gamma_r$ such that
\[
w_s(x-y) \geq \min_{n<0} \{ v_{n,s}(x)\} \qquad (s \in I).
\]
\end{lemma}
\begin{proof}
Apply Proposition~\ref{P:semiunit pres} to produce a semiunit presentation
$\sum_i u_i \pi^i$ of $x$. 
Then we can choose $m > 0$ such that
$w_s(u_i \pi^i) > \min_{n<0} \{ v_{n,s}(x)\}$ for $s \in I$
and $i > m$. Put $y = \sum_{i=0}^m u_i \pi^i$; then the desired
inequality follows from Lemma~\ref{L:sum}(a).
\end{proof}

\begin{cor} \label{C:unit}
A nonzero element $x$ of $\Gamma_{I}$ is a unit in $\Gamma_{I}$
if and only if it has no slopes; if $I = (0,r]$, this happens if and only
if $x$ is a unit in $\Gamma_r[\fp]$.
\end{cor}
\begin{proof}
If $x$ is a unit in $\Gamma_I$, it has no slopes by
Lemma~\ref{L:additivity}. Conversely, suppose that $x$ has no slopes;
then there exists a single $m$ which minimizes $v_{m,s}(x)$ for all $s \in I$.
Without loss of generality we may assume that $m=0$;
we may then apply Lemma~\ref{L:approx} to produce $y \in \Gamma_r$
such that $w_s(x-y) \geq \min_{n<0} \{v_{n,s}(x)\}$ for all $s \in I$.
Since $\Gamma$ has enough $r$-units, we can choose a unit $z \in \Gamma_r$
such that $w(y-z) > 0$; then $w_s(1 - xz^{-1}) > 0$ for all $s \in I$.
Hence the series $\sum_{i=0}^\infty (1 - xz^{-1})^i$ converges
in $\Gamma_I$, and its limit $u$ satisfies $uxz^{-1} = 1$. This
proves that $x$ is a unit.

In case $I = (0,r]$, $x$ has no slopes if and only if there is a unique
$m$ which minimizes $v_{m,s}(x)$ for all $s \in (0,r]$; this is 
only possible if $v_n(x) = \infty$ for $n<m$. By
Corollary~\ref{C:subring}, this implies $x \in \Gamma_r[\fp]$; by the
same argument, $x^{-1} \in \Gamma_r[\fp]$.
\end{proof}

\subsection{Principality}

In Remark~\ref{R:annulus}, the annulus of which $\Gamma^{\naive}_r$ is the
rigid of rigid analytic functions is affinoid (in the sense of Berkovich
in case the endpoints are not rational) and one-dimensional, and so
$\Gamma^{\naive}_r$ is a principal ideal domain. This can be established
more generally.

Before proceeding further, we mention a useful ``positioning lemma'',
which is analogous to but not identical with \cite[Lemma~3.24]{me-local}.
\begin{lemma} \label{L:positioning}
For $r \in (0, r_0)$ and $x \in \Gamma_{[r,r]}$ nonzero, 
there exists a unit $u \in \Gamma_{r_0}$ and an integer $i$ such that,
if we write $y = u \pi^i x$, then:
\begin{enumerate}
\item[(a)]
$w_r(y) = 0$;
\item[(b)]
$v_0(y-1) > 0$;
\item[(c)]
$v_{n,r}(y) > 0$ for $n < 0$.
\end{enumerate}
\end{lemma}
\begin{proof}
Define  $i$ to be the largest integer minimizing
$v_{-i,r}(x)$.
Apply Lemma~\ref{L:approx} to find $z \in \Gamma_r$ such that 
$w_r(\pi^i x - z) \geq \min_{n<0} \{v_{n,r}(\pi^i x) \}$.
Since $\Gamma$ has enough $r_0$-units, we can choose a unit
$u$ of $\Gamma_{r_0}$ such that $u^{-1} \equiv z \pmod{\pi}$;
then $u$ and $i$ have the desired properties.
\end{proof}

\begin{defn}
For $x \in \Gamma_r$ nonzero, define the \emph{height} of $x$
as the largest $n$ such that $w_r(x) = v_{n,r}(x)$; it can also be 
described as the $p$-adic valuation of $x$ plus the total multiplicity of $x$.
By convention, we say $0$ has height $-\infty$.
\end{defn}

\begin{lemma}[Division algorithm] \label{L:divalg}
For $r \in (0,r_0)$ and
$x, y \in \Gamma_r$ with $x$ nonzero, there exists $z \in \Gamma_r$
such that $y-z$ is divisible by $x$, and $z$
has height less than that of $x$. Moreover, we can ensure that
$w_r(z) \geq w_r(y)$.
\end{lemma}
\begin{proof}
Let $m$ be the height of $x$.
Apply Proposition~\ref{P:semiunit pres} to choose a semiunit
presentation $\sum_{i=0}^\infty u_i \pi^i$ of $x$,
and put $x' = x - \sum_{i=0}^{m-1} u_i \pi^i$; then $x' \pi^{-m}$ is a unit in
$\Gamma_r$, and by
Lemma~\ref{L:sum},
\[
w_r(x - x') \geq w_r(x) + (1-r/r_0).
\]

Define a sequence $\{y_l\}_{l=0}^\infty$ as follows. Put $y_0 = y$.
Given $y_l$ with $y_l - y$ divisible by $x$ and $w_r(y_l) 
\geq w_r(y)$, if $y_l$ has height less 
than $m$, we may take $z = y_l$ and be done 
with the proof of the lemma.
So we may assume that $y_l$ has height at least $m$, which
means that $\min_n \{w_{n,r}(y_l)\}$ is achieved by at least one
$n \geq m$.

Pick $y'_l\in \Gamma_{r_0}$ with $w_r(y'_l-y_l) \geq w_r(y_l) + (1-r/r_0)$,
and apply Lemma~\ref{L:approx} to $y'_l \pi^{-m}$ to produce
$z_l \in \Gamma_{r_0}$ such that $w_{r_0}(z_l - y'_l \pi^{-m}) \geq
\min_{n<0} \{ v_{n,r_0}(y'_l \pi^{-m})\}$.
Put
\begin{align*}
y_{l+1} &= y_l - z_l (\pi^m/x') x \\
&= (y_l - y'_l) + (y'_l - z_l \pi^m) + z_l \pi^m (1 - x'/x).
\end{align*}
By construction, we have $w_r(y_l-y'_l) \geq w_r(y_l) + (1-r/r_0)$
and $w_r(z_l \pi^m (1-x'/x)) \geq w_r(y_l) + (1-r/r_0)$.
Moreover, for $n \geq m$, we have
\begin{align*}
v_{n,r}(y'_l - z_l \pi^m) &= (r/r_0) v_{n,r_0}(y'_l - z_l \pi^m) + (1-r/r_0)n \\
&\geq (r/r_0) w_{r_0}(y'_l - z_l \pi^m) + (1-r/r_0)m \\
&\geq (r/r_0) \min_{j < m} \{v_{j,r_0}(y'_l)\} + (1-r/r_0)m \\
&= \min_{j < m} \{ v_{j,r}(y'_l) + (1-r/r_0)(m-j)\}  \\
&\geq \min_{j < m} \{ v_{j,r}(y'_l)\} + (1-r/r_0) \\
&\geq w_r(y'_l) + (1-r/r_0) = w_r(y_l) + (1-r/r_0).
\end{align*}
It follows that for $n \geq m$, we have
$v_{n,r}(y_{l+1}) \geq w_r(y_l) + (1-r/r_0)$.
We may assume that $y_{l+1}$ also has height at least $m$,
in which case $w_r(y_{l+1}) \geq w_r(y_l) + (1-r/r_0)$.
Hence (unless the process stops at some finite $l$, in which case we
already know that we win) the $y_l$ converge to zero under $w_r$, and
\[
y = x \sum_{l=0}^\infty (y_l - y_{l+1})/x \in \Gamma_r
\]
is divisible by $x$, so we may take $z=0$.
\end{proof}

\begin{remark}
Note how we used the fact that 
$\Gamma^K$ has enough $r_0$-units, not just enough $r$-units.
Also, note that the discreteness of the valuation on $K$ was essential
to ensuring that the sequence $\{y_l\}$ converges to zero.
\end{remark}

This division algorithm has the usual consequence.
\begin{prop} \label{P:bounded PID1}
For $r \in (0,r_0)$, $\Gamma_r$ is a principal ideal domain.
\end{prop}
\begin{proof}
Let $J$ be a nonzero ideal of $\Gamma_r$, and pick $x \in J$ of minimal height. 
Then for any $y \in J$, apply Lemma~\ref{L:divalg} to produce
$z$ of height less than $x$ with $y-z$ divisible by $x$.
Then $z \in J$, so we must have $z = 0$ by the minimality in the choice of $x$.
In other words, every $y \in J$ is divisible by $x$, as claimed.
\end{proof}
\begin{remark}
Here is one of the roadblocks mentioned in Remark~\ref{R:no spherical0}:
if $\calO$ is not discretely valued, then it is not even a PID itself,
so the analogue of $\Gamma_r$ cannot be one either.
\end{remark}

To extend Proposition~\ref{P:bounded PID1}
to more $\Gamma_I$, we use the following factorization
lemma (compare \cite[Lemma~3.25]{me-local}).
We will refine this lemma a bit later; see Lemma~\ref{L:factor unit2}.
\begin{lemma} \label{L:factor unit}
For $I = [r',r] \subseteq [0,r_0)$ and $x \in \Gamma_I$, there exists
a unit $u$ of $\Gamma_I$ such that $ux \in \Gamma_r$,
and all of the slopes of $ux$ in $[0,r]$ belong to $I$.
\end{lemma}
\begin{proof}
By applying Lemma~\ref{L:positioning}, we may reduce to the case where 
$w_{r'}(x) = 0$, $v_0(x-1) > 0$, and $v_{n,r'}(x) > 0$ for $n < 0$; then
for $n < 0$, we must have $v_n(x) > 0$ and so $v_{n,s}(x) > 0$ for all 
$s \in I$. Put 
\[
c = \min_{s \in I} \{ \min_{n \leq 0} \{v_{n,s}(x-1)\} \} > 0.
\]

Define the sequence $u_0, u_1, \dots$ of units of $\Gamma_I$ as follows.
First set $u_0 = 1$. Given $u_l$ such that
$\min_{n \leq 0} \{v_{n,s}(u_l x - 1)\} \geq c$ for all $s \in I$, 
apply Lemma~\ref{L:approx} to
produce $y_l \in \Gamma_r$ such that 
$w_s(y_l - u_l x) \geq \min_{n<0} \{v_{n,s}(u_l x)\}$ for all $s \in I$.
We may thus take $u_{l+1} = u_l(1 - y_l + u_l x)$, because 
$w_s(y_l - u_l x) \geq c$;
moreover, for $n < 0$,
\begin{align*}
v_{n,r'}(u_{l+1}x)
&= v_{n,r'}(y_l - u_{l+1} x) \\
&= v_{n,r'}((y_l - u_l x)(1-u_lx)) \\
&\geq \min_{m} \{v_{m,r'}(y_l - u_l x) + v_{n-m,r'}(1-u_lx)\}.
\end{align*}
This last minimum is at least $\min_{n<0} \{v_{n,r'}(u_l x)\}$.
Moreover, if it is ever less than $\min_{n<0} \{v_{n,r'}(u_l x)\}
+ c$, then the smallest value of $n$
achieving $\min_n \{v_{n,r'}(u_{l+1}x)\}$ is strictly greater than
the smallest value of $n$ achieving
$\min_n \{v_{n,r'}(u_l x)\}$ (since in that case,
terms in the last minimum above with $m \leq 0$ cannot affect the
minimum of the left side over all $n<0$).

In other words, for every $l$, there exists $l' > l$ such that
\[
\min_{n<0} \{v_{n,r'}(u_{l'}x) \}
\geq \min_{n<0} \{v_{n,r'}(u_lx)\} + c.
\]
Hence in case $s = r'$, we have $\min_{n < 0} \{ v_{n,s}(u_l x)\} \to \infty$ 
as $l \to \infty$; consequently, the same also holds for $s \in I$.
It follows that the sequence $\{u_l\}$ converges to a limit $u \in \Gamma_I$,
and that $v_{n,s}(ux) = \infty$ for $n<0$, so that $ux \in \Gamma_r$
by Corollary~\ref{C:subring}.
Moreover, by construction,
$\min_n \{v_{n,r'}(ux)\}$ is achieved by $n=0$, so all of the slopes
of $ux$ are at least $r'$.
\end{proof}

\begin{prop} \label{P:bounded PID}
Let $I$ be a closed subinterval of $[0,r_0)$. Then 
$\Gamma_I$ is a principal ideal domain.
\end{prop}
\begin{proof}
Put $I = [r',r]$, and let $J$ be a nonzero ideal of $\Gamma_I$.
By Lemma~\ref{L:factor unit}, each element $x$
of $J$ can be written (nonuniquely) as a unit $u$
of $\Gamma_I$ times an element $y$ of $\Gamma_r$. Let $J'$ be the ideal
of $\Gamma_r$ generated by all such $y$; by Proposition~\ref{P:bounded
PID}, $J'$ is principal, generated by some $z$. Since $J' \subseteq
J \cap \Gamma_r$, we have $z \in J$; on the other hand, each $x \in J$
has the form $uy$ with $u \in \Gamma_I$ and $y \in \Gamma_r$, and
$y$ is a multiple of $z$ in $\Gamma_r$, so $x$ is a multiple of $z$
in $\Gamma_I$. Hence $z$ generates $J$, as desired.
\end{proof}

\begin{remark}
Proposition~\ref{P:bounded PID} generalizes
Lazard's \cite[Corollaire de Proposition~4]{lazard}.
\end{remark}

\subsection{Matrix approximations and factorizations}

We need  a matrix approximation lemma similar to \cite[Lemma~6.2]{me-local};
it is in some sense a matricial analogue of Lemma~\ref{L:positioning}.
\begin{lemma} \label{L:approximate matrix}
Let $I$ be a closed subinterval of $[0,r]$ for some $r \in (0,r_0)$,
and let $M$ be an invertible $n \times n$ matrix over $\Gamma_I$.
Then there exists an invertible $n \times n$ matrix $U$ over $\Gamma_r[\fp]$
such that $w_s(MU - I_n) > 0$ for $s \in I$. Moreover, if
$w_s(\det(M) - 1) > 0$, we can ensure that $\det(U) = 1$.
\end{lemma}
\begin{proof}
By applying Lemma~\ref{L:positioning} to $\det(M)$ (and then multiplying
a single row of $U$ by the resulting unit), we can ensure that
$w_s(\det(M) - 1) > 0$ for $s \in I$. With this extra hypothesis,
we proceed by induction on $n$.

Let $C_i$ denote the cofactor of $M_{ni}$ in $M$, so that
$\det(M) = \sum_{i=1}^n C_i M_{ni}$, and in fact
$C_i = (M^{-1})_{in} \det(M)$.
Put $\alpha_i = \det(M)^{-1} M_{ni}$, so that
$\sum_{i=1}^n \alpha_i C_i = 1$.

Choose $\beta_1, \dots, \beta_{n-1}, \beta'_n \in \Gamma_r[\fp]$
such that for $s \in I$ and $i=1,\dots,n-1$,
\[
w_s(\beta_i - \alpha_i) > -w_s(C_i)
\qquad
w_s(\beta'_n - \alpha_n) > -w_s(C_n).
\]
Note that for $c \in \calO$ with $w(c)$ sufficiently large, 
$\beta_n = \beta'_n + c$ satisfies
$w_s(\beta_n - \alpha_n) > -w_s(C_n)$ for $s \in I$.
Moreover, by Proposition~\ref{P:bounded PID},
we can find $\gamma \in \Gamma_r[\fp]$ generating the ideal 
generated by $\beta_1, \dots, \beta_{n-1}$;
then the $\beta'_n+c$ are pairwise coprime for different $c \in \calO$,
so only finitely many of them can have a nontrivial common factor
with $\gamma$. In particular, for $w(c)$ sufficiently large,
$\beta_1, \dots, \beta_n$ generate the unit ideal in $\Gamma_r[\fp]$.

With $\beta_1, \dots, \beta_n$ so chosen, we can 
choose a matrix $A$ over $\Gamma_r[\fp]$ of determinant 1 such that
$A_{ni} = \beta_i$ for $i=1, \dots, n$ (because $\Gamma_r[\fp]$ is a
PID, again by Proposition~\ref{P:bounded PID}). 
Put $M' = MA^{-1}$, and let $C'_n$ be the cofactor of
$M'_{nn}$ in $M'$. Then
\begin{align*}
C'_n &= (AM^{-1})_{nn} \det(M) \\
&= \sum_{i=1}^n A_{ni} (M^{-1})_{in} \det(M) = \sum_{i=1}^n \beta_i C_i,
\end{align*}
so that
\[
C'_n = 1 + \sum_{i=1}^n (\beta_i - \alpha_i) C_i
\]
and so $w_s(C'_n - 1) > 0$ for $s \in I$.
In particular, $C'_n$ is a unit in $\Gamma_I$.

Apply the induction hypothesis to the upper left $(n-1)\times (n-1)$
submatrix of $M'$, and extend the resulting $(n-1) \times (n-1)$
matrix $V$ to
an $n \times n$ matrix by setting $V_{ni} = V_{in} = 0$
for $i=1, \dots, n-1$ and $V_{nn} = 1$. Then we have
$\det(M'V) = \det(M)$, so $w_s(\det(M'V) - 1) > 0$ for $s \in I$, and 
\[
w_s((M'V - I_n)_{ij}) > 0 \qquad (i=1,\dots,n-1;\, j=1,\dots, n-1;\, s \in I).
\]

We now perform an ``approximate Gaussian elimination'' over $\Gamma_I$ 
to transform
$M'V$ into a new matrix $N$ with $w_s(N-I_n) > 0$ for $s \in I$.
First, define a sequence of matrices $\{X^{(h)}\}_{h=0}^\infty$ by
$X^{(0)} = M'V$ and 
\[
X_{ij}^{(h+1)} = \begin{cases} X_{ij}^{(h)} & i<n \\
X_{nj}^{(h)} - \sum_{m=1}^{n-1} X_{nm}^{(h)} X_{mj}^{(h)} & i=n;
                 \end{cases}
\]
note that $X^{(h+1)}$ is obtained from $X^{(h)}$ by subtracting
$X_{nm}^{(h)}$ times the $m$-th row from the $n$-th row for 
$m=1,\dots,n-1$ in succession. At each step, for each $s \in I$,
$\min_{1 \leq j \leq n-1} \{w_s(X_{nj}^{(h)})\}$
increases by at least
$\min_{1 \leq i,j \leq n-1}
\{w_s((M'V - I_n)_{ij})\}$; the latter is bounded away from zero over
all $s \in I$, because
$I$ is closed and $w_s(x)$ is a continuous
function of $s$.
Thus for $h$ sufficiently large, we have
\[
w_s(X_{nj}^{(h)}) > \max\{0, \max_{1 \leq i \leq n-1} \{-w_s(X_{in}^{(h)})
\}\} \qquad (s \in I;\, j = 1, \dots, n-1).
\]
Pick such an $h$ and set $X = X^{(h)}$; note that $\det(X) = \det(M'V)$,
so $w_s(\det(X) - 1) > 0$ for $s \in I$.
For $s \in I$,
\begin{align*}
w_s((X-I_n)_{ij}) &> 0 \qquad (i=1,\dots, n;\, j =1, \dots, n-1) \\
w_s(X_{in} X_{nj}) &> 0 \qquad (i=1,\dots, n-1;\, j =1, \dots, n-1) 
\end{align*}
and hence also $w_s(X_{nn} - 1) > 0$.

Next, we perform ``approximate backsubstitution''.
Define a sequence of matrices $\{W^{(h)}\}_{h=0}^\infty$ by
setting $W^{(0)} = X$ and
\[
W_{ij}^{(h+1)} = \begin{cases} W_{ij}^{(h)} - W_{in}^{(h)} W_{nj}^{(h)} &
i<n \\
W_{ij}^{(h)} & i=n; \end{cases}
\]
note that $W^{(h+1)}$ is obtained from $W^{(h)}$ by subtracting
$W_{in}^{(h)}$ times the $n$-th row from the $i$-th row for
$i=1, \dots, n-1$. At each step, for $s \in I$,
$w_s(W_{in}^{(h)})$ increases by at least
$w_s(X_{nn}-1)$; again, the latter is bounded away from zero over
all $s \in I$ because $I$ is closed and $w_s(x)$ is continuous in $s$.
Thus for $h$ sufficiently large,
\[
w_s(W^{(h)}_{in}) > 0 \qquad (s \in I;\, 1 \leq i \leq n-1).
\]
Pick such an $h$ and set $W = W_h$; then $w_s(W-I_n) > 0$ for $s \in I$.
(Note that the inequality $w_s(X_{in} X_{nj}) > 0$ for
$i=1, \dots, n-1$ and $j=1, \dots, n-1$ ensures that
the second set of row operations does not disturb the fact that
$w_s(W^{(h)}_{ij}) > 0$ for $i=1,\dots, n-1$ and $j=1,\dots, n-1$.)

To conclude, note that by construction, $(M'V)^{-1} W$ is a product
of elementary matrices over $\Gamma_I$, each consisting of the diagonal
matrix plus one off-diagonal entry. By suitably approximating the
off-diagonal entry of each
matrix in the product by an element of $\Gamma_r$, we get an invertible
matrix $Y$ over $\Gamma_r$ such that $w_s(M'VY - I_n) > 0$ for $s \in I$.
We may thus take $U = A^{-1}VY$ to obtain the desired result.
\end{proof}

We also need a factorization lemma in the manner
of \cite[Lemma~6.4]{me-local}.
\begin{lemma} \label{L:factor matrix}
Let $I = [a,b]$ and $J = [c,d]$ be subintervals
of $[0,r_0)$ bounded away from $r_0$, with $a \leq c \leq b \leq d$,
and let $M$ be an $n \times n$ matrix over $\Gamma_{I \cap J}$
with $w_s(M-I_n) > 0$ for $s \in I \cap J$. Then
there exist invertible $n \times n$ matrices $U$ over $\Gamma_I$
and $V$ over $\Gamma_J$ such that $M = UV$.
\end{lemma}
\begin{proof}
We construct sequences of matrices $U_l$ and $V_l$ over $\Gamma_I$
and $\Gamma_J$, respectively, with
\begin{align*}
w_s(U_l - I_n) &\geq (s/c) w_c(M - I_n) \qquad (s \in [a,c]) \\
w_s(V_l - I_n) &\geq (s/b) w_b(M - I_n) \qquad (s \in [b,d]) \\
\min\{w_s(U_l - I_n), w_s(V_l - I_n)\} &\geq w_s(M - I_n) \qquad
(s \in [c,b]) \\
w_s(U_l^{-1}M V_l^{-1} - I_n) &\geq 2^l w_s(M - I_n) \qquad (s \in [c,b]),
\end{align*}
as follows. Start with $U_0 = V_0 = I_n$.
Given $U_l, V_l$, put $M_l = U_l^{-1} M V_l^{-1}$.
Apply Corollary~\ref{C:splitting} to split $M_l - I_n = Y_l + Z_l$
with $Y_l$ defined over $\Gamma_I$, $Z_l$ defined over $\Gamma_J$,
and
\begin{align*}
w_s(Y_l) &\geq (s/c)w_c(M_l - I_n) \geq (s/c)w_c(M - I_n) \qquad (s \in [a,c]) \\
w_s(Z_l) &\geq (s/b)w_b(M_l - I_n) \geq (s/b)w_b(M - I_n) \qquad (s \in [b,d]) \\
\min\{w_s(Y_l), w_s(Z_l)\} &\geq w_s(M_l - I_n) \geq 2^l w_s(M - I_n) \qquad (s \in [c,b]).
\end{align*}
Put $U_{l+1} = U_l(I+Y_l)$ and $V_{l+1} = (I + Z_l)V_l$;
then one calculates that $w_s(M_{l+1} - I_n) 
\geq 2^{l+1} w_s(M - I_n)$ for $s \in [c,b]$.

We deduce that the sequences $\{U_l\}$ and $\{V_l\}$ each converge
under $w_s$ for $s \in [c,b]$, and the limits $U$ and $V$ satisfy
$\min\{w_s(U - I_n), w_s(V - I_n)\} \geq w_s(M - I_n)$ for $s \in [c,b]$,
and $M = UV$. However, the subset $x \in \Gamma_I$
on which 
\[
w_s(x) \geq \begin{cases} (s/c)w_c(M-I_n) & s \in [a,c] \\
w_s(M-I_n) & s \in [c,b]
             \end{cases}
\]
is complete under any one $w_s$,
so $U$ has entries in $\Gamma_I$ and 
$w_s(U-I_n) \geq (s/c) w_c(M - I_n)$ for $s \in [a,c]$.
Similarly, $V$ has entries in $\Gamma_J$ and
$w_s(V-I_n) \geq (s/b) w_b(M - I_n)$ for $s \in [b,d]$.
In particular,
$U$ and $V$ are invertible over $\Gamma_I$ and $\Gamma_J$,
and $M = UV$, yielding the desired factorization.
\end{proof}

\subsection{Vector bundles}

Over an open rigid analytic annulus, one specifies a vector bundle by
specifying a vector bundle (necessarily freely generated by global sections)
on each closed subannulus and providing glueing data; if the field of 
coefficients is spherically complete, it can be shown that the result is 
again freely generated by global sections. Here we generalize
the discretely valued case of this result
to analytic rings. (For rank 1, the annulus statement can be
extracted from results of \cite{lazard}; the general case can be found in
\cite[Theorem~3.4.3]{me-monoover}. In any case, it follows from our
Theorem~\ref{T:vector bundles} below.)

\begin{defn} \label{D:vector bundle}
Let $I$ be a subinterval of $[0,r_0)$ bounded away from $r_0$,
and let $S$ be a collection of closed subintervals of $I$ closed
under finite intersections, whose union is all of $I$.
Define an \emph{$S$-vector bundle} over $\Gamma_I$ to be a 
collection consisting of one 
finite free $\Gamma_J$-module $M_J$ for each $J \in S$,
plus isomorphisms
\[
\iota_{J_1,J_2}: M_{J_1} \otimes_{\Gamma_{J_1}} \Gamma_{J_2} \cong M_{J_2}
\]
whenever $J_2 \subseteq J_1$, satisfying the compatibility condition
$\iota_{J_2,J_3} \circ \iota_{J_1,J_2} = \iota_{J_1,J_3}$
whenever $J_3 \subseteq J_2 \subseteq J_1$.
These may be viewed as forming a category
in which a morphism between the collections $\{M_J\}$ and $\{N_J\}$ consists
of a collection of morphisms $M_J \to N_J$ of $\Gamma_J$-modules which
commute with the isomorphisms $\iota_{J_1,J_2}$.
\end{defn}

This definition obeys the analogue of the usual glueing property for
coherent sheaves on an affinoid space (i.e., the theorem of
Kiehl-Tate).

\begin{lemma} \label{L:vector bundles}
Let $I$ be a subinterval of $[0,r_0)$ bounded away from $r_0$,
and let $S_1 \subseteq S_2$ be two collections of closed subintervals
of $I$ as in Definition~\ref{D:vector bundle}.
Then the natural functor from the category of $S_2$-vector bundles
over $\Gamma_I$ to $S_1$-vector bundles over $\Gamma_I$ is an equivalence.
\end{lemma}
\begin{proof}
We define a quasi-inverse functor as follows.
Given $J \in S_2$, by compactness we can choose $J_1, \dots, J_m \in S_1$
with $J \subseteq J' = J_1 \cup \cdots \cup J_m$; it is enough to consider
the case where $m=2$ and $J_1 \cap J_2 \neq \emptyset$, as we can
repeat the construction to treat the general case.

Define $M_{J'}$ to be the $\Gamma_{J'}$-submodule of
$M_{J_1} \oplus M_{J_2}$ consisting of those pairs
$(\bv_1, \bv_2)$ such that 
\[
\iota_{J_1,J_1 \cap J_{2}}(\bv_1) = \iota_{J_{2},J_1 \cap 
J_{2}}(\bv_{2}).
\]
Let $\bv_1, \dots, \bv_n$ be a basis of $M_{J_1}$ and let
$\bw_1, \dots, \bw_n$ be a basis of $M_{J_2}$. Then
there is an invertible $n \times n$ matrix $A$ over $M_{J_1 \cap J_2}$
given by $\bw_j = \sum_i A_{ij} \bv_i$. By
Lemma~\ref{L:factor matrix}, $A$ can be factored as $UV$, where
$U$ is invertible over $\Gamma_{J_1}$ and $V$ is invertible over
$\Gamma_{J_2}$. Set 
\[
\be_j = \left(\sum_i U_{ij} \bv_i, \sum_i (V^{-1})_{ij} \bw_i \right);
\]
then $\be_1, \dots, \be_n$ form a basis of $M_{J'}$, since the first
components form a basis of $M_{J_1}$, the second components form
a basis of $M_{J_2}$, and the intersection of $\Gamma_{J_1}$ and
$\Gamma_{J_2}$ within $\Gamma_{J_1 \cap J_2}$ equals $\Gamma_{J'}$
(by Corollary~\ref{C:splitting}). In particular, the natural maps
$M_{J'} \otimes_{\Gamma_{J'}} \Gamma_{J_i} 
\to M_{J_i}$ for $i=1,2$ are isomorphisms.
We may thus set $M_J = M_{J'} \otimes_{\Gamma_{J'}} \Gamma_J$.
\end{proof}
\begin{defn}
By Lemma~\ref{L:vector bundles},
the category of $S$-vector bundles over $\Gamma_I$
is canonically independent of the
choice of $S$. We thus refer to its elements simply as 
\emph{vector bundles over $\Gamma_I$}.
\end{defn}

It follows that for $I$ closed, any vector bundle over $\Gamma_I$
is represented by a free module; a key result for us is that one has
a similar result over $\Gamma_{\an,r}$.
\begin{theorem} \label{T:vector bundles}
For $r \in (0,r_0)$,
the natural functor from finite free $\Gamma_{\an,r}$-modules to
vector bundles over $\Gamma_{\an,r} = \Gamma_{(0,r]}$ is an equivalence.
\end{theorem}
\begin{proof}
To produce a quasi-inverse functor,
let $J_1 \subseteq J_2 \subseteq \cdots$ be an increasing 
sequence of closed intervals with right endpoints $r$,
whose union is $(0,r]$; for ease of notation,
write $\Gamma_i$ for $\Gamma_{J_i}$.
We can specify a vector bundle over $\Gamma_i$ by specifying a finite
free $\Gamma_i$-module $E_i$ for each $i$, plus identifications
$E_{i+1} \otimes_{\Gamma_{i+1}} \Gamma_i \cong E_i$.

Choose a basis $\bv_{1,1}, \dots, \bv_{1,n}$ of $E_1$.
Given a basis $\bv_{i,1}, \dots, \bv_{i,n}$ of $E_i$, we
choose a basis $\bv_{i+1,1}, \dots, \bv_{i+1,n}$ of $E_{i+1}$ as 
follows. Pick any basis $\be_1, \dots, \be_n$ of $E_{i+1}$, and define
an invertible $n \times n$ matrix $M_i$ over $\Gamma_i$ by writing
$\be_l = \sum_{j} (M_i)_{jl} \bv_{i,j}$. Apply Lemma~\ref{L:approximate
matrix} to produce an invertible $n \times n$ matrix $U_i$ over $\Gamma_r$
such that $w_s(M_iU_i-I_n) > 0$ for $s \in J_i$.
Apply Lemma~\ref{L:approx} to produce an $n \times n$ matrix $V_i$
over $\Gamma_r$ with $w_s(M_iU_i - I_n - V_i) 
\geq \min_{m<0} \{v_{m,s}(M_iU_i - I_n)\}$
for $s \in J_i$; then $w_r(V_i) > 0$, so $I_n + V_i$ is invertible
over $\Gamma_r$.
Put $W_i = M_iU_i(I_n + V_i)^{-1}$,
and define $\bv_{i+1,1}, \dots, \bv_{i+1,n}$
by $\bv_{i+1,l} = \sum_j (W_i)_{jl} \bv_{i,j}$; these
form another basis of $E_{i+1}$ because we changed basis over
$\Gamma_r$.

If we write $J_i = [r_i,r]$, then for any fixed $s \in (0,r]$, 
we have 
\begin{align*}
w_s(W_i - I_n) 
&= w_s((M_i U_i - I_n - V_i)(I_n + V_i)^{-1}) \\
&\geq \min_{m<0} \{ v_{m,s}(M_iU_i - I_n) \} \\
&= \min_{m<0} \{ (s/r_i) v_{m,r_i}(M_i U_i - I_n) + (1-s/r_i)m \} \\
&\geq \min_{m<0} \{ v_{m,r_i}(M_i U_i - I_n)\} + (s/r_i - 1) \\
&> (s/r_i - 1),
\end{align*}
which tends to $\infty$ as $i \to \infty$. Thus the product
$W_1 W_2 \cdots$ converges to an invertible matrix $W$ over $\Gamma_{\an,r}$,
and the basis $\be_1, \dots, \be_n$ of $E_1$ defined by
\[
\be_l = \sum_j W_{jl} \bv_{1,j}
\]
actually forms a basis of each $E_i$. Hence the original vector bundle
can be reconstructed from the free $\Gamma_{\an,r}$-module generated by
$\be_1, \dots, \be_n$; this yields the desired quasi-inverse.
\end{proof}
\begin{cor} \label{C:closed free}
For $r \in (0,r_0)$,
let $M$ be a finite free $\Gamma_{\an,r}$-module. Then every closed
submodule of $M$ is free; in particular, every closed ideal of
$\Gamma_{\an,r}$ is principal.
\end{cor}
\begin{proof}
A submodule is closed if and only if it gives rise to a sub-vector 
bundle of the vector bundle associated to $M$; thus
the claim follows from Theorem~\ref{T:vector bundles}.
\end{proof}

\begin{remark} \label{R:vector bundle}
One might expect that more generally every vector bundle over $\Gamma_I$
is represented by a finite free $\Gamma_I$-module; we did not verify this.
\end{remark}

\subsection{The B\'ezout property}

One pleasant consequence of Theorem~\ref{T:vector bundles} is the
fact that the ring $\Gamma_{\an,r}$ has the B\'ezout property,
as we verify in this section.
We start by refining the conclusion of Lemma~\ref{L:factor unit}
(again, compare \cite[Lemma~3.25]{me-local}).
\begin{lemma} \label{L:factor unit2}
For $r,s,s' \in (0,r_0)$ with $s'<s < r$, 
and $f \in \Gamma_{[s',r]}$, there exists $g \in \Gamma_r[\fp]$ with the
following properties.
\begin{enumerate}
\item[(a)]
The ideals generated by $f$ and $g$ in $\Gamma_{[s,r]}$ coincide.
\item[(b)]
The slopes of $g$ in $[0,r]$ are all contained in $[s,r]$.
\end{enumerate}
Moreover, any such $g$ also has the following property.
\begin{enumerate}
\item[(c)]
$f$ is divisible by $g$ in $\Gamma_{[s',r[}$.
\end{enumerate}
\end{lemma}
\begin{proof}
By Lemma~\ref{L:factor unit}, we can find a unit $u$ of $\Gamma_{[s',r]}$
such that $uf \in \Gamma_r[\fp]$ and the slopes of $uf$ in $[0,r]$
are all contained in $[s,r]$. We may thus take $g=uf$ to obtain at least
one $g \in \Gamma_r[\fp]$ 
satisfying (a) and (b); hereafter, we let $g$ be any element of
$\Gamma_r[\fp]$ satisfying (a) and (b). Then the multiplicity of each
element of $[s,r]$ as a slope of $g$ is equal to its
multiplicity as a slope of $f$.

Since $\Gamma_r[\fp]$ is a PID by Proposition~\ref{P:bounded
PID}, we can find an element $h \in \Gamma_r[\fp]$ generating the
ideal generated by $uf$ and $g$ in $\Gamma_r[\fp]$; in particular,
the multiplicity of each element of $[s,r]$ as a slope of $h$ is less
than or equal to its multiplicity as a slope of $g$. 
However, $h$ must also generate the ideal generated by $f$ and $g$
in $\Gamma_{[s,r]}$, which is generated already by $f$ alone;
in particular, the multiplicity of each element of $[s,r]$ as a slope
of $f$ is  equal to its multiplicity as a slope of $h$.

We conclude that each element of $[s,r]$ occurs as a slope of $f,g,h$
all with the same multiplicity. Since $g$ only has slopes in $[s,r]$,
$g/h$ must be a unit in $\Gamma_r[\fp]$;
hence $uf$ is already divisible by $g$
in $\Gamma_r[\fp]$, so $f$ is divisible by $g$ in $\Gamma_{[s',r]}$
as desired.
\end{proof}

\begin{lemma} \label{L:principal parts1}
Given $r \in (0,r_0)$ and $x \in \Gamma_r[\fp]$ with greatest slope $s_0 < r$,
choose $r' \in (s_0,r)$. Then for any $y \in \Gamma_{r}[\fp]$ and any $c>0$,
there exists $z \in \Gamma_{r}[\fp]$ with $y-z$ divisible by $x$ in
$\Gamma_{\an,r}$, such that $w_s(z) > c$ for $s \in [r',r]$.
\end{lemma}
\begin{proof}
As in the proof of Lemma~\ref{L:positioning}, we can find a unit
$u \in \Gamma_r$ and an integer $i$ such that
$\min_n \{v_{n,r}(u x \pi^i)\}$
is achieved by $n=0$ but not by any $n>0$, and that
$v_0(ux\pi^i - 1) > 0$. Since $s_0 < r$, in fact
$\min_n \{v_{n,r}(ux \pi^i)\}$ is only achieved by $n=0$,
so $w_r(ux \pi^i - 1) > 0$. Similarly,
$w_s(ux \pi^i - 1) > 0$ for $s \in [r',r]$; since $[r',r]$
is a closed interval, we can choose $d>0$ such that
$w_s(ux \pi^i - 1) \geq d$ for $s \in [r',r]$.
Now simply take
\[
z = y(1 - ux \pi^i)^N
\]
for some integer $N$ with $Nd + w_s(y) > c$ for $s \in [r',r]$.
\end{proof}

We next introduce a ``principal parts lemma'' (compare 
\cite[Lemma~3.31]{me-local}).
\begin{lemma} \label{L:principal parts2}
For $r \in (0,r_0)$, let $I_1 \subset I_2 \subset \cdots$ be an increasing
sequence of closed subintervals of $(0,r]$ with right endpoints $r$,
whose union is all of $(0,r]$, and put $\Gamma_i = \Gamma_{I_i}$.
Given $f \in \Gamma_{\an,r}$ and $g_i \in \Gamma_i$ such that
for each $i$, $g_{i+1}-g_i$ is divisible by $f$ in $\Gamma_i$,
there exists $g \in \Gamma_{\an,r}$ such that for each $i$,
$g-g_i$ is divisible by $f$ in $\Gamma_i$.
\end{lemma}
\begin{proof}
Apply Lemma~\ref{L:factor unit2} 
to produce $f_i \in \Gamma_r$ 
dividing $f$ in $\Gamma_{\an,r}$, such that $f$ and $f_i$ generate
the same ideal in $\Gamma_i$, $f_i$ only has slopes in $I_i$,
and $f/f_i$ has no slopes in
$I_i$; put $f_0 = 1$. 
By Lemma~\ref{L:factor unit2} again (with $s'$ varying), $f_i$ is divisible
by $f_{i-1}$ in $\Gamma_{\an,r}$, hence also in $\Gamma_r[\fp]$; 
put $h_i = f_i/f_{i-1} \in \Gamma_r[\fp]$ and $h_0 = 1$.

Set $x_0 = 0$. Given $x_i \in \Gamma_r[\fp]$ 
with $x_i - g_i$ divisible by $f_i$ in $\Gamma_i$,
note that  the ideal generated by $h_{i+1}$ and $f_i$ in
$\Gamma_r[\fp]$ is principal by Proposition~\ref{P:bounded PID}.
Moreover, any generator has no slopes by Lemma~\ref{L:additivity} and
so must be a unit in $\Gamma_r[\fp]$ by Corollary~\ref{C:unit}.
That is, we can find $a_{i+1},b_{i+1} \in \Gamma_r[\fp]$ with
$a_{i+1} h_{i+1} + b_{i+1} f_i = 1$. Moreover, by applying
Lemma~\ref{L:principal parts1}, we may choose $a_{i+1},b_{i+1}$
with $w_s(b_{i+1} (g_{i+1} - x_i) f_i) \geq i$ for $s \in I_i$.
(More precisely, apply Lemma~\ref{L:principal parts1} with the roles
of $x$ and $y$ therein played by $h_{i+1}$ and
$b_{i+1}$, respectively; this is valid because $h_{i+1}$ has greatest
slope less than any element of $I_i$.)

Now put $x_{i+1} = x_i + b_{i+1} (g_{i+1} -x_i) f_i$;
then $x_{i+1} - g_i$ is divisible by $f_i$ in $\Gamma_i$,
as then is $x_{i+1} - g_{i+1}$ since
$g_{i+1} - g_i$ is divisible by $f_i$ in $\Gamma_i$.
By Lemma~\ref{L:factor unit2}, $x_{i+1} - g_{i+1}$ is divisible by
$f_i$ also in $\Gamma_{i+1}$.
Since $x_{i+1} - g_{i+1}$ is also
divisible by $h_{i+1}$ in $\Gamma_{i+1}$,
$x_{i+1} - g_{i+1}$ is divisible by $f_{i+1}$ in $\Gamma_{i+1}$.

For any given $s$, we have $w_s(x_{i+1} - x_i) \geq i$ for $i$ large,
so the $x_i$ converge to a limit $g$ in $\Gamma_{\an,r}$.
Since $g - g_i$ is divisible by $f_i$ in $\Gamma_i$,
it is also
divisible by $f$ in $\Gamma_i$. This yields the desired result.
\end{proof}

\begin{remark}
The use of the $f_i$ in the proof of Lemma~\ref{L:principal parts2}
is analogous to the use of ``slope factorizations'' in 
\cite{me-local}. Slope factorizations (convergent products of pure
elements converging to a specified element of $\Gamma_{\an,r}$), 
which are inspired by the comparable construction in \cite{lazard},
will not be used explicitly here; see
\cite[Lemma~3.26]{me-local} for their construction.
\end{remark}

We are finally ready to analyze the B\'ezout property.
\begin{defn}
A \emph{B\'ezout ring/domain} is a ring/domain 
in which every finitely generated ideal
is principal.
Such rings look like principal ideal rings from the point of view
of finitely generated modules.
For instance:
\begin{itemize}
\item
Every $n$-tuple of elements of a B\'ezout ring which generate the
unit ideal is unimodular, i.e., it occurs as the first row of a matrix
of determinant 1 \cite[Lemma~2.3]{me-local}.
\item
The saturated span of any subset of a finite free module over a B\'ezout
domain is a direct summand 
\cite[Lemma~2.4]{me-local}.
\item
Every finitely generated locally free module over
a B\'ezout domain is free \cite[Proposition~2.5]{me-local},
as is every finitely presented torsion-free module
\cite[Proposition~4.9]{crew2}.
\item
Any finitely generated submodule of a finite free module over a B\'ezout
ring is free (straightforward).
\end{itemize}
\end{defn}

\begin{theorem} \label{T:bezout}
For $r \in (0,r_0)$, the ring $\Gamma_{\an,r}$ is a B\'ezout ring
(as then is $\Gancon$). More
precisely, if $J$ is an ideal of $\Gamma_{\an,r}$, the following are
equivalent.
\begin{enumerate}
\item[(a)]
The ideal $J$ is closed.
\item[(b)]
The ideal $J$ is finitely generated.
\item[(c)]
The ideal $J$ is principal.
\end{enumerate}
\end{theorem}
\begin{proof}
Clearly (c) implies both (a) and (b). Also, (a) implies (c) by
Theorem~\ref{T:vector bundles}. It thus suffices to show that
(b) implies (a); by induction, it is enough to check in case $J$
is generated by two nonzero elements $x,y$. Moreover, we may form
the closure of $J$, find a generator $z$, and then divide $x$ and $y$
by $z$; in other words, we may assume that $1$ is in the closure of $J$,
and then what we are to show is that $1 \in J$.

Let $I_1 \subset I_2 \subset \cdots$ be an increasing sequence of
closed subintervals of $(0,r]$, with right endpoints $r$, whose
union is all of $(0,r]$. Then $x$ and $y$ generate the unit ideal in
$\Gamma_i = \Gamma_{I_i}$ for each $i$; that is, we can choose $a_i,
b_i \in \Gamma_i$ with $a_i x + b_i y = 1$.
Note that $b_{i+1} - b_i$ is divisible by $x$ in $\Gamma_i$; by
Lemma~\ref{L:principal parts2}, we can choose 
$b \in \Gamma_{\an,r}$ with $b - b_i$ divisible by $x$ in
$\Gamma_i$ for each $i$. Then $by-1$ is divisible by $x$ in each
$\Gamma_i$, hence also in $\Gamma_{\an,r}$ (by Corollary~\ref{C:subring2}); 
that is, $x$ and $y$
generate the unit ideal in $\Gamma_{\an,r}$, as desired.

We have thus shown that (a), (b), (c) are equivalent, proving that
$\Gamma_{\an,r}$ is a B\'ezout ring. Since $\Gancon$ is the union of 
the $\Gamma_{\an,r}$ for $r \in (0,r_0)$, it is also a B\'ezout ring
because any finitely generated ideal is generated by elements of
some $\Gamma_{\an,r}$.
\end{proof}

\begin{remark} \label{R:no spherical1}
In Lazard's theory (in which $\Gamma_I$ becomes the ring of rigid
analytic functions on the annulus $\log_{|\pi|} |u| \in I$), 
the implication (b)$\implies$(a) is
\cite[Proposition~11]{lazard}, and holds without restriction on
the coefficient field. The implication (a)$\implies$(b)
is equivalent to spherical completeness of the coefficient field
\cite[Th\'eor\`eme~2]{lazard};
however, the analogue here would probably require $K$ also to
be spherically complete (compare Remark~\ref{R:vector bundle}),
which is an undesirable restriction. For instance, it would complicate the
process of descending the slope filtration in 
Chapter~\ref{sec:descents}.
\end{remark}

\section{$\sigma$-modules}
\label{sec:sigma-mod}

We now introduce modules equipped with a semilinear endomorphism
($\sigma$-modules) and study their properties, specifically over
$\Gancon$. In order to highlight the parallels between this theory and the
theory of stable vector bundles (see for instance \cite{shatz}), we have
shaped our presentation along the lines of that of Hartl and 
Pink \cite{hartl-pink}; they study vector bundles with a
Frobenius structure on a punctured disc over a complete nonarchimedean
field of \emph{equal} characteristic $p$, and prove results
very similar to our results over $\Galgancon$.

Beware that our overall sign convention is ``arithmetic''
and not ``geometric''; it thus agrees with the sign conventions of
\cite{katz} (and of \cite{me-local}),
but disagrees with the sign convention of \cite{hartl-pink}
and with the usual convention in the vector bundle setting.

\setcounter{equation}{0}
\begin{remark}
We retain all notation from Chapter~\ref{sec:basic rings}, except that
we redefine the term ``slope''; see Definition~\ref{D:degree slope}.
In particular, $K$ is a field complete with respect to the valuation $v_K$,
and $\Gamma^K$ is assumed to have enough $r_0$-units for some $r_0 > 0$.
Remember that $v_K$ is allowed to be trivial unless otherwise
specified; that means any result
about $\Gancon$ also applies to $\Gamma[\fp]$, unless its statement
explicitly requires $v_K$ to be nontrivial.
\end{remark}

\subsection{$\sigma$-modules}

\begin{defn}
For a ring $R$ containing $\calO$ in which $\pi$ is not a zero divisor,
equipped with a ring endomorphism $\sigma$,
a \emph{$\sigma$-module} 
over a ring $R$ is a finite locally free
$R$-module $M$ equipped with a map $F: \sigma^* M \to M$ (the
\emph{Frobenius action}) which becomes
an isomorphism after inverting $\pi$. (Here $\sigma^* M = M \otimes_{R,\sigma}
R$; that is, view $R$ as a module over itself via $\sigma$, and tensor
it with $M$ over $R$.)
We can view $M$ as a left module over the twisted polynomial
ring $R\{\sigma\}$; we can also view
$F$ as a $\sigma$-linear additive endomorphism of $M$.
A homomorphism of $\sigma$-modules is a module homomorphism equivariant
with respect to the Frobenius actions.
\end{defn}
\begin{remark}
We will mostly consider $\sigma$-modules over B\'ezout rings like
$\Gancon$, in which case there is no harm in replacing
``locally free'' by ``free'' in the definition of a $\sigma$-module.
\end{remark}

\begin{remark}
The category of $\sigma$-modules is typically not abelian (unless
$v_K$ is trivial), because
we cannot form cokernels thanks to the requirement that the underlying modules be locally free.
\end{remark}

\begin{remark}
For any positive integer $a$, $\sigma^a$ is also a Frobenius lift, so
we may speak of $\sigma^a$-modules. This will be relevant when we
want to perform ``restriction of scalars'' in Section~\ref{subsec:pull
push}. However, there is no loss of generality in stating definitions
and theorems in the case $a=1$, i.e., for $\sigma$-modules.
\end{remark}

\begin{defn} \label{D:twist}
Given a $\sigma$-module $M$ of rank $n$ and an integer $c$ (which must be nonnegative if $\pi^{-1} \notin R$), define the
\emph{twist} $M(c)$ of $M$ by $c$ to be the 
module $M$ with the Frobenius action multiplied by $\pi^c$.
(Beware that this definition reflects an earlier choice of normalization,
as in Remark~\ref{R:normalize}, and a choice of a sign convention.)
If $\pi$ is invertible in $R$,
define the \emph{dual} $M^\dual$ of $M$ to be the $\sigma$-module
$\Hom_R(M, R) \cong (\wedge^{n-1} M) \otimes (\wedge^n M)^{\otimes -1}$
and the \emph{internal hom} of $M,N$ as $M^\dual \otimes N$.
\end{defn}

\begin{defn}
Given a $\sigma$-module $M$ over a ring $R$, let
$H^0(M)$ and $H^1(M)$ denote the kernel and cokernel, respectively,
of the map $F-1$ on $M$; note that if $N$ is another $\sigma$-module,
then there is a natural bilinear map $H^0(M) \times H^1(N) \to H^1(M \otimes N)$. Given two $\sigma$-modules $M_1$ and $M_2$ over $R$,
put $\Ext(M_1,M_2) = H^1(M_1^\dual \otimes M_2)$; a standard homological
calculation
(as in \cite[Proposition~2.4]{hartl-pink}) shows that $\Ext(M_1,M_2)$ 
coincides with the Yoneda $\Ext^1$ in this category. That is,
$\Ext(M_1,M_2)$ classifies
short exact sequences $0 \to M_2 \to M \to M_1 \to 0$ of $\sigma$-modules
over $R$, up to isomorphisms 
\[
\xymatrix{
0 \ar[r] & M_2 \ar[r] \ar[d] & M \ar[r] \ar[d] & M_1 \ar[r] \ar[d] & 0 \\
0 \ar[r] & M_2 \ar[r] & M \ar[r] & M_1 \ar[r] & 0
}
\]
which induce the identity maps on $M_1$ and $M_2$.
\end{defn}

\subsection{Restriction of Frobenius}
\label{subsec:pull push}

We now introduce two functors analogous to those induced by the
``finite maps'' in \cite[Section~7]{hartl-pink}. Beware that the analogy is 
not perfect; see Remark~\ref{R:bad analogy}.
\begin{defn} \label{D:push pull}
Fix a ring $R$ equipped with an endomorphism $\sigma$.
For $a$ a positive integer,
let $[a]: R\{\sigma^{a}\} \to R\{\sigma\}$ be the
natural inclusion homomorphism of twisted polynomial rings.
Define the \emph{$a$-pushforward functor}
$[a]_*$, from $\sigma$-modules to $\sigma^{a}$-modules, to be the
restriction functor along $[a]$.
Define the \emph{$a$-pullback functor} $[a]^*$, from $\sigma^{a}$-modules
to $\sigma$-modules, to be the extension of scalars functor
\[
M \mapsto R\{\sigma\} \otimes_{R\{\sigma^{a}\}} M.
\]
Note that $[a]^*$ and $[a]_*$ are left and right adjoints of each other.
Also, $[a]^* [a']^* = [aa']^*$ and $[a]_* [a']_* = [aa']_*$.
Furthermore, $[a]_* (M(c)) = ([a]_* M)(ac)$.
\end{defn}
\begin{remark} \label{R:bad analogy}
There are some discrepancies in the analogy with \cite{hartl-pink},
due to the fact that there the corresponding map $[a]$ is actually
a homomorphism of the underlying ring, rather than a change of Frobenius.
The result is that some (but not all!) of the properties of the pullback and 
pushforward are swapped between here and \cite{hartl-pink}.
For an example of this mismatch in action, see
Proposition~\ref{P:pull basic}.
\end{remark}

\begin{remark}
The functors $[a]$ will ultimately serve to rescale the slopes of a
$\sigma$-module; using them makes it possible to
avoid the reliance in \cite[Chapter~4]{me-local} on making extensions
of $\calO$. Among other things, this lets us get away with normalizing
$w$ in terms of the choice of $\calO$, since we will not have to change
that choice at any point.
\end{remark}

\begin{lemma} \label{L:push pull}
For any positive integer $a$ and any integer $c$,
$[a]_* [a]^* (R(c)) \cong R(c)^{\oplus a}$.
\end{lemma}
\begin{proof}
We can write
\[
[a]_* [a]^* (R(c)) \cong
\oplus_{i=0}^{a-1} \{\sigma^{i}\} (R(c)),
\]
where on the right side $R\{\sigma^{a}\}$ acts separately
on each factor. Hence the claim follows.
(Compare \cite[Proposition~7.4]{hartl-pink}.)
\end{proof}

\begin{lemma} \label{L:pull isom}
Suppose that the residue field of $\calO$ contains an algebraic
closure of $\FF_q$.
For $i$ a positive integer, let $L_i$ be the fixed field of $\calO[\fp]$
under $\sigma^i$.
\begin{enumerate}
\item[(a)]
For any $\sigma$-module $M$ over $\Gancon$ and any positive integer $a$,
$H^0([a]_* M) = H^0(M) \otimes_{L_1} L_a$.
\item[(b)]
For any $\sigma$-modules $M$ and $N$ over $\Gancon$ and any positive
integer $a$, $M \cong N$ if and only if $[a]_* M \cong
[a]_* N$.
\end{enumerate}
\end{lemma}
\begin{proof}
\begin{enumerate}
\item[(a)]
It suffices to show that $H^0([a]_* M)$ admits a basis invariant under
the induced action of $\sigma$.
Since $\sigma$ generates $\Gal(L_a/L_1)$, this follows from Hilbert's
Theorem 90.

\item[(b)]
A morphism $[a]_* M \to [a]_* N$ corresponds to
an element of $V = H^0(([a]_* M)^\dual \otimes [a]_* N) \cong
H^0([a]_* (M^\dual \otimes N))$, which by (a) coincides with
$H^0(M^\dual \otimes N) \otimes_{L_1} L_a$. If there is an isomorphism
$[a]_* M \cong [a]_* N$, then the determinant locus on $V$ is not all of $V$;
hence the same is true on $H^0(M^\dual \otimes N)$.
We can thus find an $F$-invariant element of $M^\dual \otimes N$
corresponding to an isomorphism $M \cong N$.
(Compare \cite[Propositions~7.3 and~7.5]{hartl-pink}.)
\end{enumerate}
\end{proof}

\subsection{$\sigma$-modules of rank 1}

In this section, we analyze some $\sigma$-modules of rank 1 over $\Gancon$;
this amounts to solving some simple equations involving $\sigma$, as in 
\cite[Proposition~3.19]{me-local} (compare also
\cite[Propositions~3.1 and~3.3]{hartl-pink}).

\begin{defn}
Define the twisted powers $\pi^{\{m\}}$ of $\pi$ by the two-way recurrence
\[
\pi^{\{0\}} = 1, \qquad \pi^{\{m+1\}} = (\pi^{\{m\}})^\sigma \pi.
\]
\end{defn}

First, we give a classification result.
\begin{prop} \label{P:rank 1 class}
Let $M$ be a $\sigma$-module of rank $1$ over $R$, for
$R$ one of $\Galgcon, \Galgcon[\fp], \Galgancon$.
Then there exists a unique integer $n$, which is 
nonnegative if $R = \Galgcon$, such that $M \cong R(n)$.
\end{prop}
\begin{proof}
Let $\bv$ be a generator of $M$, and write $F\bv = x \bv$. Then $x$
must be a unit in $R$, so by Corollary~\ref{C:unit}, $x \in \Galgcon[\fp]$
and so $w(x)$ is defined. If $M \cong R(n)$, we must then have $n = w(x)$;
hence $n$ is unique if it exists.

Since the residue field of $\Galgcon$ is algebraically closed,
we can find a unit $u$ in $\Galgcon$ such that
$u^\sigma \pi^{-n} x \equiv u \pmod{\pi}$.
Choose $r>0$ such that $w_r(u^\sigma \pi^{-n} x/u - 1) > 0$; then there exists
a unit $y \in \Galg$ with $u^\sigma \pi^{-n} x/u = y^\sigma/y$,
and a direct calculation (by induction on $m$) shows that
$v_{m,r}(y) > 0$ for all $m>0$.
(For details, see the proof of Lemma~\ref{L:negative eigen}.)
Hence $y$ is a unit in $\Galgcon$, and so
$\bw = (u/y) \bv$ is a generator of $M$ satisfying
$F\bw = \pi^n \bw$. Thus there exists an isomorphism $M \cong R(n)$.
\end{proof}

We next compute some instances of $H^0$.
\begin{lemma} \label{L:h0h1}
Let $n$ be a nonnegative integer. If $x \in \Gancon$ and
$x - \pi^n x^\sigma \in \Gcon[\fp]$, then $x \in \Gcon[\fp]$.
\end{lemma}
\begin{proof}
Suppose the contrary; put $y = x - \pi^n x^\sigma$.
We can find $m$ with $v_m(y) = \infty$ and $0 < v_m(x) < \infty$,
since both hold for $m$ sufficiently small by
Corollary~\ref{C:subring}. Then
\[
v_m(x) > q^{-1} v_m(x) 
= q^{-1} v_m(\pi^n x^\sigma)
= q^{-1} v_{m-n}(x^\sigma)
= v_{m-n}(x) 
\geq v_m(x),
\]
contradiction. 
\end{proof}

\begin{prop} \label{P:h0}
Let $n$ be an integer.
\begin{enumerate}
\item[(a)]
If $n=0$, then $H^0(\Gcon[\fp](n)) = H^0(\Gancon(n)) \neq 0$;
moreover, any nonzero element of $H^0(\Gancon)$ is a unit in $\Gancon$.
\item[(b)]
If $n > 0$, then $H^0(\Gancon(n)) = 0$.
\item[(c1)]
If $n<0$ and $v_K$ is trivial, then $H^0(\Gancon(n)) = 0$.
\item[(c2)]
If $n<0$, $v_K$ is nontrivial, and $K$ is perfect, then
$H^0(\Gancon(n)) \neq 0$.
\end{enumerate}
\end{prop}
\begin{proof}
The group $H^0(R(n))$ consists of those $x \in R$ with
\begin{equation} \label{eq:h0}
\pi^n x^\sigma = x,
\end{equation}
so our assertions are all really about the solvability of this equation.
\begin{enumerate}
\item[(a)] 
If $n=0$, then $H^0(\Gcon[\fp](n))= H^0(\Gancon(n))$ by
Lemma~\ref{L:h0h1}, and the former equals the 
fixed field of $\calO[\fp]$ under $\sigma$.
\item[(b)] 
By Lemma~\ref{L:h0h1}, any solution $x$ of \eqref{eq:h0} over $\Gancon$
actually belongs to $\Gcon[\fp]$. In particular, $w(x) = w(x^\sigma)$
is well-defined. But 
\eqref{eq:h0} yields $w(x^\sigma) + n = w(x)$, which for $n>0$
forces $x=0$.
\item[(c1)]
If $n<0$ and $v_K$ is trivial, then
$w(x) = w(x^\sigma)$ is well-defined, but
\eqref{eq:h0} yields $w(x^\sigma) + n = w(x)$, which forces $x=0$.
\item[(c2)]
If $n<0$, $v_K$ is nontrivial, and $K$ is perfect,
we may pick $\overline{u} \in K$ with
$v_K(\overline{u}) > 0$ (since $v_K$ is nontrivial), and then set
$x$ to be the limit of the convergent series
\[
\sum_{m \in \ZZ} (\pi^{\{m\}})^n [\overline{u}^{q^m}]
\]
to obtain a nonzero solution of \eqref{eq:h0}.
\end{enumerate}
\end{proof}

\begin{remark}
If $n<0$, $v_K$ is nontrivial, and $K$ is not perfect, then
the size of $H^0(\Gancon(n))$ depends on the particular choice of
the Frobenius lift $\sigma$ on $\Gancon$. For instance, in
the notation of Section~\ref{subsec:robba}, if $\sigma$ is a so-called
``standard Frobenius lift'' sending $u$ to $u^q$, then 
any solution $x = \sum x_i u^i$ of \eqref{eq:h0} must have
$x_i = 0$ whenever $i$ is not divisible by $q$. By the same token,
$x_i = 0$ whenever $i$ is not divisible by $q^2$, or by $q^3$, and so on;
hence we must have $x \in \calO[\fp]$, which as in (b) above is impossible
for $n > 0$. On the other hand, if $u^\sigma = (u+1)^q - 1$, then
$x = \log(1+u) \in \Gancon$ satisfies $x^\sigma = qx$; indeed, the existence of
such an $x$ is a backbone of the theory of $(\Phi, \Gamma)$-modules
associated to $p$-adic Galois representations, as in 
\cite{berger-weak}.
\end{remark}

We next consider $H^1$. 
\begin{prop} \label{P:h1}
Let $n$ be an integer.
\begin{enumerate}
\item[(a)]
If $n=0$ and $K$ is separably closed, then $H^1(\Gcon[\fp](n)) = H^1(\Gancon(n)) = 0$.
\item[(b1)]
If $n \geq 0$, then the map $H^1(\Gcon[\fp](n)) \to H^1(\Gancon(n))$ is injective.
\item[(b2)]
If $n > 0$ and $v_K$ is trivial, then $H^1(\Gancon(n)) = 0$.
\item[(b3)]
If $n > 0$, $v_K$ is nontrivial, and $K$ is perfect,
then $H^1(\Gancon(n)) \neq 0$, with a nonzero element given by
$[\overline{x}]$ for any $\overline{x} \in K$ with
$v_K(\overline{x}) < 0$.
\item[(c)]
If $n < 0$ and $K$ is perfect, then $H^1(\Gcon[\fp](n)) = 
H^1(\Gancon(n)) = 0$.
\end{enumerate}
\end{prop}
\begin{proof}
The group $H^1(R(n))$ consists of the quotient of the additive
group of $R$ by the subgroup of those $x \in R$ for which the
equation
\begin{equation} \label{eq:h1}
x = y - \pi^n y^\sigma
\end{equation}
has a solution $y \in R$,
so our assertions are all really about the solvability of this equation.
\begin{enumerate}
\item[(a)]
If $n = 0$ and $K$ is separably closed, then for each $x \in \Gamma$,
there exists $y \in \Gamma$ such that $x \equiv y - y^q \pmod{\pi}$.
By iterating this construction, we can produce for any $x \in \Gamma[\fp]$
an element $y \in \Gamma[\fp]$ satisfying \eqref{eq:h1}, such that
\[
v_m(y) \geq \min\{v_m(x), v_m(x)/q\};
\]
in particular, if $x \in \Gcon[\fp]$, then $y \in \Gcon[\fp]$. Moreover,
given $x \in \Gancon$, we can write $x$ as a convergent series of elements
of $\Gcon[\fp]$, and thus produce a solution of \eqref{eq:h1}.
Hence $H^1(\Gcon(n)) = H^1(\Gancon(n)) = 0$.
\item[(b1)]
This follows at once from Lemma~\ref{L:h0h1}.
\item[(b2)]
If $n > 0$ and $v_K$ is trivial, then
for any $x \in \Gancon$, the series
\[
y = \sum_{m=0}^\infty (\pi^{\{m\}})^n x^{\sigma^m}
\]
converges $\pi$-adically to a solution of \eqref{eq:h1}.
\item[(b3)]
By (b1), it suffices to show that
$x = [\overline{x}]$ represents a nonzero element of $H^1(\Gcon[\fp](n))$.
By (b2), there exists a unique $y \in \Gamma[\fp]$ satisfying
\eqref{eq:h1}; however, we have
$v_{mn}(y) = q^m v_K(\overline{x})$ for all $m \geq 0$, 
and so $y \notin \Gcon[\fp]$.
\item[(c)]
Let $x = \sum_i [\overline{x_i}] \pi^i$ be the Teichm\"uller presentation
of $x$. Pick $c>0$, and 
let $z_1$ and $z_2$ be the sums of $[\overline{x_i}] \pi^i$
over those $i$ for which $v_K(\overline{x_i}) < c$ and $v_K(\overline{x_i})
\geq c$, respectively. Then the sums
\begin{align*}
y_1 &= \sum_{m=0}^\infty -(\pi^{\{-m-1\}})^n z_1^{\sigma^{-m-1}} \\
y_2 &= \sum_{m=0}^\infty (\pi^{\{m\}})^n z_2^{\sigma^m}
\end{align*}
converge to solutions of $z_1 = y_1 - \pi^n y_1^\sigma$ and
$z_2 = y_2 - \pi^n y_2^\sigma$, respectively. Hence
$y = y_1 + y_2$ is a solution of \eqref{eq:h1}. Moreover, if
$x \in \Gcon[\fp]$, then $\overline{x_i} = 0$ for $i$ sufficiently small,
so we can choose $c$ to ensure $z_2 = 0$; then $y = y_1 \in \Gcon[\fp]$.
\end{enumerate}
\end{proof}

\subsection{Stability and semistability}

As in \cite{hartl-pink},
we can set up a formal analogy between the study of $\sigma$-modules
over $\Gancon$ and the study of stability of vector bundles.

\begin{defn} \label{D:degree slope}
For $M$ a $\sigma$-module of rank 1 over 
$\Gancon$ generated by some $\bv$, define the
\emph{degree} of $M$, denoted $\deg(M)$,
to be the unique integer $n$ such that $M \otimes \Galgancon \cong
\Galgancon(n)$, as provided by Proposition~\ref{P:rank 1 class};
concretely, $\deg(M)$ is the valuation of the unit via which
$F$ acts on a generator of $M$.
For a $\sigma$-module $M$ over 
$\Gancon$ of rank $n$, define
$\deg(M) = \deg(\wedge^n M)$. Define $\mu(M) = \deg(M)/\rank(M)$;
we refer to $\mu(M)$ as the
\emph{slope} of $M$ (or as the \emph{weight} of $M$, per
the terminology of \cite[Section~6]{hartl-pink}).
\end{defn}

\begin{lemma} \label{L:same rank}
Let $M$ be a $\sigma$-module over $\Gancon$, and let $N$ be
a $\sigma$-submodule of $M$ with $\rank(M) = \rank(N)$. Then
$\deg(N) \geq \deg(M)$, with equality if and only if $M=N$;
moreover, equality must hold if $v_K$ is trivial.
\end{lemma}
\begin{proof}
By taking exterior powers, it suffices to check this for
$\rank M = \rank N = 1$; also, there is no harm in 
assuming that $K$ is algebraically closed.
By Proposition~\ref{P:rank 1 class},
$M \cong \Galgancon(c)$ and $N \cong \Galgancon(d)$ for some integers
$c,d$. By twisting, we may reduce to the case $d=0$.
Then by Proposition~\ref{P:h0}, we have $0 \geq c$, with
equality forced if $v_K$ is trivial; moreover,
if $c=0$, then $N$ contains a generator which also belongs to
$H^0(M)$. But every nonzero element of the latter also generates $M$,
so $c=0$ implies $M=N$.
(Compare \cite[Proposition~6.2]{hartl-pink}.)
\end{proof}

\begin{lemma} \label{L:degree adds}
If $0 \to M_1 \to M \to M_2 \to 0$ is a short exact sequence
of $\sigma$-modules over $\Gancon$, then 
$\deg(M) = \deg(M_1) + \deg(M_2)$.
\end{lemma}
\begin{proof}
Put $n_1 = \rank(M_1)$ and $n_2 = \rank(M_2)$. Then
the claim follows from the existence of the isomorphism
\[
\wedge^{n_1+n_2} M \cong \left( \wedge^{n_1} M_1 \right)
\otimes \left( \wedge^{n_2} M_2 \right)
\]
of $\sigma$-modules.
\end{proof}

\begin{prop} \label{P:pull basic}
Let $a$ be a positive integer, let $M$ be a $\sigma$-module over $\Gancon$,
and let $N$ be a $\sigma^a$-module over $\Gancon$.
\begin{enumerate}
\item[(a)] $\deg([a]_* M) = a \deg(M)$ and $\deg([a]^* N) = \deg(N)$.
\item[(b)] $\rank([a]_* M) = \rank(M)$ and $\rank([a]^* N) = a \rank(N)$.
\item[(c)] $\mu([a]_* M) = a \mu(M)$ and $\mu([a]^* N) = \frac{1}{a} \mu(N)$.
\end{enumerate}
\end{prop}
\begin{proof}
Straightforward (compare \cite[Proposition~7.1]{hartl-pink}, but note
that the roles of the pullback and pushforward are interchanged here).
\end{proof}

\begin{defn} \label{D:semistable}
We say a $\sigma$-module $M$ over $\Gancon$ is \emph{semistable} if
$\mu(M) \leq \mu(N)$ for any nonzero $\sigma$-submodule
$N$ of $M$. We say $M$ is \emph{stable} if
$\mu(M) < \mu(N)$ for any nonzero proper $\sigma$-submodule
$N$ of $M$. Note that the direct sum of semistable $\sigma$-modules
of the same slope is also semistable.
By Proposition~\ref{P:pull basic}, for any positive integer $a$,
if $[a]_* M$ is (semi)stable, then
$M$ is (semi)stable.
\end{defn}
\begin{remark}
As noted earlier, the inequalities are reversed from the usual
definitions of stability and semistability for vector bundles, because
of an overall choice of sign convention.
\end{remark}

\begin{remark}
Beware that this use of the term ``semistable'' is only distantly related
to its use to describe $p$-adic Galois representations!
\end{remark}

\begin{lemma} \label{L:semistable}
For any integer $c$ and any positive integer $n$,
the $\sigma$-module $\Gancon(c)^{\oplus n}$ is semistable of 
slope $c$.
\end{lemma}
\begin{proof}
There is no harm 
in assuming that $K$ is algebraically closed,
and that $v_K$ is nontrivial.
Let $N$ be a nonzero $\sigma$-submodule of $M$ of rank $d'$
and degree $c'$. Then
\[
\Galgancon(c') \cong \wedge^{d'} N \subseteq \wedge^{d'} \Galgancon(c)^{\oplus n}
\cong \Galgancon(cd')^{\oplus \binom{n}{d'}}.
\]
In particular, $H^0(\Galgancon(cd' - c')) \neq 0$;
by Proposition~\ref{P:h0}, this implies $c' \geq cd'$,
yielding semistability.
(Compare \cite[Proposition~6.3(b)]{hartl-pink}.)
\end{proof}

\begin{lemma} \label{L:push semi}
For any positive integer $a$ and any integer $c$,
the $\sigma$-module $[a]^* (\Gancon(c))$ over $\Gancon$
is semistable of rank $a$, degree $c$, and 
slope $c/a$. Moreover, if $a$ and $c$ are coprime, then
$[a]^*(\Gancon(c))$ is stable.
\end{lemma}
\begin{proof}
By Lemma~\ref{L:push pull} and Lemma~\ref{L:semistable},
$[a]_* [a]^*(\Gancon(c))$ is semistable of rank $a$ and slope $c$;
as noted in Definition~\ref{D:semistable}, it follows that
$[a]^* (\Gancon(c))$ is semistable.

Let $M$ be a nonzero $\sigma$-submodule of $[a]^* (\Gancon(c))$.
If $a$ and $c$ are coprime, then $\deg(M) = (c/a)
\rank(M)$ is an integer; since $\rank(M) \leq a$, this is only
possible for $\rank(M) = a$.
But then Lemma~\ref{L:same rank} implies that $\mu(M) > c/a$
unless $M = [a]^*(\Gancon(c))$.
We conclude that $[a]^*(\Gancon(c))$ is stable.
(Compare \cite[Proposition~8.2]{hartl-pink}.)
\end{proof}

\subsection{Harder-Narasimhan filtrations}

Using the notions of degree and slope, we can make the usual formal construction of Harder-Narasimhan filtrations, with its usual properties.

\begin{defn}
Given a multiset $S$ of 
$n$ real numbers, define the \emph{Newton polygon} of $S$
to be the graph of the piecewise linear function on $[0,n]$
sending 0 to 0, whose slope on $[i-1,i]$ is the $i$-th smallest element
of $S$; we refer to the point on the graph corresponding to the image of $n$ as the \emph{endpoint} of the polygon. Conversely, given such a graph, define its \emph{slope multiset} to be the
slopes of the piecewise linear function on $[i-1,i]$ for $i=1, \dots, n$.
We say that the Newton polygon of $S$ \emph{lies above} the Newton polygon
of $S'$ if no vertex of the
polygon of $S$ lies below the polygon of $S'$, and the two polygons
have the same endpoint.
\end{defn}

\begin{defn}
Let $M$ be a $\sigma$-module over $\Gancon$. A \emph{semistable filtration}
of $M$ is an exhaustive filtration
$0 = M_0 \subset M_1 \subset \cdots \subset M_l = M$ of $M$
by saturated $\sigma$-submodules, such that each successive quotient 
$M_i/M_{i-1}$ is semistable of some slope $s_i$.
A \emph{Harder-Narasimhan
filtration} (or \emph{HN-filtration}) of $M$ is a semistable filtration
with $s_1 < \cdots < s_l$.
An HN-filtration is unique if it exists, as $M_1$ can then be characterized
as the unique maximal $\sigma$-submodule of $M$ of minimal slope,
and so on. 
\end{defn}

\begin{defn}
Let $M$ be a $\sigma$-module over $\Gancon$. 
Given a semistable
filtration $0 = M_0 \subset M_1 \subset \cdots \subset M_l = M$ of $M$,
form the multiset consisting of, for $i=1, \dots, l$, the slope
$\mu(M_i/M_{i-1})$ with multiplicity $\rank(M_i/M_{i-1})$. We call this
the \emph{slope multiset} of the filtration, and we call the associated
Newton polygon the \emph{slope polygon} of the filtration.
If $M$ admits a Harder-Narasimhan filtration, we refer to the slope multiset as the \emph{Harder-Narasimhan slope multiset} (or \emph{HN-slope multiset}) of
$M$, and to the Newton polygon as the \emph{Harder-Narasimhan polygon}
(or \emph{HN-polygon}) of $M$.
\end{defn}

\begin{prop} \label{P:filtrations}
Let $M$ be a $\sigma$-module over $\Gancon$ admitting an HN-filtration.
Then the HN-polygon lies above the slope polygon of any semistable
filtration of $M$. 
\end{prop}
\begin{proof}
Let $0 = M_0 \subset M_1 \subset \cdots \subset M_l = M$ be an HN-filtration,
and let $0 = M'_0 \subset M'_1 \subset \cdots \subset M'_m = M$ be a semistable
filtration. To prove the inequality, it 
suffices to prove that for each of $i = 1,\dots, l$, we can choose
$\rank(M/M_i)$ slopes from the slope multiset of the semistable filtration
whose sum is greater than or equal to the
sum of the greatest $\rank(M/M_i)$ HN-slopes of $M$; note that the latter
is just $\deg(M/M_i)$.

For $l = 1, \dots, m$, put
\[
d_l = \rank(M'_l + M_i) - \rank(M'_{l-1} + M_i) \leq \rank(M'_l/M'_{l-1}).
\]
Since $M'_l/M'_{l-1}$ is semistable and $(M'_l + M_i)/(M'_{l-1}+M_i)$
is a quotient of $M'_l/M'_{l-1}$, we have
\[
\mu(M'_l/M'_{l-1}) \geq \mu((M'_l + M_i)/(M'_{l-1} + M_i)).
\]
However, we also have
\[
\sum_{i=1}^l d_l \mu((M'_l + M_i)/(M'_{l-1} + M_i)) = \deg(M/M_i),
\]
yielding the desired inequality.
\end{proof}

\begin{remark}
We will ultimately prove that 
every $\sigma$-module over $\Gancon$ admits a 
Harder-Narasimhan filtration (Proposition~\ref{P:HN exists}),
and that the successive quotients become isomorphic over
$\Galgancon$ to direct sums of ``standard'' $\sigma$-modules
of the right slope (Theorems~\ref{T:isoclinic lattice} and~\ref{T:descend special}).
Using the formalism of
Harder-Narasimhan filtrations makes it a bit more convenient to
articulate the proofs of these assertions.
\end{remark}

\subsection{Descending subobjects}

We will ultimately be showing that the formation of a
Harder-Narasimhan filtration of a $\sigma$-module over $\Gancon$
commutes with base change. In order to prove this sort of statement,
it will be useful to have a bit of terminology.
\begin{defn} \label{D:descent}
Given an injection $R \hookrightarrow S$ of integral domains equipped with 
compatible
endomorphisms $\sigma$, a $\sigma$-module $M$ over $R$,
and a saturated $\sigma$-submodule $N_S$ of $M_S = M \otimes_R S$, we say that 
\emph{$N_S$ descends to $R$} if there is a saturated $\sigma$-submodule
$N$ of $M$ such that $N_S = N \otimes_R S$; note that $N$ is 
unique if it exists,
because it can be characterized as $M \cap N_S$. Likewise, given a 
filtration of $M_S$,
we say the filtration descends to $R$ if it is induced by a filtration of $M$.
\end{defn}

The following lemma lets us reduce most descent questions to consideration of
submodules of rank 1.
\begin{lemma} \label{L:wedge descent}
With notation as in Definition~\ref{D:descent}, suppose that $R$ is
a B\'ezout domain, and put $d = \rank N_S$.
Then $N_S$ descends to $R$ if and only if $\wedge^d N_S \subseteq 
(\wedge^d M)_S$
descends to $R$.
\end{lemma}
\begin{proof}
If $N_S = N \otimes_R S$, then $\wedge^d N_S = (\wedge^d N) \otimes_R S$ descends to
$R$. Conversely, if $\wedge^d N_S = (N') \otimes_R S$, let $N$ be the $\sigma$-submodule
of $M$ consisting of those $\bv \in M$ such that $\bv \wedge \bw = 0$ for all
$\bw \in N'$. Then $N$ is saturated and $N \otimes_R S = N_S$, since $N$ is defined
by linear conditions which in $M_S$ cut out precisely $N_S$. Since $R$ is a B\'ezout
domain, this suffices to ensure that $N$ is free; since $N$ is visibly
stable under $F$, $N$ is in fact a $\sigma$-submodule of $M$, and so
$N_S$ descends to $R$.
\end{proof}

\section{Slope filtrations of $\sigma$-modules}
\label{sec:dm}

In this chapter, we give a classification of $\sigma$-modules over
$\Galgancon$, as in \cite[Chapter~4]{me-local}. However, this
presentation looks somewhat different, mainly because of the formalism
introduced in the previous chapter. We also have integrated into a
single presentation the cases where $v_K$ is nontrivial and where $v_K$
is trivial; these are presented separately in 
\cite{me-local} (in Chapters~4 and~5 respectively).
These two cases do have different flavors, which we will point out as 
we go along.

Beware that although we have mostly made the exposition self-contained, there remains one notable exception: we do not repeat the key calculation made in
\cite[Lemma~4.12]{me-local}. 
See Lemma~\ref{L:pairing} for the relevance of this calculation.

\setcounter{equation}{0}
\begin{convention}
To lighten the notational load, we write $\calR$ for
$\Galgancon$. Whenever working over $\calR$, we also make the harmless
assumption that  $\pi$ is
$\sigma$-invariant.
\end{convention}

\subsection{Standard $\sigma$-modules}

Following \cite[Section~8]{hartl-pink}, we introduce the standard
building blocks into which we will decompose $\sigma$-modules over
$\calR$.

\begin{defn} \label{D:standards}
Let $c,d$ be coprime integers with $d > 0$.
Define the $\sigma$-module $M_{c,d} = [d]^*(\calR(c))$ over $\calR$; that is,
$M_{c,d}$ is freely generated by $\be_1, \dots, \be_d$ with
\[
F\be_1 = \be_2, \quad \dots, \quad F\be_{d-1} = \be_d, \quad
F\be_d = \pi^{c} \be_1.
\]
This $\sigma$-module is stable of slope $c/d$ by Lemma~\ref{L:push semi}.
We say a $\sigma$-module $M$ is \emph{standard} if it is isomorphic
to some $M_{c,d}$; in that case, we say
a basis $\be_1, \dots, \be_d$ as above is a \emph{standard basis} of $M$.
\end{defn}

\begin{lemma} \label{L:standard arith}
\begin{enumerate}
\item[(a)]
$M_{c,d} \otimes M_{c',d'} \cong M_{c'',d''}^{\oplus dd'/d''}$,
where $c/d + c'/d' = c''/d''$ in lowest terms.
\item[(b)]
$M_{c,d}(c') \cong M_{c+c'd,d}$.
\item[(c)]
$M_{c,d}^\dual \cong M_{-c,d}$.
\end{enumerate}
\end{lemma}
\begin{proof}
To verify (a), it is enough to do so after applying $[dd']_*$
thanks to Lemma~\ref{L:pull isom}. Then the desired isomorphism
follows from Lemma~\ref{L:push pull}. Assertion (b) follows from (a),
and (c) follows from the explicit description of $M_{c,d}$ given above.
(Compare \cite[Proposition~8.3]{hartl-pink}.)
\end{proof}

\begin{prop} \label{P:pushforward}
Let $c,d$ be coprime integers with $d > 0$.
\begin{enumerate}
\item[(a)]
The group $H^0(M_{c,d})$ is nonzero if and only if $v_K$ is nontrivial and
$c/d \leq 0$, or $v_K$ is trivial and $c/d = 0$.
\item[(b)]
The group $H^1(M_{c,d})$ is nonzero if and only if $v_K$ is nontrivial
and $c/d > 0$.
\end{enumerate}
\end{prop}
\begin{proof}
These assertions follow from
Propositions~\ref{P:h0} and~\ref{P:h1},
plus the fact that $H^i([d]^* M) \cong H^i(M)$ for $i=0,1$.
\end{proof}
\begin{cor} \label{C:hom ext2}
Let $c,c',d,d'$ be integers, with $d,d'$ positive and 
$\gcd(c,d) = \gcd(c',d') = 1$.
\begin{enumerate}
\item[(a)]
We have $\Hom(M_{c',d'}, M_{c,d}) \neq 0$ if and only if $v_K$ is
nontrivial and $c'/d' \geq c/d$, or $v_K$ is trivial and $c'/d' = c/d$.
\item[(b)]
We have $\Ext(M_{c',d'}, M_{c,d}) \neq 0$ if and only if $v_K$ is
nontrivial and $c'/d' < c/d$.
\end{enumerate}
\end{cor}

\begin{remark}
One can show that $\End(M_{c,d})$ is a division algebra (this can be deduced
from the fact that $M_{c,d}$ is stable) and even 
describe it explicitly, as in 
\cite[Proposition~8.6]{hartl-pink}. For our purposes, it will be enough
to check that $\End(M_{c,d})$ is a division algebra after establishing
the existence of Dieudonn\'e-Manin decompositions; see
Corollary~\ref{C:irreducible}.
\end{remark}

\subsection{Existence of eigenvectors}

In classifying $\sigma$-modules over $\calR$, it is useful to employ
the language of ``eigenvectors''.
\begin{defn}
Let $d$ be a positive integer.
A \emph{$d$-eigenvector} (or simply \emph{eigenvector} if $d=1$)
of a $\sigma$-module $M$ over $\calR$
is a nonzero element $\bv$ of $M$ such that $F^d \bv = \pi^c \bv$ for some
integer $c$. We refer to the quotient $c/d$ as the \emph{slope} of $\bv$.
\end{defn}

\begin{prop} \label{P:eigenvector}
Suppose that $v_K$ is nontrivial. Then every nontrivial $\sigma$-module
over $\calR$ contains an eigenvector.
\end{prop}
\begin{proof}
The calculation is basically that of \cite[Proposition~4.8]{me-local}:
use the fact that $F$ makes things with ``very positive partial valuations''
converge better whereas $F^{-1}$ makes things with
``very negative partial valuations'' converge better. However, 
one can simplify the final analysis a bit, as is done in
\cite[Theorem~4.1]{hartl-pink}.

We first set some notation as in \cite[Proposition~4.8]{me-local}.
Let $M$ be a nontrivial $\sigma$-module over $\calR$.
Choose a basis $\be_1, \dots, \be_n$ of $M$, and define the 
invertible $n \times n$ matrix $A$ over $\calR$ by the equation
$F\be_j = \sum_i A_{ij} \be_i$.
Choose $r>0$ such that $A$ and its inverse have entries in $\Galg_{\an,r}$.
Choose $\epsilon > 0$, and choose
an integer $c$ with 
$c \leq \min\{w_r(A), w_r((A^{-1})^{\sigma^{-1}})\} - \epsilon$.
Choose an integer $m > c$ such that
the interval
\[
\left(
\frac{-c + m}{(q-1)r},
\frac{q(c + m)}{(q-1)r}
\right)
\]
is nonempty (true for $m$ sufficiently large), 
and choose $d$ in the intersection of that interval
with the image of $v_K$. (The choice of $d$ is
possible because $v_K$ is nontrivial and
the residue field of $\Galgcon$
is algebraically closed, so the image of $v_K$ contains a copy
of $\QQ$ and hence is dense in $\RR$.)

For an interval $I \subseteq (0,r]$,
define the functions $a,b: \Gamma_{I} \to \Gamma_{I}$ as follows.
For $x \in \Gamma_{[r,r]}$, let $\sum_{j \in \ZZ} [\overline{x_j}] \pi^j$
be the Teichm\"uller presentation of $x$ (as in
Definition~\ref{D:teich}). Let $a(x)$ and $b(x)$ be the sums of
$[\overline{x_j}] \pi^j$ over those $j$ with $v_K(\overline{x_j}) < d$
and $v_K(\overline{x_j}) \geq d$, respectively.
We may think of $a$ and $b$ as splitting $x$ into ``negative'' and ``positive''
terms.
(This decomposition shares its canonicality with the corresponding
decomposition in \cite[Theorem~4.1]{hartl-pink} but not with the one in 
\cite[Proposition~4.8]{me-local}.)

For $I \subseteq (0,r]$, let $M_I$ denote the 
$\Gamma_I$-span of the $\be_i$.
For $\bv \in M_{I}$, 
write $\bv = \sum_i x_i \be_i$, and for $s \in I$, define
\begin{align*}
v_{n,s}(\bv) &= \min_i\{v_{n,s}(x_i)\} \\
w_s(\bv) &= \min_i \{w_s(x_i)\}.
\end{align*}
Put $a(\bv) = \sum_i a(x_i) \be_i$ and $b(\bv) = \sum_i b(x_i) \be_i$;
then by the choice of $d$, we have for $\bv \in M_{(0,r]}$,
\begin{align*}
w_r(\pi^m F^{-1}(a(\bv))) &\geq w_r(a(\bv)) + \epsilon \\
w_r(\pi^{-m} F(b(\bv))) &\geq w_r(b(\bv)) + \epsilon.
\end{align*}
Put
\[
f(\bv) = \pi^{-m} b(\bv) - F^{-1}(a(\bv)).
\]
If $\bv \in M_{(0,r]}$ is such that 
$\bw = F\bv - \pi^m \bv$ also lies in $M_{(0,r]}$, then
\begin{align*}
F(\bv + f(\bw)) - \pi^m (\bv + f(\bw)) &=
F(f(\bw)) - \pi^m f(\bw) + \bw \\
&= F(\pi^{-m} b(\bw)) - a(\bw)
- b(\bw) + \pi^m F^{-1}(a(\bw)) + \bw \\
&= \pi^{-m} F(b(\bw)) + \pi^m F^{-1}(a(\bw))
\end{align*}
lies in $M_{(0,r]}$ as well, and
\begin{align*}
w_r(f(\bw)) &\geq w_r(\pi^{-m}\bw) \\
w_r(F(\bv + f(\bw)) - \pi^m (\bv + f(\bw)))
&\geq w_r(\bw) + \epsilon.
\end{align*}

Now define a sequence $\{\bv_l\}_{l=0}^\infty$ in $M_{(0,r]}$ as follows. Pick
$\overline{x} \in K$ with $v_K(\overline{x}) = d$, and set
\[
\bv_0 = \pi^{-m} [\overline{x}]\be_1 + [\overline{x}^{1/q}]
F^{-1}\be_1.
\]
Given $\bv_l \in M_{(0,r]}$, set
\[
\bw_l = F\bv_l - \pi^m \bv_l, \qquad
\bv_{l+1} = \bv_l + f(\bw_l).
\]
We calculated above that $\bw_l \in M_{(0,r]}$ implies $\bw_{l+1} \in M_{(0,r]}$;
since  $\bw_0 \in M_{(0,r]}$ evidently, we have 
$\bw_l \in M_{(0,r]}$ for all $l$, so 
$\bv_{l+1} \in M_{(0,r]}$ and the iteration continues. Moreover,
\begin{align*}
w_r(\bv_{l+1} - \bv_l) &\geq w_r(\pi^{-m} \bw_l) \\
w_r(\bw_l) &= w_r(F(\bv_{l-1} + f(\bw)) - \pi^m (\bv_{l-1} + f(\bw))) \\
&\geq w_r(\bw_{l-1}) + \epsilon.
\end{align*}
Hence $\bw_l \to 0$ and $f(\bw_l) \to 0$ under $w_r$ as $l \to \infty$, 
so the $\bv_l$ converge to a limit $\bv \in M_{[r,r]}$.

We now check that $\bv \neq 0$.
Since $w_r(\bw_{l+1}) \geq w_r(\bw_l) + \epsilon$ for all $l$,
we certainly have $w_r(\bw_l) \geq w_r(\bw_0)$ for all $l$,
and hence $w_r(\bv_{l+1} - \bv_l) \geq w_r(\pi^{-m} \bw_0)$.
To compute $w_r(\bw_0)$, note that
$w_r(\pi^{-m} [\overline{x}] \be_1) = dr - m$,
whereas 
\begin{align*}
w_r([\overline{x}^{1/q}]F^{-1}\be_1) &\geq
dr/q + c \\
&>dr - m
\end{align*}
by the choice of $d$. Hence $w_r(\bv_0) = dr-m$.
On the other hand,
\begin{align*}
w_r(\pi^{-m} \bw_0) &= 
w_r(\pi^{-m} F\bv_0 - \bv_0) \\
&= w_r(\pi^{-2m} [\overline{x}^q] F\be_1 - [\overline{x}^{1/q}] 
F^{-1}\be_1) \\
&\geq \min\{dqr + c - 2m, dr/q + c\}.
\end{align*}
The second term is strictly greater than $dr-m$ as above, while the
first term equals $dr-m$ plus $dr(q-1)+c-m$, and the latter is 
positive again by the choice of $d$.
Thus $w_r(\pi^{-m} \bw_0) > w_r(\bv_0)$, and so
$w_r(\bv_{l+1} - \bv_l) > w_r(\bv_0)$; in particular,
$w_r(\bv) = w_r(\bv_0)$ and so $\bv \neq 0$.

Since the $\bv_l$ converge to $\bv$ in $M_{[r,r]}$, the
$F\bv_l$ converge to $F\bv$ in $M_{[r/q,r/q]}$. On the other hand,
$F\bv_l = \bw_l + \pi^m \bv_l$, so the $F\bv_l$ converge to $\pi^m \bv$ in
$M_{[r,r]}$. In particular, the $F\bv_l$ form a Cauchy sequence under $w_s$
for $s=r/q$ and $s=r$, hence also for all $s \in [r/q,r]$, and the 
limit in $M_{[r/q,r]}$ must equal both $F\bv$ and $\pi^m \bv$. Therefore $\bv
\in M_{[r/q,r]}$ and $F\bv = \pi^m \bv$ in $M_{[r/q,r/q]}$. But now by
induction on $i$, $\bv \in M_{[r/q^i,r]}$ for all $i$, so $\bv \in 
M_{(0,r]} \subset M$ is a nonzero eigenvector, as desired.
\end{proof}

\begin{lemma} \label{L:into constants}
Suppose that $v_K$ is nontrivial.
Let $M$ be a $\sigma$-module of rank $n$ over $\calR$.
\begin{enumerate}
\item[(a)]
There exists an integer $c_0$ such that for any integer $c \geq c_0$, there
exists an injection $\calR(c)^{\oplus n} \hookrightarrow M$.
\item[(b)]
There exists an integer $c_1$ such that for any integer $c \leq c_1$, there
exists an injection $M \hookrightarrow \calR(c)^{\oplus n}$.
\end{enumerate}
\end{lemma}
\begin{proof}
By taking duals, we may reduce (b) to (a).
We prove (a) by induction on $n$, with empty base case $n=0$.
By Proposition~\ref{P:eigenvector},
there exists an eigenvector of $M$; the saturated span of this
eigenvector is a rank 1 $\sigma$-submodule of $M$, necessarily
isomorphic to some $\calR(m)$ by Proposition~\ref{P:rank 1 class}.
By the induction hypothesis,
we can choose $c_0 \geq m$ so that $\calR(c_0)^{n-1}$ injects into
$M/\calR(m)$. Let $N$ be the preimage of $\calR(c_0)^{n-1}$
in $M$; then there exists an exact sequence
\[
0 \to \calR(m) \to N \to \calR(c_0)^{n-1} \to 0,
\]
which splits by Corollary~\ref{C:hom ext2}.
Thus $N \cong \calR(m) \oplus \calR(c_0)^{n-1} \subseteq M$,
and $\calR(m)$ contains a copy of $\calR(c_0)$ by
Corollary~\ref{C:hom ext2}. This yields the desired result.
\end{proof}

\begin{prop} \label{P:HN}
For any nonzero $\sigma$-module $M$ over $\Gancon$, the slopes of 
all nonzero $\sigma$-submodules of $N$ are bounded below. Moreover,
there is a nonzero $\sigma$-submodule of $N$ of minimal slope,
and any such $\sigma$-submodule is semistable.
\end{prop}
\begin{proof}
To check the first assertion, we may assume (by enlarging $K$ as needed)
that $K$ is algebraically closed, so that $\Gancon = \calR$, and that
$v_K$ is nontrivial.
By Lemma~\ref{L:into constants}, there exists an injection
$M \hookrightarrow \calR(c)^{\oplus n}$ for some $c$, where $n = \rank M$.
By Lemma~\ref{L:semistable}, it follows that $\mu(N) \geq c$ for any
$\sigma$-submodule $N$ of $M$, yielding the first assertion. 

As for the second assertion,
the  slopes of $\sigma$-submodules of $M$ form a discrete
subset of $\QQ$, because their denominators are bounded above by $n$.
Hence this set has a least element, yielding the remaining assertions.
\end{proof}

\begin{prop} \label{P:HN exists}
Every $\sigma$-module over $\Gancon$ admits a Harder-Narasimhan
filtration.
\end{prop}
\begin{proof}
Let $M$ be a nontrivial $\sigma$-module over $\Gancon$.
By Proposition~\ref{P:HN}, the set of 
slopes of nonzero $\sigma$-submodules
of $M$ has a least element $s_1$.
Suppose that $N_1, N_2$ are $\sigma$-submodules of $M$ with
$\mu(N_1) = \mu(N_2) = s_1$; then the internal sum $N_1 + N_2$ is a
quotient of the direct sum $N_1 \oplus N_2$.
By Proposition~\ref{P:HN}, each of $N_1$ and $N_2$ is semistable, as then
is $N_1 \oplus N_2$. Hence
$\mu(N_1 + N_2) \leq s_1$; by the minimality of $s_1$,
we have $\mu(N_1 + N_2) = s_1$. Consequently, the set of
$\sigma$-submodules of $M$ of  slope $s_1$ has a maximal element
$M_1$. Repeating this argument with $M$ replaced by $M/M_1$, and so on,
yields a Harder-Narasimhan filtration.
\end{proof}

\begin{remark} \label{R:no spherical2}
Running this argument is a severe obstacle to working
with spherically complete coefficients (as suggested by
Remark~\ref{R:no spherical0}), as it is no longer clear
that there exists a minimum slope among the $\sigma$-submodules of
a given $\sigma$-module.
\end{remark}

\begin{remark}
It may be possible to simplify the calculations in this section by
using Lemma~\ref{L:frob turnover}; we have not looked thoroughly into this possibility.
\end{remark}

\subsection{More eigenvectors}

We now give a crucial refinement of the conclusion of
Proposition~\ref{P:eigenvector} by extracting an eigenvector of a
specific slope in a key situation. This is essentially
\cite[Proposition~4.15]{me-local};
compare also \cite[Proposition~9.1]{hartl-pink}.
Beware that we are omitting one particularly unpleasant part of the
calculation; see Lemma~\ref{L:pairing} below.

We start by identifying $H^0(\calR(-1))$.
\begin{lemma}
The map
\[
\overline{y} \mapsto \sum_{i\in \ZZ} [\overline{y}^{q^{-i}}] \pi^i
\]
induces a bijection $\gothm_K \to H^0(\calR(-1))$,
where $\gothm_K$ denotes the subset of $K$ on which $v_K$ is positive.
\end{lemma}
\begin{proof}
On one hand, for $v_K(\overline{y}) > 0$,
the sum $y = \sum_{i\in \ZZ} [\overline{y}^{q^{-i}}] \pi^i$
converges and satisfies $y^\sigma = \pi y$.
Conversely, if $y^\sigma = \pi y$, comparing the Teichm\"uller presentations
of $y^\sigma$ and of $\pi y$ forces $y$ to assume the desired form.
\end{proof}

We next give a ``positioning argument'' for elements of $H^1(\calR(m))$,
following \cite[Lemmas~4.13 and~4.14]{me-local}.
\begin{lemma} \label{L:pos h1}
For $m$ a positive integer,
every nonzero element of
$H^1(\calR(m))$ is represented by some $x \in \Galgcon$
with $v_n(x) = v_{m-1}(x)$ for $n \geq m$.
Moreover, we can ensure that for each $n \geq 0$, either
$v_n(x) = \infty$ or $v_n(x) < 0$.
\end{lemma}
\begin{proof}
We first verify that each element of $H^1(\calR(m))$ is represented
by an element of $\Galgcon[\fp]$.
If $v_n(x) \geq 0$ for all $n \in \ZZ$, then the sum
$y = \sum_{i=0}^\infty x^{\sigma^i} \pi^{mi}$ converges in $\calR$
and satisfies $y - \pi^m y^\sigma = x$, so $x$ represents the zero 
class in $H^1(\calR(m))$. In other words, if $x \in \calR$ has
plus-minus-zero representation $x_+ + x_- + x_0$, then $x$ and $x_-$
represent the same class in $H^1(\calR(m))$, and visibly $x_- \in 
\Galgcon[\fp]$.

We next verify that each element of $H^1(\calR(m))$ is represented
by an element $x$ of $\Galgcon[\fp]$ with $\inf_n \{v_n(x)\} > -\infty$.
Given any $x \in \Galgcon$, 
let $x = \sum_i [\overline{x_i}] \pi^i$
be the Teichm\"uller representation of $x$. For each $i$,
let $c_i$ be the smallest nonnegative integer such that
$q^{-c_i} v_K(\overline{x_i}) \geq -1$, and put
\[
y_i = \sum_{j=1}^{c_i} [\overline{x_i}^{q^{-j}}] \pi^{i-mj};
\]
since $c_i$ grows only logarithmically in $i$, $i-mc_i \to \infty$
and the sum $\sum_i y_i$ converges $\pi$-adically.
Moreover, 
\[
\liminf_{i \to \infty}
\min_{1 \leq j \leq c_i}
\{-q^{-j} v_0(\overline{x_i})/(i-mj) \}
\]
is finite, because the same is true for each of
$q^{-j}$ (clear), $-v_0(\overline{x_i})/i$ (by the definition of $\calR$),
and $i/(i-mj)$ (because $c_i$ grows logarithmically in $i$).
Hence the sum $y = \sum_i y_i$ actually converges in $\calR$;
if we set $x' = x - \pi^m y^\sigma + y$, then $x$ and $x'$ represent
the same class in $H^1(\calR(m))$. However,
\[
x' = \sum_i [\overline{x_i}^{q^{-c_i}}] \pi^{i-mc_i}
\]
satisfies $v_n(x') \geq -1$ for all $n$.

Since $x$ and $\pi^m x^\sigma$
represent the same element of $H^1(\calR(m))$, we can also say that
each element of $H^1(\calR(m))$ is represented by an element 
$x$ of $\Galgcon$ with $\inf_n \{v_n(x)\} > -\infty$.
Given such an $x$, put $h = \inf_n \{v_n(x)\}$; we attempt to construct
a sequence $x_0, x_1, \dots$ with the following properties:
\begin{enumerate}
\item[(a)] $x_0 = x$;
\item[(b)] each $x_l$ generates the same element of $H^1(\calR(m))$
as does $x$;
\item[(c)] $x_l \equiv 0 \pmod{\pi^{lm}}$;
\item[(d)] for each $n$ and $l$, $v_n(x_l) \geq h$.
\end{enumerate}
We do this as follows.
Given $x_l$, let $\sum_{i=lm}^\infty [\overline{x_{l,i}}] \pi^i$ be the
Teichm\"uller presentation of $x_l$, and put
\[
u_l = \sum_{i=lm}^{lm+m-1} [\overline{x_{l,i}}] \pi^i.
\]
If $v_{lm+m-1}(x_l) \geq h/q$,
put $x_l - u_l + \pi^m u_l^\sigma$;
otherwise, leave $x_{l+1}$ undefined.

If $x_l$ is defined for each $l$, then set $y = \sum_{l=0}^\infty
u_l$; this sum converges in $\calR$, and its limit
satisfies $x - y + \pi^m y^\sigma = 0$. Hence in this case, $x$ represents
the trivial class in $H^1(\calR(m))$.

On the other hand, if we are able to define $x_l$ but not $x_{l+1}$,
put 
\begin{align*}
y &= x_l - u_l \\
z &= \pi^{-(l+1)m} y^{\sigma^{-l-1}} + \pi^{-lm} u_l^{\sigma^{-l}};
\end{align*}
then $x$ and $z$ represent the same class in $H^1(\calR(m))$. Moreover,
$h/q > v_n(u)$ and 
$v_n(y) \geq h$ for all $n$, so
$hq^{-l-1} > v_n(\pi^{-lm} u_l^{\sigma^{-l}})$
and
$v_n(\pi^{-(l+1)m} y^{\sigma^{-l-1}}) \geq hq^{-l-1}$; consequently
$v_n(z) = v_{m-1}(z)$ for $n \geq m$.
In addition, we can pass from $z$ to the minus part $z_-$
of its plus-minus-zero representation (since $z$ and $z_-$
represent the same class in $H^1(\calR(m))$, as noted above)
to ensure that $v_n(z_-)$ is either
infinite or negative for all $n \geq 0$.
We conclude that every nonzero class in $H^1(\calR(m))$ has a representative
of the desired form.
\end{proof}

\begin{lemma} \label{L:pairing}
Let $d$ be a positive integer.
For any $x \in H^1([d]^* (\calR(d+1)))$, there exists $y \in H^0(\calR(-1))$
nonzero such that $x$ and $y$ pair to zero in $H^1(\calR(-1) \otimes
[d]^* (\calR(d+1))) = H^1([d]^* (\calR(1)))$.
\end{lemma}
\begin{proof}
Identify $H^1([d]^* (\calR(d+1)))$ and $H^1([d]^* (\calR(1)))$ with $H^1(\calR(d+1))$ and $H^1(\calR(1))$, respectively, where the
modules in the latter cases are $\sigma^d$-modules. Then the question
can be stated as follows: for any $x \in \calR$, there exist $y,z \in \calR$
with $y$ nonzero, such that
\[
y^\sigma = \pi y, \qquad xy = z - \pi z^{\sigma^d}.
\]
By Lemma~\ref{L:pos h1}, we may assume that $x \in \Galgcon$ and that
$v_n(x) = v_d(x) < 0$ for 
$n > d$. In this case, the claim follows from a rather involved
calculation \cite[Lemma~4.12]{me-local} which we will not repeat here.
(For a closely related calculation, see \cite[Proposition~9.5]{hartl-pink}.)
\end{proof}

\begin{prop} \label{P:raise Newton}
Assume that $v_K$ is nontrivial.
For $d$ a positive integer, suppose that
\[
 0 \to M_{1,d} \to M \to M_{-1,1} \to 0
\]
is a short exact sequence of $\sigma$-modules over $\calR$. Then $M$
contains an eigenvector of slope $0$.
\end{prop}
\begin{proof}
The short exact sequence corresponds to a class in
\begin{align*}
\Ext(M_{-1,1}, M_{1,d}) &\cong H^1(M_{1,1} \otimes M_{1,d}) \\
&\cong H^1(M_{d+1,d}) \\
&\cong H^1([d]^* (\calR(d+1))).
\end{align*}
From the snake lemma, we obtain an exact sequence
\[
H^0(M) \to H^0(M_{-1,1}) = H^0(\calR(-1)) \to H^1(M_{1,d}) = H^1([d]^* \calR(1)),
\]
in which the second map (the connecting homomorphism)
coincides with the pairing by the given class in
$H^1([d]^*(\calR(d+1)))$. By Lemma~\ref{L:pairing}, this homomorphism
is not injective; hence $H^0(M) \neq 0$, as desired.
\end{proof}
\begin{cor} \label{C:raise Newton1}
Assume that $v_K$ is nontrivial.
For $d$ a positive integer, suppose that
\[
0 \to M_{1,1} \to M \to M_{-1,d} \to 0
\]
is a short exact sequence of $\sigma$-modules over $\calR$. Then $M^\dual$
contains an eigenvector of slope $0$.
\end{cor}
\begin{proof}
Dualize Proposition~\ref{P:raise Newton}.
\end{proof}

\begin{cor} \label{C:raise Newton2}
Assume that $v_K$ is nontrivial.
For $c,c',c''$ integers with $c+c' \leq 2c''$, suppose that
\[
0 \to M_{c,1} \to M \to M_{c',1} \to 0
\]
is a short exact sequence of $\sigma$-modules over $\calR$.
Then $M$ contains an eigenvector of slope $c''$.
\end{cor}
\begin{proof}
By twisting, we may reduce to the case $c'' = 0$.
If $c \leq 0$, then already $M_{c,1}$ contains an eigenvector of
slope 0 by Corollary~\ref{C:hom ext2}, so we may assume
$c \geq 0$. Since $c+c' \leq 0$, by Corollary~\ref{C:hom ext2}, we can
find a copy of $M_{-c,1}$ within $M_{c',1}$; taking the preimage of
$M_{-c,1}$ within $M$ allows us to reduce to the case $c' = -c$.

We treat the cases $c \geq 0$
by induction, with base case $c=0$ already treated.
If $c > 0$, then twisting yields an exact sequence
\[
0 \to M_{c-1,1} \to M(-1) \to M_{-c-1,1} \to 0.
\]
By Corollary~\ref{C:hom ext2}, we can choose a submodule of $M_{-c-1,1}$
isomorphic to $M_{-c+1,1}$; let $N$ be its inverse image in $M(-1)$.
Applying the induction hypothesis to the sequence
\[
0 \to M_{c-1,1} \to N \to M_{-c+1,1} \to 0
\]
yields an eigenvector of $N$ of slope 0, and hence an eigenvector of $M$
of slope 1. Let $P$ be the saturated span of that eigenvector;
it is isomorphic to $M_{m,1}$ for some $m$ by
Proposition~\ref{P:rank 1 class}, and we must have
$m \leq 1$ by Corollary~\ref{C:hom ext2}.
If $m \leq 0$, $M$ has an eigenvector of slope 0,
so suppose instead that $m = 1$.
We then have an exact sequence
\[
0 \to P \cong M_{1,1} \to M \to M/P \to 0
\]
in which $M/P$, which has rank 1 and degree $-1$ (by Lemma~\ref{L:degree adds}),
is isomorphic to $M_{-1,1}$ by 
Proposition~\ref{P:rank 1 class}.
Applying Proposition~\ref{P:raise Newton} now yields the desired result.
(Compare \cite[Corollary~9.2]{hartl-pink}.)
\end{proof}

\subsection{Existence of standard submodules}

We now run the induction setup of \cite[Theorem~11.1]{hartl-pink}
to produce standard submodules of a $\sigma$-module of small slope.

\begin{defn} \label{D:induct}
For a given integer $n \geq 1$, let (A$_n$), (B$_n$) denote the following 
statements about $n$.
\begin{enumerate}
\item[(A$_n$)]
Let $a$ be any positive integer, and let $M$ be any $\sigma^a$-module over 
$\calR$. If $\rank(M)\leq n$ and $\deg(M) \leq 0$,
then $M$ contains an eigenvector of slope 0.
\item[(B$_n$)]
Let $a$ be any positive integer, and let $M$ be any $\sigma^a$-module over 
$\calR$. 
If $\rank(M) \leq n$, then $M$ contains a saturated $\sigma^a$-submodule
which is standard of  slope $\leq \mu(M)$.
\end{enumerate}
Note that if $v_K$ is nontrivial, both
(A$_1$) and (B$_1$) hold thanks to
Proposition~\ref{P:rank 1 class}.
\end{defn}

\begin{lemma} \label{L:induct1}
Assume that $v_K$ is nontrivial.
For $n \geq 2$, if (A$_{n-1}$) and (B$_{n-1}$)
hold, then (A$_n$) holds.
\end{lemma}
\begin{proof}
It suffices to show that if $M$ is a $\sigma^a$-module over 
$\calR$ with $\rank(M) = n$ and $\deg(M) \leq 0$, then
$M$ contains an eigenvector of slope 0; after twisting, we
may reduce to the case $1-n \leq \deg(M) \leq 0$. Suppose on the contrary
that no such eigenvector exists; by Corollary~\ref{C:hom ext2},
$M$ then contains no eigenvector of any nonpositive slope.

On the other hand, by Proposition~\ref{P:eigenvector},
$M$ contains an eigenvector; in particular, $M$ contains a saturated
$\sigma^a$-submodule of rank 1. By Proposition~\ref{P:rank 1 class},
we thus have an exact sequence
\[
0 \to M_{c,1} \to M \to N \to 0
\]
for some integer $c$, which by hypothesis must be positive. 

Choose $c$ as small as possible; then
$\deg(N) = \deg(M) - c \leq 0$ by
Lemma~\ref{L:degree adds}, so by (A$_{n-1}$), $N$ contains an
eigenvector of slope 0. That is, $N$ contains a $\sigma^a$-submodule
isomorphic to $M_{0,1}$; let $M'$ be the preimage in $M$
of that submodule. We then have an exact sequence
\[
0 \to M_{c,1} \to M' \to M_{0,1} \to 0.
\]
By Corollary~\ref{C:raise Newton2},
$M'$ contains an eigenvector of slope $\lceil c/2 \rceil$.
By the minimality of the choice of $c$, we must have
$c \leq \lceil c/2 \rceil$, or $c = 1$.

Put $c' = \deg(M) - 1 = \deg(N)$, so that $c' < 0$.
By (B$_{n-1}$), $N$ contains a saturated
$\sigma^a$-submodule $P$ of average slope
$\leq c'/(n-1)$; let $P'$ be the preimage in $M$ of that
$\sigma^a$-submodule.
If $\rank(P) < n-1$, then 
$\deg(P') \leq 1 + \rank(P)(c'/(n-1)) < 1$; since
$\deg(P')$ is an integer, we have
$\deg(P') \leq 0$. By (A$_{n-1}$), $P'$
contains an eigenvector of slope 0, contradicting our hypothesis.

If $\rank(P) = n-1$, 
we have an exact sequence
\[
0 \to M_{1,1} \to P' \to P \cong M_{c'',n-1} \to 0
\]
for some $c'' \leq -1$. 
By Corollary~\ref{C:hom ext2}, $M_{c'',n-1}$ contains a copy of 
$M_{-1,n-1}$; choose such a copy and let $P''$ be its inverse image in $P'$.
By Corollary~\ref{C:raise Newton1}, $(P'')^\dual$ contains an eigenvector of
slope $0$, and hence a primitive eigenvector of slope at most $0$; this
eigenvector corresponds to a rank $n-1$ submodule of $P''$ of slope at most 
$0$. By (A$_{n-1}$), $P''$ contains an eigenvector of slope $0$, contradicting
our hypothesis.

In any case, our hypothesis that $M$ contains no eigenvector of
slope 0 has been contradicted, yielding the desired result.
\end{proof}

\begin{lemma} \label{L:induct2}
Assume that $v_K$ is nontrivial.
For $n \geq 2$, if (A$_{n}$) and (B$_{n-1}$)
hold, then (B$_n$) holds.
\end{lemma}
\begin{proof}
Let $M$ be a $\sigma^a$-module of rank $n$ and degree $c$.
Put $b = n/\gcd(n,c)$; then by
(A$_n$) applied after twisting, $[b]_* M$ (which has rank $n$
and degree $bc$, by Proposition~\ref{P:pull basic}) has an eigenvector of slope $bc/n$. That is,
$M$ has a $b$-eigenvector $\bv$ of slope $c/n$; this gives a
nontrivial map $f: M_{bc/n,b} \to M$ sending a standard basis to
$\bv, F\bv, \dots, F^{b-1} \bv$.
Let $N$ be the saturated span of the image of $f$, and put
$m = \rank N$. Then $\wedge^m N$ admits a $b$-eigenvector of
slope $cm/n$, so by Corollary~\ref{C:hom ext2}, the average slope of $N$ is at most $c/n$.
If $m < n$, we may apply (B$_{n-1}$) to $N$ to obtain
the desired result.

Suppose instead that $m = n$, which also implies $b=n$ since
necessarily $m \leq b \leq n$.
Then the map $f$ is injective, so its image has average slope
$c/n$. By Lemma~\ref{L:same rank}, $f$ must in fact be surjective;
thus $M \cong M_{c,n}$, as desired.
\end{proof}

\subsection{Dieudonn\'e-Manin decompositions}

\begin{defn}
A \emph{Dieudonn\'e-Manin decomposition} of a $\sigma$-module $M$ over
$\calR$ is a direct sum decomposition $M = \oplus_{i=1}^m M_{c_i,d_i}$
of $M$ into standard $\sigma$-submodules.
The \emph{slope multiset} of such a decomposition is the union
of the multisets consisting of $c_i/d_i$ with multiplicity $d_i$ for
$i=1, \dots, m$.
\end{defn}
\begin{remark}
If $M$ admits a Dieudonn\'e-Manin decomposition, then
$M$ admits a basis of $n$-eigenvectors for $n = (\rank M)!$;
more precisely, any basis of $H^0([n]_* M)$ over the fixed field of
$\sigma^n$ gives a basis of $[n]_* M$ over $\calR$. The slopes of these
$n$-eigenvectors coincide with the slope multiset of the decomposition.
\end{remark}

\begin{prop} \label{P:Dieudonne-Manin}
Assume that $v_K$ is nontrivial. Then every $\sigma$-module $M$
over $\calR$ admits a Dieudonn\'e-Manin decomposition.
\end{prop}
\begin{proof}
We first show that every semistable $\sigma$-module $M$ over $\calR$
is isomorphic to a direct sum of standard $\sigma$-submodules of slope
$\mu(M)$. We see this by induction on $\rank(M)$; by
Lemmas~\ref{L:induct1} and~\ref{L:induct2}, we have (B$_n$) for
all $n$, so $M$ contains a saturated $\sigma$-submodule $N$ which
is standard of some slope $\leq \mu(M)$. Since $M$ has been assumed
semistable, we have $\mu(N) \geq \mu(M)$; hence 
$\mu(N) = \mu(M)$, and $M/N$ is also semistable. By the induction
hypothesis, $M/N$ splits as a direct sum of standard $\sigma$-submodules
of slope $\mu(M)$; then by Corollary~\ref{C:hom ext2}, the exact sequence
\[
0 \to N \to M \to M/N \to 0
\]
splits. This yields the desired result.

In the general case, by Proposition~\ref{P:HN exists},
$M$ has an HN-filtration $0 = M_0 \subset M_1 \subset \cdots \subset M_l = M$,
with each successive quotient $M_i/M_{i-1}$ semistable of slope $s_i$,
and $s_1 < \cdots < s_l$. By the above, $M_i/M_{i-1}$ admits
a Dieudonn\'e-Manin decomposition with all slopes $s_i$; 
the filtration then splits
thanks to Corollary~\ref{C:hom ext2}. Hence $M$ admits
a Dieudonn\'e-Manin decomposition.
\end{proof}

We will use the case of $v_K$ nontrivial to establish the existence of
Dieudonn\'e-Manin decompositions also when $v_K$ is trivial.
\begin{defn}
For $R$ a (commutative) ring, let
$R((t^{\QQ}))$ denote the Hahn-Mal'cev-Neumann algebra of generalized
power series $\sum_{i \in \QQ} c_i t^i$, where each $c_i\in R$
and the set of $i \in \QQ$ with $c_i \neq 0$ is well-ordered (has no
infinite decreasing subsequence); these series form a ring under 
formal series multiplication, with a natural valuation $v$ given
by $v(\sum_i c_i t^i) = \min\{i: c_i \neq 0\}$.
For $R$ an algebraically closed
field, $R((t^{\QQ}))$ is also algebraically closed;
see \cite[Chapter~13]{passman} for this and other properties of these
algebras.
\end{defn}

\begin{lemma} \label{L:gen power}
Suppose that $k$ is algebraically closed and that
$K = k((t^{\QQ}))$.
Let $M$ be a $\sigma$-module over $\calO[\fp]$ such that
$M \otimes \calR$ admits a basis of eigenvectors. Then any such basis
is a basis of $M$.
\end{lemma}
\begin{proof}
We may identify $\Gamma^K$ with the $\pi$-adic completion of 
$\calO((t^{\QQ}))$. In so doing, elements of $\calR$ can be viewed
as formal sums $\sum_{i \in \QQ} c_i t^i$ with $c_i \in \calO[\fp]$.

Suppose $\bv \in M \otimes \calR$ 
nonzero satisfies $F\bv = \pi^m \bv$. We can then formally write
$\bv = \sum_{i \in \QQ} \bv_i t^i$ with $\bv_i \in M$,
and then we have $F\bv_i = \pi^m \bv_{qi}$ for each $i$.
If $\bv_i \neq 0$ for some $i < 0$, we then have
$\bv_{q^l i} = \pi^{-lm} F^l \bv_i$, but this violates the
convergence condition defining $\calR$. Hence $\bv_i = 0$ for $i < 0$.

Let $\calR^+$ be the subring of $\calR$ consisting of series
$\sum_i c_i t^i$ with $c_i = 0$ for $i<0$. Now if
$M \otimes \calR$ admits a basis of eigenvectors, then we have just
shown that each 
basis element belongs to $M \otimes \calR^+$, and likewise for the dual
basis of $M^\dual \otimes \calR^+$. 
We can then reduce modulo the ideal of $\calR^+$ consisting
of series with constant coefficient zero, to produce a basis of eigenvectors
of $M$.
\end{proof}

\begin{remark}
Beware that in the proof of Lemma~\ref{L:gen power}, 
there do exist nonzero eigenvectors
in $M \otimes \calR^+$ with constant coefficient zero; however, these
eigenvectors cannot be part of a basis.
\end{remark}

\begin{theorem} \label{T:Dieudonne-Manin}
Let $M$ be a $\sigma$-module over $\calR$.
\begin{enumerate}
\item[(a)]
There exists a Dieudonn\'e-Manin decomposition of $M$.
\item[(b)]
For any Dieudonn\'e-Manin decomposition $M = \oplus_{j=1}^m M_{c_j,d_j}$
of $M$, let $s_1 < \cdots < s_l$ be the distinct elements of the slope
multiset of the decomposition. For $i=1, \dots, l$, let
$M_i$ be the direct sum of $M_{c_j,d_j}$ over all $j$ for which $c_j/d_j
\leq s_i$. Then the filtration $0 \subset M_1 \subset \cdots \subset M_l = M$
coincides with the HN-filtration of $M$.
\item[(c)]
The slope multiset of any Dieudonn\'e-Manin decomposition of $M$
consists of the HN-slopes of $M$. In particular, the slope multiset does not depend on the choice of the decomposition.
\end{enumerate}
\end{theorem}
\begin{proof}
\begin{enumerate}
\item[(a)]
For $v_K$ nontrivial, this is Proposition~\ref{P:Dieudonne-Manin},
so we need only treat the case of $v_K$ trivial. Another way to
say this is every $\sigma$-module $M$ over $\calO[\fp]$ is isomorphic
to a direct sum of standard $\sigma$-submodules.

Set notation as in Lemma~\ref{L:gen power}. By
Proposition~\ref{P:Dieudonne-Manin}, $M \otimes \calR$ is
isomorphic to a direct sum of standard $\sigma$-modules. That
direct sum has the form $N \otimes \calR$, where $N$ is the corresponding
direct sum of standard $\sigma$-modules over $\calO[\fp]$. The isomorphism
$M \otimes \calR \to N \otimes \calR$ corresponds to an element of
$H^0((M^\dual \otimes N) \otimes \calR)$, which extends to a basis of
eigenvectors of
$[n]_* (M^\dual \otimes N) \otimes \calR$ for some $n$.
By Lemma~\ref{L:gen power}, this basis consists of elements of
$[n]_* (M^\dual \otimes N)$; hence $M \cong N$, as desired.

\item[(b)]
By Lemma~\ref{L:push semi}, any standard $\sigma$-module is stable;
hence a direct sum of standard $\sigma$-modules of a single slope is semistable.
Thus the described filtration is indeed an HN-filtration.

\item[(c)]
This follows from (b).
\end{enumerate}
\end{proof}

\begin{remark}
The case of $v_K$ trivial in Theorem~\ref{T:Dieudonne-Manin}(a)
is precisely the standard
Dieudonn\'e-Manin classification of $\sigma$-modules over a complete
discretely valued field with algebraically closed residue field.
It is more commonly derived on its own, as in
\cite{dieudonne}, \cite{manin}, \cite{katz}, or
\cite[Theorem~5.6]{me-local}.
\end{remark}

\begin{cor} \label{C:irreducible}
For any coprime integers $c,d$ with $d>0$,
$\End(M_{c,d})$ is a division algebra.
\end{cor}
\begin{proof}
Suppose $\phi \in \End(M_{c,d})$ is nonzero.
Decompose $\image(\phi)$ according to Theorem~\ref{T:Dieudonne-Manin}; 
then each
standard summand of $\image(\phi)$ must have slope $\leq c/d$ by
Corollary~\ref{C:hom ext2}. On the other hand, each summand is a 
$\sigma$-submodule of
$M_{c,d}$, so must have slope $\geq c/d$ again by Corollary~\ref{C:hom ext2}.
Thus each standard summand of $\image(\phi)$ must have slope exactly $c/d$.
In particular, there can be only one such summand, it must have rank $d$, and
by Lemma~\ref{L:same rank}, $\image(\phi) = M_{c,d}$. Hence
$\phi$ is surjective; since $\phi$
is a linear map between free modules of the same
finite rank, it is also injective. We conclude that $\End(M_{c,d})$ is indeed
a division algebra, as desired.
\end{proof}

\begin{prop} \label{P:which stable}
A $\sigma$-module $M$ over $\calR$ is semistable (resp.\ stable)
if and only if $M \cong M_{c,d}^{\oplus n}$ for some $c,d,n$
(resp.\ $M \cong M_{c,d}$ for some $c,d$).
\end{prop}
\begin{proof}
This is an immediate corollary of Theorem~\ref{T:Dieudonne-Manin}.
(Compare \cite[Corollary~11.6]{hartl-pink}.)
\end{proof}

\begin{remark}
By Theorem~\ref{T:Dieudonne-Manin}, every $\sigma$-module over 
$\calR$ decomposes as a direct sum of semistable $\sigma$-modules,
i.e., the Harder-Narasimhan filtration splits. However,
when $v_K$ is nontrivial, 
this decomposition/splitting is not canonical, so it does not make
sense to try to prove any descent results for such decompositions.
(When $v_K$ is trivial, the splitting is unique by virtue of
Corollary~\ref{C:hom ext2}.)
Of course, the number and type of summands in a Dieudonn\'e-Manin
decomposition are unique, since they are determined by the HN-polygon; 
indeed, they constitute complete invariants
for isomorphism of $\sigma$-modules over $\calR$ 
(compare \cite[Corollary~11.8]{hartl-pink}).
\end{remark}

\begin{prop} \label{P:wedge first step}
Let $M$ be a $\sigma$-module over $\calR$,
let $M_1$ be the first step in the Harder-Narasimhan filtration,
and put $d = \rank(M_1)$. Then $\wedge^d M_1$ 
is the first step in the Harder-Narasimhan
filtration of $\wedge^d M$.
\end{prop}
\begin{proof}
Decompose $M$ according to Theorem~\ref{T:Dieudonne-Manin}, so that 
$M_1$ is the direct
sum of the summands of minimum slope $s_1$. Take the $d$-th exterior power 
of this 
decomposition (i.e., apply the K\"unneth formula); 
by Lemma~\ref{L:standard arith},
the minimum slope among the new summands is $d s_1$, achieved only by 
$\wedge^d M_1$.
\end{proof}
\begin{remark}
More generally, the first step in the Harder-Narasimhan filtration of 
$\wedge^i M$
is $\wedge^i M_j$, for the smallest $j$ such that $\rank(M_j) \geq i$; the
argument is similar.
\end{remark}

\begin{prop} \label{P:which homs}
Let $M$ be a $\sigma$-module over $\calR$, and let $M \cong
\oplus_{i=1}^l M_{c_i,d_i}$ be a Dieudonn\'e-Manin decomposition of $M$.
\begin{enumerate}
\item[(a)]
If $v_K$ is nontrivial, then 
there exists a nonzero homomorphism
$f: M_{c,d} \to M$  of $\sigma$-modules if and only if
$c/d \geq \min_i \{c_i/d_i\}$, and there exists a
nonzero homomorphism
$f: M \to M_{c,d}$  of $\sigma$-modules if and only if
$c/d \leq \max_i \{c_i/d_i\}$.
\item[(b)]
If $v_K$ is trivial, then 
there exists a nonzero homomorphism
$f: M_{c,d} \to M$ or $f: M \to M_{c,d}$ of $\sigma$-modules if and only if
$c/d \in \{c_1/d_1, \dots, c_l/d_l\}$.
\end{enumerate}
\end{prop}
\begin{proof}
Apply Corollary~\ref{C:hom ext2}.
\end{proof}

\subsection{The calculus of slopes}

Theorem~\ref{T:Dieudonne-Manin} affords a number of consequences for the calculus of slopes.
\begin{defn}
Let $M$ be a $\sigma$-module over $\Gancon$. Define
the \emph{absolute HN-slopes} and \emph{absolute HN-polygon}
of $M$ to be the HN-slopes and HN-polygon of $M \otimes \calR$,
and denote the latter by $P(M)$.
We say $M$ is \emph{pure} (or \emph{isoclinic})
of slope $s$ if the absolute HN-slopes of $M$ are all equal to $s$.
By Proposition~\ref{P:which stable}, $M$ is isoclinic if and only if 
$M \otimes \calR$ is semistable. We use the adjective \emph{unit-root} to mean
``isoclinic of slope 0''.
\end{defn}
\begin{remark}
We will show later (Theorem~\ref{T:descend special}) that the HN-filtration of $M \otimes \calR$ coincides with the base extension of the HN-filtration of $M$, which will mean that the absolute HN-slopes of $M$ coincide with the HN-slopes of $M$.
\end{remark}
\begin{prop} \label{P:slope arith}
Let $M$ and $M'$ be $\sigma$-modules over $\Gancon$.
Let $c_1, \dots, c_m$ and $c'_1, \dots, c'_n$ be the absolute HN-slopes of 
$M$ and $M'$, respectively.
\begin{enumerate}
\item[(a)]
The absolute HN-slopes of $M \oplus M'$ are $c_1, \dots, c_m, c'_1, \dots, c'_n$.
\item[(b)]
The absolute HN-slopes of $M \otimes M'$ are $c_i c'_j$ for $i=1, \dots, m$, 
$j=1,\dots, n$.
\item[(c)]
The absolute HN-slopes of $\wedge^d M$ are $c_{i_1} + \cdots + c_{i_d}$ for all
$1 \leq i_1 < \cdots < i_d \leq m$.
\item[(d)]
The absolute HN-slopes of $[a]_* M$ are $ac_1, \dots, ac_m$.
\item[(e)]
The absolute HN-slopes of $M(b)$ are $c_1 + b, \dots, c_m + b$.
\end{enumerate}
\end{prop}
\begin{proof}
There is no harm in tensoring up to $\calR$, or in
applying $[a]_*$ for some positive
integer $a$. 
In particular, using Theorem~\ref{T:Dieudonne-Manin},
we may reduce to the case where $M$ and $M'$ admit bases of eigenvectors,
whose slopes must be the $c_i$ and the $c'_j$. Then we obtain
bases of eigenvectors of $M \oplus M'$, $M \otimes M'$, $\wedge^d M$,
$[a]_* M$, $M(b)$, and thus may read off the claims.
\end{proof}

\begin{prop} \label{P:hom ext3}
Let $M_1, M_2$ be $\sigma$-modules over $\Gancon$ such that
each absolute HN-slope of $M_1$ is less than each absolute HN-slope of $M_2$. Then $\Hom(M_1, M_2) = 0$.
\end{prop}
\begin{proof}
Tensor up to $\calR$,
then apply Theorem~\ref{T:Dieudonne-Manin} and 
Proposition~\ref{P:which homs}.
\end{proof}

\subsection{Splitting exact sequences}

As we have seen already (e.g., in Proposition~\ref{P:raise Newton}), a 
short exact
sequence of $\sigma$-modules over $\Gancon$ 
may or may not split. Whether 
or not
it splits depends very much on the Newton polygons involved. For starters, we
have the following.

\begin{defn}
Given Newton polygons $P_1, \dots, P_m$, define the \emph{sum} $P_1 + 
\cdots + P_m$ of these
polygons to be the Newton polygon whose slope multiset is the union of the
slope multisets of $P_1, \dots, P_m$. Also, write $P_1 \geq P_2$ 
to mean that $P_1$
lies above $P_2$.
\end{defn}

\begin{prop} \label{P:split exact}
Let $0 \to M_1 \to M \to M_2 \to 0$ be an exact sequence of $\sigma$-modules
over $\Gancon$. 
\begin{enumerate}
\item[(a)]
We have $P(M) \geq P(M_1) + P(M_2)$.
\item[(b)]
We have $P(M) = P(M_1) + P(M_2)$ if and only if
the exact sequence splits over $\calR$.
\end{enumerate}
\end{prop}
\begin{proof}
For (a), note that from the HN-filtrations
 of $M_1 \otimes \calR$ and $M_2 \otimes \calR$, we obtain a semistable
filtration of $M \otimes \calR$ whose
Newton polygon is $P(M_1) + P(M_2)$. The claim
now follows from Proposition~\ref{P:filtrations}.

For (b), note that if the sequence splits, then $P(M) = P(M_1) + P(M_2)$ by 
Proposition~\ref{P:slope arith}. Conversely, suppose that $P(M) = P(M_1) + 
P(M_2)$; we prove by induction on rank that if $\Gancon = \calR$,
then the exact sequence splits.
Our base case is where
$M_1$ and $M_2$ are standard. If $\mu(M_1) \leq \mu(M_2)$, then the exact 
sequence
splits by Corollary~\ref{C:hom ext2}, so assume that $\mu(M_1) > \mu(M_2)$.
By Theorem~\ref{T:Dieudonne-Manin}, we have an
isomorphism $M \cong M_1 \oplus M_2$, which gives a map $M_2 \to M$. 
By Corollary~\ref{C:irreducible}, the composition $M_2 \to M \to M_2$ is 
either
zero or an isomorphism; in the former case, by exactness the image of 
$M_2 \to M$
must land in $M_1 \subseteq M$. But that violates 
Corollary~\ref{C:hom ext2}, 
so the composition $M_2 \to M \to M_2$ is an isomorphism, and the exact 
sequence
splits.

We next treat the case of $M_1$ nonstandard. 
Apply Theorem~\ref{T:Dieudonne-Manin}
to obtain a decomposition $M_1 \cong N \oplus N'$ with $N$ standard. We have
\begin{align*}
P(M) &\geq P(N) + P(M/N) \qquad \mbox{[by (a)]} \\
&\geq P(N) + P(M_1/N) + P(M_2) \qquad \mbox{[by (a)]} \\
&= P(M_1) + P(M_2) \qquad \mbox{[because $N$ is a summand of $M_1$]} \\
&= P(M) \qquad \mbox{[by hypothesis]}.
\end{align*}
Hence all of the inequalities must be equalities; in particular,
$P(M/N) = P(M_1/N) + P(M_2)$. By the induction hypothesis, the exact sequence
$0 \to M_1/N \to M/N \to M_2 \to 0$ splits; consequently, the exact sequence
$0 \to N' \to M \to M/N' \to 0$ splits. But we have an exact sequence
$0 \to N \to M/N' \to M_2 \to 0$ and as above, we have $P(M/N') = P(N) + P(M_2)$, so
this sequence also splits by the induction hypothesis. This yields the claim.

To conclude, note that the case of $M_2$ nonstandard follows from the case of $M_1$
nonstandard by taking duals. Hence we have covered all cases.
\end{proof}

\begin{cor}
Let $0 \to M_1 \to M \to M_2 \to 0$ be an exact sequence of
$\sigma$-modules over $\calR$, such that every slope 
of $M_1$ is less than or equal to every slope of 
$M_2$. Then the exact sequence splits; in particular,
the HN-multiset of $M$ is the union of the HN-multiset
of $M_1$ and $M_2$.
\end{cor}
\begin{proof}
With the assumption on the slopes, the filtration induced by the HN-filtrations
of $M_1$ and $M_2$ becomes an HN-filtration
after possibly removing one redundant step in the middle (in case
the highest slope of $M_1$ coincides with the lowest slope of $M_2$).
Thus its Newton polygon coincides with the HN-polygon, so
Proposition~\ref{P:split exact} yields the claim.
\end{proof}
\begin{cor} \label{C:filt splits}
Let $0 = M_0 \subset M_1 \subset \cdots \subset M_l = M$ be a
filtration of a $\sigma$-module $M$ over $\calR$ by saturated
$\sigma$-submodules with isoclinic quotients. Suppose that the Newton polygon of the filtration coincides with the HN-polygon of $M$. Then
the filtration splits.
\end{cor}

\begin{remark}
In certain contexts, one can obtain stronger splitting theorems;
for instance, the key step in \cite{me-full} is a splitting theorem
for $\sigma$-modules with connection over $\Gcon$
(in the notation of Section~\ref{subsec:robba}).
\end{remark}

\section{Generic and special slope filtrations}
\label{sec:generic special}

Given a $\sigma$-module over $\Gcon[\fp]$, we have two paradigms
for constructing slopes and HN-polygons: the ``generic'' paradigm, in which
we pass to $\Gamma[\fp]$ as if $v_K$ were trivial, and the ``special'' 
paradigm, in which we pass to $\Gancon$. (See Section~\ref{subsec:gen spec2} for an
explanation of the use of these adjectives.)
In this chapter, we compare these
paradigms: our main results are that the special HN-polygon lies above the
generic one (Proposition~\ref{P:comparison}), and that when the two polygons
coincide, one obtains a common HN-filtration over $\Galgcon[\fp]$
(Theorem~\ref{T:coincide}).
This last result is a key tool for constructing slope filtrations in 
general.

\setcounter{equation}{0}
\begin{convention}
We continue to retain notations as in Chapter~\ref{sec:basic rings}.
We again point out that when working over $\Gcon[\fp]$, the adjective
``generic'' will imply passage to $\Gamma[\fp]$, while the adjective
``special'' will imply passage to $\Gancon$. We also abbreviate such expressions as ``generic absolute HN-slopes'' to ``generic HN-slopes''. (Keep in mind that
the modifier ``absolute'' will ultimately be rendered superfluously
anyway by Theorem~\ref{T:descend special}.)
\end{convention}

\subsection{Interlude: lattices}

Besides descending subobjects, we will also have need to descend entire
$\sigma$-modules; this matter is naturally discussed in terms of lattices.

\begin{defn}
Let $R \hookrightarrow S$ be an injection of domains, and let $M$
be a finite locally free $S$-module. An \emph{$R$-lattice} in $M$ is an
$R$-submodule $N$ of $M$ such that the induced map $N \otimes_R S \to M$
is a bijection. If $M$ is a $\sigma$-module, an $R$-lattice in the category
of $\sigma$-modules is a module-theoretic $R$-lattice which is stable under
$F$.
\end{defn}

The existence of a $\Gamma$-lattice for a $\sigma$-module defined over
$\Gamma[\fp]$ is closely tied to nonnegativity of the slopes.
\begin{prop} \label{P:lattice}
Let $M$ be a $\sigma$-module over $\Gcon[\fp]$ with
nonnegative generic slopes. Then $M$ contains an $F$-stable
 $\Gcon$-lattice $N$.
 Moreover, if the generic slopes of $M$ are all
zero, then $N$ can be chosen so that $F: \sigma^* M \to M$ is an 
isomorphism.
\end{prop}
\begin{proof}
Put $M' = M \otimes \Galg[\fp]$; then 
by Theorem~\ref{T:Dieudonne-Manin},
we can write $M'$ as a direct sum
of standard submodules, whose slopes by hypothesis are nonnegative.
From this presentation, we immediately obtain a 
$\Galg$-lattice of $M'$ (generated by standard basis vectors of
the standard submodules);
its intersection with $M$ gives the desired lattice.
\end{proof}

Proposition~\ref{P:lattice} also has the following converse. 
\begin{prop} \label{P:nonnegative}
Let $M$ be a $\sigma$-module over $\Gamma$. Then the generic
HN-slopes of $M$ are all nonnegative; moreover, they are all zero
if and only if $F: \sigma^* M \to M$ is an
isomorphism.
\end{prop}
\begin{proof}
Let $\be_1, \dots, \be_n$ be a basis of $M$, and
define the $n \times n$ matrix $A$ over $\Gamma$
by $F\be_j = \sum_i A_{ij} \be_i$. 
Suppose $\bv$ is an eigenvector of $M$,
and write $\bv = \sum_i c_i \be_i$
and $F\bv = \sum_i d_i \be_i$; then $\min_i \{w(d_i)\} - \min_i\{w(c_i)\}$
is the slope of $\bv$. But $w(A_{ij}) \geq 0$ for all $i,j$, so
$\min_i\{w(d_i)\} \geq \min_i\{w(c_i)\}$. This yields the first claim.

For the second claim, note on one hand that if
$F: \sigma^* M \to M$ is an isomorphism,
then $M^\dual$ is also a $\sigma$-module over $\Gamma$, and so
the generic HN-slopes of both $M$ and $M^\dual$ are nonnegative.
Since these slopes are negatives of each other
by Proposition~\ref{P:slope arith}, they must all be zero.
On the other hand, if the generic HN-slopes of $M$ are all zero, 
then by Theorem~\ref{T:Dieudonne-Manin},
$M \otimes \Galg[\fp]$ admits a basis 
$\be_1, \dots, \be_n$ of eigenvectors of slope 0.
Put $M' = M \otimes \Galg$; for $i=0, \dots, n$, let $M'_i$ be the intersection
of $M'$ with the $\Galg[\fp]$-span of $\be_1, \dots, \be_i$.
Then each $M'_i$ is $F$-stable;
moreover, $M'_i/M'_{i-1}$ is spanned by the image of $\pi^{c_i} \be_i$
for some $c_i$, and so $F: \sigma^*(M'_i/M'_{i-1}) \to (M'_i/M'_{i-1})$
is an isomorphism for each $i$. It follows that $F: \sigma^* M' \to M'$
is an isomorphism, as then is $F: \sigma^* M \to M$.
\end{proof}

\begin{remark}
The results in this section can also be proved using cyclic vectors, as in
\cite[Proposition~5.8]{me-local}; compare
Lemma~\ref{L:use cyclic vector} below.
\end{remark}

\subsection{The generic HN-filtration}

Since the distinction between $v_K$ trivial and nontrivial was not
pronounced in the previous chapter, it is worth taking time out to
clarify some phenomena specific to the ``generic'' ($v_K$ trivial)
setting.

\begin{prop} \label{P:DM decomp}
For any $\sigma$-module $M$ over $\Galg[\fp]$, there is a unique decomposition
$M = P_1 \oplus \cdots \oplus P_l$, where each $P_i$ is isoclinic, and the
generic slopes $\mu(P_1), \dots, \mu(P_l)$ are all distinct.
\end{prop}
\begin{proof}
The existence of such a decomposition follows from
Theorem~\ref{T:Dieudonne-Manin}; the uniqueness follows from
repeated application of Corollary~\ref{C:hom ext2}.
\end{proof}

\begin{defn}
Let $M$ be a $\sigma$-module over $\Galg[\fp]$.
Define the \emph{slope decomposition} of $M$ to be the decomposition
$M = P_1 \oplus \cdots \oplus P_l$ given by Proposition~\ref{P:DM decomp}.
\end{defn}

For the rest of this section, we catalog some routine 
methods for identifying
the generic slopes of a $\sigma$-module.
\begin{defn}
Let $M$ be a $\sigma$-module over $\Gamma[\fp]$ of rank $n$.
A \emph{cyclic vector} of $M$ is an element $\bv \in M$ such that
$\bv, F\bv, \cdots, F^{n-1}\bv$ form a basis of $M$.
\end{defn}

\begin{lemma} \label{L:use cyclic vector}
Let $M$ be a $\sigma$-module over $\Gamma[\fp]$ of rank $n$.
Let $\bv$ be a cyclic vector of $M$, and define $a_0, \dots, 
a_{n-1} \in \Gamma[\fp]$ by the equation
\[
F^n \bv + a_{n-1} F^{n-1} \bv + \cdots + a_0 \bv = 0.
\]
Then the generic HN-polygon of $M$ coincides with the Newton polygon
of the polynomial $x^n + a_{n-1} x^{n-1} + \cdots + a_0$.
\end{lemma}
\begin{proof}
There is no harm in assuming that $K$ is algebraically closed, that $\pi$
is fixed by $\sigma$, or that
$\calO$ is large enough that the slopes of $M$ are all integers.
Then $M$ admits a basis $\be_1, \dots, \be_n$ with
$F\be_i = \pi^{c_i} \be_i$ for some integers $c_i$.

Put $\bv_0 = \bv$. Given $\bv_l$, write
$\bv_l = x_{l,1}\be_1 + \cdots + x_{l,n} \be_n$, put
$b_l = \pi^{c_l} x_{l,l}^\sigma/x_{l,l}$,
and put $\bv_{l+1} = F\bv_l - b_l \bv_l$.
Then $\bv_l$ lies in the span of $\be_{l+1}, \dots, \be_n$;
in particular, $\bv_n = 0$.

We then have
\[
(F-b_{n-1})\cdots(F-b_1)(F-b_0)\bv = 0;
\]
since $\bv$ is a cyclic vector, there is a unique way to write $F^n \bv$
as a linear combination of $\bv, F\bv, \cdots, F^{n-1}\bv$. Hence we have
an equality of operators
\[
(F-b_{n-1})\cdots(F-b_1)(F-b_0) = F^n + a_{n-1} F^{n-1} + \cdots + a_0,
\]
from which the equality of polygons may be read off directly.
\end{proof}

\begin{remark}
One can turn Lemma~\ref{L:use cyclic vector} around and use it to
prove the existence of Dieudonn\'e-Manin decompositions in the case
of $v_K$ trivial; for instance, this is
the approach in \cite[Theorem~5.6]{me-local}.
One of the essential difficulties in \cite{me-local} is that there is
no analogous way to ``read off'' the HN-polygon of a $\sigma$-module over
$\Gancon$;
this forces the approach to constructing the slope filtration over
$\Gancon$ to be somewhat indirect.
\end{remark}

Lemma~\ref{L:use cyclic vector} is sometimes inconvenient to apply, because
the calculus of cyclic vectors is quite ``nonlinear''. The following criterion
will prove to be more useful for our purposes.

\begin{lemma} \label{L:identify generic}
Let $M$ be a $\sigma$-module over $\Galg[\fp]$. Suppose that there
exists a basis $\be_1, \dots, \be_n$ of $M$ with the property that 
the matrix $A$ given by $F \be_j = \sum_i A_{ij} \be_i$ satisfies
$w(A D^{-1} - I_n)>0$ for some $n \times n$
diagonal matrix $D$ over $\Galg[\fp]$. 
Then the generic slopes of $M$ are equal to the valuations
of the diagonal entries of $D$. Moreover, there exists an invertible
matrix $U$ over $\Galg$ with $w(U-I_n) > 0$, $w(D U^\sigma D^{-1} - I_n) > 0$,
and $U^{-1} A U^\sigma = D$.
\end{lemma}
\begin{proof}
One can directly solve for $U$;
see \cite[Proposition~5.9]{me-local} for this calculation. 
Note that it does not matter whether
the $D^{-1}$ appears to the left or to the right of $A$, as a change of basis
will flip it over to the other side; the entries of $D$ will get hit by
$\sigma$ or its inverse, but their valuations will not change.
Alternatively, the existence of $U$ also follows from
Proposition~\ref{P:quant reverse} below.
\end{proof}

\begin{remark}
Lemmas~\ref{L:use cyclic vector} and~\ref{L:identify generic} suggest that one can read off the generic HN-polygon of a $\sigma$-module over $\Gamma[\fp]$
by computing the slopes of eigenvectors of a matrix via which $F$ acts
on some basis. This does not work in general, as observed by 
Katz \cite{katz}.
\end{remark}

\subsection{Descending the generic HN-filtration}

In the generic setting ($v_K$ trivial), we have the following
descent property for Harder-Narasimhan filtrations.

\begin{prop} \label{P:descend generic}
Let $M$ be a $\sigma$-module over $\Gamma[\fp]$. Then the 
Harder-Narasimhan filtration of
$M$, tensored up to $\Galg[\fp]$, gives the Harder-Narasimhan filtration
of $M \otimes \Galg[\fp]$.
\end{prop}
\begin{proof}
One can prove this by Galois descent, as in 
\cite[Proposition~5.10]{me-local}; here is an alternate argument.
It suffices to check that the first step $M'_1$ of the Harder-Narasimhan
filtration of $M' = M \otimes \Galg[\fp]$ descends to
$\Gamma[\fp]$; by Lemma~\ref{L:wedge descent} and 
Proposition~\ref{P:wedge first step}, 
we may by taking exterior powers reduce
to the case where $M'_1$ has rank 1.
In particular, the least slope $s_1$ must be an integer. By
twisting, we may assume that $s_1 = 0$.

By Proposition~\ref{P:lattice}, we can find a $\sigma$-stable 
$\Gamma$-lattice $N$ of $M$; put $N' = N \otimes \Galg$.
Then $N' \cap M'_1$ may be characterized
as the set of limit points, for the $\pi$-adic topology,
of sequences of the form $\{F^l \bv_l\}_{l=0}^\infty$ with
$\bv_l \in N'$ for each $l$. (This may be verified on a basis of
$d$-eigenvectors for appropriate $d$ thanks to Theorem~\ref{T:Dieudonne-Manin},
where it is evident.)

The characterization of $N' \cap M'_1$ we just gave is linear, so it cuts out
a rank one submodule of $N$ already over $\Gamma$. This yields
the desired result.
\end{proof}

\subsection{de Jong's reverse filtration}

We now consider the case of $\sigma$-modules over $\Gcon[\fp]$, in which
case we have a ``generic'' HN-filtration defined over $\Galg[\fp]$, and
a ``special'' HN-filtration defined over $\Galgancon$. These two filtrations
are not directly comparable, because they live over incompatible overrings
of $\Gcon[\fp]$. 
To compare them, we must use a ``reverse filtration'' that meets both halfway;
the construction is due to de Jong \cite[Proposition~5.8]{dejong}, but our
presentation follows 
\cite[Proposition~5.11]{me-local} (wherein it is called the
``descending generic filtration'').

We first need a lemma that descends some eigenvectors from $\Galg$ to
$\Galgcon$; besides de Jong's \cite[Proposition~5.8]{dejong}, this
generalizes a lemma of Tsuzuki \cite[4.1.1]{tsuzuki-etale}.
\begin{lemma} \label{L:negative eigen}
Let $M$ be a $\sigma$-module over $\Galgcon[\fp]$ all of whose 
generic slopes are nonpositive. Let $\bv \in M \otimes \Galg[\fp]$
be an eigenvector of slope $0$. Then $\bv \in M$.
\end{lemma}
\begin{proof}
By Proposition~\ref{P:lattice}, we can find an $F$-stable $\Galgcon$-lattice
$N^\dual$ of $M^\dual$; the dual lattice $N$ is an $F^{-1}$-stable
$\Galgcon$-lattice of $M$. Let $\be_1, \dots, \be_n$ be a basis of $N$,
define the $n \times n$ matrix $A$ over $\Galgcon$ by the equation
$F^{-1} \be_j = \sum_i A_{ij} \be_i$, and choose $r>0$ such that
$A$ is invertible over $\Galg_r$.

Let $\bv \in M \otimes \Galg[\fp]$ be an eigenvector of slope 0;
in showing that $\bv \in M$, there is no harm 
in assuming (by multiplying by a power of $\pi$ as needed)
that $\bv \in N \otimes \Galg$.
Write $\bv = \sum x_i \be_i$ with $x_i \in \Galg$.
Then for each $l \geq 0$, we have
\begin{align*}
\min_i \min_{m \leq l} \{ v_{m,r}(x_i) \}
&\geq w_r(A) + \min_i \min_{m \leq l} \{ v_{m,r}(x_i^{\sigma^{-1}})\} \\
&\geq w_r(A) + q^{-1} \min_i \min_{m \leq l} \{ v_{m,r}(x_i)\}.
\end{align*}
It follows that
\[
\min_i \min_{m \leq l} \{ v_{m,r}(x_i) \} \geq q w_r(A)/(q-1),
\]
and so $x_i \in \Galgcon$. Hence $\bv \in M$, as desired.
\end{proof}

We now proceed to construct the reverse filtration.
\begin{defn}
Let $M$ be a $\sigma$-module over $\Galg[\fp]$
with slope decomposition
$P_1 \oplus \cdots \oplus P_l$,
labeled so that $\mu(P_1) > \cdots > \mu(P_l)$. Define the
\emph{reverse filtration} of $M$ as the semistable filtration
$0 = M_0 \subset M_1 \subset \cdots \subset M_l = M$ with
$M_i = P_1 \oplus \cdots \oplus P_i$ for $i=1, \dots, l$.
By construction, its Newton polygon coincides with the generic
Newton polygon of $M$.
\end{defn}

\begin{prop} \label{P:reverse filtration}
Let $M$ be a $\sigma$-module over $\Galgcon[\fp]$. Then the reverse
filtration of $M \otimes \Galg[\fp]$ descends to $\Galgcon[\fp]$.
\end{prop}
\begin{proof}
Put $M' = M \otimes \Galg[\fp]$, and let
$0 = M'_0 \subset M'_1 \subset \cdots \subset M'_l = M'$ be the
reverse filtration of $M'$. It suffices to show that
$M'_1$ descends to $\Galgcon[\fp]$;
by Lemma~\ref{L:wedge descent}, we may reduce to the case where
$\rank M'_1 = 1$ by passing from $M$ to an exterior power.
By twisting, we may then reduce to the case $\mu(M'_1) = 0$.
By Proposition~\ref{P:rank 1 class}, $M'_1$ is then generated
by an eigenvector of slope 0; by Lemma~\ref{L:negative eigen},
that eigenvector belongs to $M$. Hence $M'_1$ descends to $M$,
proving the claim.
(Compare \cite[Proposition~5.11]{me-local}.)
\end{proof}

\begin{remark}
The reverse filtration actually descends all the way to 
$\Gancon[\fp]$ whenever $K$ is perfect, but we will not need this.
\end{remark}

It may also be useful for some applications to have a quantitative version
of Proposition~\ref{P:reverse filtration}.
\begin{prop} \label{P:quant reverse}
Let $M$ be a $\sigma$-module over $\Galgcon[\fp]$. Suppose that for some $r>0$, there
exists a basis $\be_1, \dots, \be_n$ of $M$ with the property that 
the $n \times n$ matrix $A$ given by $F \be_j = \sum_i A_{ij} \be_i$ has entries in
$\Galg_r[\fp]$. Suppose further that
$w(A D^{-1} - I_n)>0$ and $w_r(AD^{-1} - I_n)>0$ for some $n \times n$
diagonal matrix $D$ over $\calO[\fp]$, with
$w(D_{11}) \geq \cdots \geq w(D_{nn})$. 
Then there exists an invertible $n \times n$ matrix $U$ over $\Galg_{qr}$ with
$w(U-I_n) > 0$,
$w_r(U-I_n) > 0$, $w(DU^\sigma D^{-1} - I_n) > 0$, $w_r(D U^\sigma D^{-1} - I_n) > 0$, such that $U^{-1} A U^\sigma D^{-1} - I_n$ is upper triangular nilpotent.
\end{prop}
\begin{proof}
Put $c_0 = w_r(AD^{-1} - I_n)$ and $U_0 = I_n$. Given $U_l$, put
$A_l = U_l^{-1} A U_l^\sigma$, and write $A_l D^{-1} -I_n = B_l + C_l$ with $B_l$ upper triangular nilpotent and $C_l$ lower triangular. Suppose that
$w_r(A_l D^{-1} - I_n) \geq c_0$; put
$c_l = w_r(C_l)$.
Choose a matrix $X_l$ over $\Gamma_{qr}$ 
with 
\begin{gather*}
C_l + X_l = D X_l^\sigma D^{-1}, \\
\min\{w(X_l), w(D X_l^\sigma D^{-1})\} \geq w(C_l), \\
\min\{w_r(X_l), w_r(D X_l^\sigma D^{-1})\} \geq w_r(C_l).
\end{gather*}
(This amounts to solving a system of equations of the form $c + x = \lambda^{-1}
x^\sigma$ for $\lambda \in \calO$; the analysis is as in
Proposition~\ref{P:h1}.)

Put $U_{l+1} = U_l (I_n + X_l)$.
We then have
\[
A_{l+1}D^{-1} = (I_n - X_l + X_l^2 - \cdots)
(I_n + B_l + C_l)(I_n + D X_l^\sigma D^{-1}),
\]
whence $w_r(C_{l+1}) \geq c_l + c_0$. In other words, $c_l \geq (l+1)c_0$;
hence the $U_l$ converge under $w_{qr}$ to a limit $U$ such that $U^{-1}
A U^\sigma D^{-1} - I_n$ is upper triangular nilpotent, as desired.
\end{proof}

\begin{remark}
For more on de Jong's original application of the reverse filtration,
see Section~\ref{subsec:fullfaith}.
\end{remark}

\subsection{Comparison of polygons}

Using the reverse filtration, we obtain the fundamental comparison between
the generic and special Newton
polygons of a $\sigma$-module over $\Gcon[\fp]$.
\begin{prop} \label{P:comparison}
Let $M$ be a $\sigma$-module over $\Gcon[\fp]$. Then the
special HN-polygon of $M$ (i.e., the HN-polygon of $M \otimes \Galgancon$)
lies above the generic HN-polygon (i.e., the HN-polygon of $M \otimes \Galg[\fp]$).
\end{prop}
\begin{proof}
Tensor up to $\Galgcon[\fp]$,
then apply Proposition~\ref{P:filtrations} to the 
reverse filtration (Proposition~\ref{P:reverse filtration}).
\end{proof}

The case where the polygons coincide is especially pleasant, and
will be crucial to our later results.
\begin{theorem} \label{T:coincide}
Let $M$ be a $\sigma$-module over $\Gcon[\fp]$ whose generic and
special HN-polygons coincide. Then the generic and special HN-filtrations
of $M \otimes \Galg[\fp]$ and $M \otimes \Galgancon$, respectively,
are both obtained by base change from an exhaustive filtration of $M$.
\end{theorem}
\begin{proof}
It suffices to check that the first steps of the generic and
special HN-filtrations descend and coincide;
by Lemma~\ref{L:wedge descent}, we may reduce to the case where the least
slope of the common polygon occurs with multiplicity 1.
Choose a basis $\be_1, \dots, \be_n$ of $M$, let $\bv \in
M \otimes \Galg[\fp]$ be a generator
of the first step of the HN-filtration of $M \otimes \Galg[\fp]$, and write
$\bv = a_1 \bv_1 + \cdots + a_n \bv_n$. By Proposition~\ref{P:descend
generic}, $a_i/a_j \in \Gamma[\fp]$ for any $i,j$ with $a_j \neq 0$.

By Corollary~\ref{C:filt splits}, the reverse filtration splits over
$\Galgancon$; by Proposition~\ref{P:h1}(b1), it also splits over 
$\Galgcon[\fp]$. Hence $a_i/a_j \in \Galgcon[\fp]$ for any $i,j$ 
with $a_j \neq 0$. Since $\Gamma[\fp] \cap \Galgcon[\fp] = \Gcon[\fp]$,
we have $a_i/a_j \in \Gcon[\fp]$ for any $i,j$ with $a_j \neq 0$.
Hence the first step of the generic HN-filtration descends to
$\Gcon[\fp]$; let $M_1$ be the corresponding $\sigma$-submodule of
$M\otimes \Gcon[\fp]$. 

Let $M'_1$ be the first step of the HN-filtration
of $M \otimes \Galgancon$; then 
as in the proof of Proposition~\ref{P:HN exists},
$M'_1$ is the maximal $\sigma$-submodule of $M \otimes \Galgancon$ of slope $\mu(M'_1) =
\mu(M_1)$. In particular, $M_1 \otimes \Galgancon \subseteq M'_1$,
and by Lemma~\ref{L:same rank}, $M_1 \otimes \Galgancon$ and $M'_1$
actually coincide. Hence the first step of the special HN-filtration
is also a base extension of $M_1$. This yields the claim.
(Compare \cite[Proposition~5.16]{me-local}.)
\end{proof}

\begin{remark}
The conditions of Theorem~\ref{T:coincide} may look restrictive,
and indeed they are: many ``natural'' examples of $\sigma$-modules
over $\Gcon[\fp]$ do not have the same special and generic HN-polygons
(e.g., the example of Section~\ref{subsec:gen spec2}).
However, Theorem~\ref{T:coincide} represents a critical step in the
descent process for slope filtrations, as it allows us to move information
from the generic paradigm into the special paradigm: 
specifically, the descent from
algebraically closed $K$ down to general $K$ is much easier in the generic
setting. Of course, in order to use this link, we must be able to
force ourselves into the setting of Theorem~\ref{T:coincide};
this is done in the next chapter.
\end{remark}

Note that Theorem~\ref{T:coincide} implies that the generic and special
HN-polygons of $M$ coincide if and only if the generic HN-filtration
descends to $\Galgcon[\fp]$. In some applications, it may be more useful
to have a quantitative refinement of that statement; we can give one
(in imitation of the proof of Proposition~\ref{P:raise Newton})
as follows.
\begin{prop}
Let $M$ be a $\sigma$-module over $\Galgcon[\fp]$. Suppose that for some $r>0$, there
exists a basis $\be_1, \dots, \be_n$ of $M$ with the property that 
the $n \times n$ matrix $A$ given by $F \be_j = \sum_i A_{ij} \be_i$ has entries in
$\Galg_r[\fp]$. Suppose further that
$w(A D^{-1} - I_n)>0$ and $w_r(AD^{-1} - I_n)$ for some $n \times n$
diagonal matrix $D$ over $\calO[\fp]$, with
$w(D_{11}) \geq \cdots \geq w(D_{nn})$.
Let $U$ be an invertible $n \times n$ matrix over $\Galg$
such that
$w(U-I_n) > 0$, $w(DU^\sigma D^{-1} - I_n) > 0$, and $U^{-1} A U^\sigma = D$
(as in Lemma~\ref{L:identify generic}).
Then the generic and special HN-polygons of $M$ coincide if and only if
\begin{equation} \label{eq:quant bound}
w_r(U-I_n) > 0, \qquad w_r(DU^\sigma D^{-1} - I_n) > 0.
\end{equation}
\end{prop}
\begin{proof}
The conditions on $U$ imply that $M$ admits a basis of eigenvectors over
$\Galgcon[\fp]$; the slopes of these eigenvectors then give both the generic and special HN-slopes. Conversely, if the generic and special HN-polygons of $M$ coincide, then $U$ is invertible over $\Galgcon$ by Theorem~\ref{T:coincide}.
It thus remains to prove that if $U$ has entries in $\Galgcon$,
then in fact \eqref{eq:quant bound} holds.

We start with a series of reductions.
By Proposition~\ref{P:quant reverse}, we may reduce to the case where
$AD^{-1} - I_n$ is upper triangular nilpotent. Considering all matrices now as block matrices by grouping rows and columns where the diagonal entries of $D$ have the same valuation, we see that the matrix $U$ must now be block upper triangular, with diagonal blocks invertible over $\calO$.
Neither this property nor
\eqref{eq:quant bound} is disturbed by multiplying $U$ by a block diagonal matrix over $\calO$, so we may reduce to the case where $U$ is block upper triangular with identity matrices on the diagonal. Finally, note that at this point it suffices to check the case where $n=2$, $w(D_{11}) > w(D_{22})$, and
\[
AD^{-1} = \begin{pmatrix} 1 & a \\ 0 & 1 \end{pmatrix}, \qquad
U = \begin{pmatrix} 1 & u \\ 0 & 1 \end{pmatrix},
\]
as the general case follows by repeated application of this case.

Put $\lambda \in D_{11} D_{22}^{-1} \in \pi \calO$.
Then what we are to show is that given $a \in \Galg_r$, $u \in \Galgcon$ with
\begin{equation} \label{eq:quant system}
w_r(a) > 0, \qquad
w(u) > 0, \qquad
\lambda a + u = \lambda u^\sigma,
\end{equation}
we must have $u \in \Galg_r$ and $w_r(u) > 0$. 
Before verifying this,
we make one further simplifying reduction, using the fact that there
is no harm in replacing $a$ by $a-b + \lambda^{-1} b^{\sigma^{-1}}$ as long as
$w_r(b) > 0$ and $w_r(\lambda^{-1} b^{\sigma^{-1}}) > 0$.

Note that 
for $b = \pi^m [\overline{x}]$ with $\overline{x} \in K^{\alg}$,
the condition that
\[
w_r(b) \leq w_r(\lambda^{-1} b^{\sigma^{-1}} )
\]
is equivalent to
\[
v_K(\overline{x}) \leq -\frac{q}{(q-1)r}w(\lambda).
\]
Moreover, this condition and the bound $w_r(b) \geq 0$ together imply that $w(\lambda^{-1} b^{\sigma^{-1}} ) > 0$.

Now write $a = \sum_{i=0}^\infty \pi^i [\overline{a_i}]$,  and put
$c = q w(\lambda)/(q-1)r$. For
each $i$, let $j_i$ be the smallest nonnegative integer such that
$q^{-j_i} v_K(\overline{a_i}) > -c$. We may then
replace $a$ by
\[
a' = \sum_{i=0}^\infty \lambda^{-1} \lambda^{-\sigma^{-1}} \cdots
\lambda^{-\sigma^{j_i-1}} (\pi^i [\overline{a_i}])^{\sigma^{-j_i}}
\]
without disturbing the truth of \eqref{eq:quant system}.
In particular, we may reduce to the case where $v_n(a)
> -c$ for all $n \geq 0$.

Under these conditions, put $m = w(\lambda) > 0$, and note that if $v_n(u) \leq 
-c/q$ for some $u$, then $v_{n+m}(\lambda u^\sigma) = q v_n(u) \leq -c$, so the
equation $\lambda a + u = \lambda u^\sigma$ implies $v_{n+m}(u) = q v_n(u)$. By
induction on $l$, we have
\[
v_{n+lm}(u) = q^l v_n(u)
\]
for all nonnegative integers $l$, but this contradicts the hypothesis that
$u \in \Galgcon$. Consequently, we must have $v_n(u) > -c/q$ for all $n$.

Since $-c/q \geq -c(q-1)/q = -w(\lambda)/r$, the bound $v_n(u) > -c/q$
implies that $v_{n,r}(u) > 0$ for $n \geq m$. Since $u \equiv 0 \pmod{\lambda}$,
we have  $w_r(u) > 0$,
as desired.
\end{proof}

\section{Descents}
\label{sec:descents}

We now show that the formation of the Harder-Narasimhan filtration commutes
with base change, thus establishing the slope filtration theorem;
the strategy is to show that a $\sigma$-module over $\Gancon$ admits
a model over $\Gcon$ whose special and generic Newton polygons coincide,
then invoke Theorem~\ref{T:coincide}.
The material here is derived from \cite[Chapter~6]{me-local},
but our presentation here is much more streamlined.

\subsection{A matrix lemma}

The following lemma is analogous to \cite[Proposition~6.8]{me-local}, but in
our new approach, we can prove much less and still eventually get the desired
conclusion; this simplifies the matrix calculation considerably.
\begin{lemma} \label{L:frob turnover}
For $r>0$, suppose that $\Gamma = \Gamma^K$ contains enough $r_0$-units
for some $r_0 >qr$.
Let $D$ be an invertible $n \times n$ matrix over $\Gamma_{[r,r]}$, and put $h = -w_r(D) - w_r(D^{-1})$.
Let $A$ be an $n \times n$ matrix over $\Gamma_{[r,r]}$ such that 
$w_r(AD^{-1} - I_n) > h/(q-1)$. Then there exists an invertible $n \times n$ matrix
$U$ over $\Gamma_{[r,qr]}$ such that $U^{-1}AU^\sigma D^{-1} - I_n$ has entries in
$\pi \Gamma_r$ and $w_r(U^{-1} A U^\sigma D^{-1} - I_n) > 0$.
\end{lemma}
\begin{proof}
Put $c_0 = w_r(AD^{-1} -I_n) - h/(q-1)$.
Define sequences of matrices $U_0, U_1,\dots$ and $A_0, A_1, \dots$ as follows.
Start with $U_0 = I$. Given an invertible $n \times n$ matrix $U_l$ over
$\Gamma_{[r,qr]}$, put $A_l = U_l^{-1} A U_l^\sigma$. 
Suppose that $w_r(A_l D^{-1} - I_n) \geq c_0 + h/(q-1)$;
put
\[
c_l = \min_{m \leq 0} \{v_{m,r}(A_l D^{-1} - I_n) \} - h/(q-1),
\]
so that $c_l \geq c_0 > 0$.

Let $\sum u_{ijlm} \pi^m$ be a semiunit presentation
of $(A_lD^{-1} - I_n)_{ij}$. 
Let $X_l$ be the $n \times n$ matrix with
$(X_l)_{ij} = \sum_{m \leq 0} u_{ijlm} \pi^m$, and put
$U_{l+1} = U_l (I_n + X_l)$.
By Lemma~\ref{L:sum}, for $m \leq 0$ and $s \in [r,qr]$,
\begin{align*}
w_s(u_{ijlm} \pi^m) &\geq (s/r) w_r(u_{ijlm} \pi^m) \\
&\geq (s/r) \min_{k \leq m} \{v_{k,r}(A_l D^{-1} - I_n)\} \\
&\geq (s/r)(c_l + h/(q-1));
\end{align*}
hence $U_{l+1}$ is also invertible over $\Gamma_{[r,qr]}$. 
Moreover, 
\begin{align*}
w_r(DX_l^\sigma D^{-1}) &\geq w_r(D) + w_r(X_l^\sigma) + w_r(D^{-1}) \\
&= w_{qr}(X_l) - h \\
&\geq q(c_l + h/(q-1)) - h \\
&\geq qc_l + h/(q-1).
\end{align*}
Since
\[
A_{l+1} D^{-1} = (I_n + X_l)^{-1} (A_l D^{-1}) (I_n + D X_l^\sigma D^{-1}),
\]
we then have $w_r(A_{l+1} D^{-1} - I_n) \geq c_0 + h/(q-1)$, so the 
iteration may continue.

We now prove by induction that $c_l \geq (l+1)c_0$ for $l=0,1,\dots$;
this is clearly true for $l=0$. Given the claim for $l$, write
\[
(I_n + X_l)^{-1} A_l D^{-1} = 
I_n + (A_l D^{-1} - I_n - X_l) - X_l (A_l D^{-1} - I_n)
+ X_l^2 (I_n + X_l)^{-1} A_l D^{-1}.
\]
Note that 
\begin{align*}
v_m(A_l D^{-1} - I_n - X_l) &= \infty \qquad (m \leq 0) \\
w_r(X_l(A_l D^{-1} - I)) &\geq (l+2)c_0 + 2h/(q-1) \\
w_r(X_l^2 (I_n + X_l)^{-1} A_l D^{-1}) &\geq 2(l+1)c_0 + 2h/(q-1).
\end{align*}
Putting this together, this means
\[
v_{m,r}((I_n+X_l)^{-1} A_l D^{-1} - I_n)
\geq (l+2)c_l + h/(q-1) \qquad (m \leq 0).
\]
Since $w_r(DX_l^\sigma D^{-1}) \geq qc_l + h/(q-1) \geq (l+2) c_0 + h/(q-1)$,
we have $v_{m,r}(A_{l+1} D^{-1} - I_n) \geq
(l+2) c_0 + h/(q-1)$ for $m \leq 0$, i.e., $c_{l+1} \geq (l+2) c_0$.

Since $w_s(X_l) \geq (s/r)(c_l + h/(q-1))$ for $s \in [r,qr]$,
and $c_l \to \infty$ as $l \to \infty$,
the sequence $\{U_l\}$ converges to a limit $U$, which
is an invertible $n \times n$ matrix over $\Gamma_{[r,qr]}$ satisfying
$w_r(U^{-1} A U^\sigma D^{-1} - I_n) \geq c_0 + h/(q-1) > 0$.
Moreover, by construction, we have $v_m(U^{-1} A U^\sigma 
D^{-1} - I_n) = \infty$
for $m \leq 0$; by Corollary~\ref{C:subring},
$U^{-1} A U^\sigma D^{-1} - I_n$ has entries in
$\pi \Gamma_r$, as desired.
\end{proof}

\subsection{Good models of $\sigma$-modules}

We now give a highly streamlined version of some arguments from
\cite[Chapter~6]{me-local}, that produce ``good integral models''
of $\sigma$-modules over $\Gancon$.

\begin{lemma} \label{L:frob turnover2}
In Lemma~\ref{L:frob turnover}, suppose that $A$ and $D$
are both invertible over $\Gamma_{\an,r}$. Then $U$ is invertible
over $\Gamma_{\an,qr}$.
\end{lemma}
\begin{proof}
Put $B= U^{-1} A U^\sigma D^{-1}$, so that $B$ is invertible
over $\Gamma_r$. In the equation
\[
U^{-1} A U^\sigma = B D,
\]
the matrices $A, U^\sigma, B, D$ are all invertible over $\Gamma_{[r/q,r]}$,
so $U$ must be as well. Since the entries of $U$ and $U^{-1}$ already belong
to $\Gamma_{[r,qr]}$, they in fact belong to $\Gamma_{[r/q,qr]}$
by Corollary~\ref{C:subring2}. Repeating
this argument, we see that $U$ is invertible over $\Gamma_{[r/q^i,qr]}$
for all positive integers $i$, yielding the desired result.
\end{proof}

\begin{prop} \label{P:good model}
Let $M$ be a $\sigma$-module over $\GKancon$. Then there exists
a finite separable extension $L$ of $K$ and a
$\sigma$-module $N$ over $\GLcon[\fp]$ such that
$N \otimes \GLancon \cong M \otimes \GKancon$ and
the generic and special HN-polygons of $N$ coincide.
\end{prop}
\begin{proof}
Note that it is enough to do this with $L$ pseudo-finite separable,
since applying $F$ allows to pass from $L$ to $L^q$.
Also, there is no harm in assuming that the slopes of $M$ are integers,
after applying $[a]_*$ as necessary.

Let $\be_1, \dots, \be_n$ be a basis of $M$, and define the invertible
$n \times n$ matrix $A$ over $\GKancon$ by $F\be_j = \sum_i A_{ij} \be_i$.
By Theorem~\ref{T:Dieudonne-Manin}, there exists an invertible $n \times n$
matrix $U$ over $\Galgancon$ and a diagonal $n \times n$ matrix $D$ whose
diagonal entries are powers of $\pi$, such that $U^{-1} A U^\sigma = D$.
Put $h = \max_{i,j} \{w(D_{ii}) - w(D_{jj})\} = -w(D) - w(D^{-1})$.

Choose $r>0$ such that $\Gamma$ has enough $r_0$-units for some $r_0 > qr$, 
$A$ is defined and
invertible over $\GK_{\an,r}$, and $U$ is defined and invertible over
$\Galg_{\an,qr}$. 
Let $M_r$ be the $\GK_{\an,qr}$-span of the $\be_i$.
By Lemma~\ref{L:dense}, we can choose $L$ pseudo-finite
separable over $K$ such that $\Gamma^L$ has enough $r_0$-units, and such that
there exists an $n \times n$ matrix $V$ over $\Gamma^L_{[r,qr]}$ with
$w_s(V - U) > -w_s(U^{-1}) + qh/(q-1)$ for $s \in [r,qr]$.
Since
\[
V^{-1} A V^\sigma D^{-1} = (U^{-1} V)^{-1} D (U^{-1} V)^\sigma D^{-1},
\]
it follows that $w_r(V^{-1} A V^\sigma D^{-1} - I_n) > h/(q-1)$.

By Lemma~\ref{L:frob turnover}, 
there exists an invertible $n \times n$ matrix $W$
over $\Gamma^L_{[r,qr]}$ for which the matrix $B = (VW)^{-1} A (VW)^\sigma$
is such that $B D^{-1} - I_n$ 
has entries in $\pi \Gamma^L_r$ and $w_r(BD^{-1}-I_n) > 0$.
By Lemma~\ref{L:frob turnover2}, the matrix $VW$ is actually invertible
over $\Gamma^L_{\an,qr}$.

%
Let $N$ be the $\sigma$-module over $\GLcon[\fp]$-module generated by
$\bv_1, \dots, \bv_n$ with Frobenius action defined by
$F\bv_j = \sum_i B_{ij} \bv_i$; then by what we have just shown,
$N \otimes \GLancon \cong M \otimes \GLancon$.
By Lemma~\ref{L:identify generic}, the generic HN-slopes of $N$ are
the $w(D_{ii})$, which by construction are also the special HN-slopes of $N$.
(Compare \cite[Proposition~6.9]{me-local}.)
\end{proof}

\begin{remark}
We do not know whether it is possible to establish 
Proposition~\ref{P:good model} without the finite extension $L$ of $K$.
\end{remark}

\subsection{Isoclinic $\sigma$-modules}

Before proceeding to the general descent problem for HN-filtrations,
we analyze the isoclinic case.

\begin{defn}
Let $M$ be a $\sigma$-module over $\Gcon[\fp]$. 
We say that $M$ is \emph{isoclinic} (of slope $\mu(M)$) if
$M \otimes \Galg[\fp]$ is isoclinic. By
Proposition~\ref{P:comparison}, this implies that
$M \otimes \Galgancon$ is also isoclinic.
\end{defn}

\begin{prop} \label{P:isoclinic h0}
Let $M$ be a $\sigma$-module over $\Gcon[\fp]$ which is unit-root
(isoclinic of slope $0$). Then
\[
H^0(M) 
= H^0(M \otimes \Gamma[\fp]) = H^0(M \otimes \Gancon[\fp]);
\]
in particular, if $K$ is algebraically closed, then
$M$ admits a basis of eigenvectors.
\end{prop}
\begin{proof}
There is no harm (thanks to Corollary~\ref{C:intersect})
 in assuming from the outset that $K$ is algebraically closed.
In this case, the eigenvectors of $M \otimes \Galg[\fp]$ all have slope 0
by Corollary~\ref{C:hom ext2};
by Lemma~\ref{L:negative eigen}, they all belong to $M$.
Hence $M$ admits a basis of eigenvectors; in particular,
$M \cong M_{0,1}^{\oplus n}$. Then $H^0(M) = H^0(M \otimes \Gancon[\fp])$
by Proposition~\ref{P:h0}.
\end{proof}

\begin{theorem} \label{T:isoclinic lattice}
\begin{enumerate}
\item[(a)]
The base change functor, from isoclinic $\sigma$-modules over $\Gcon[\fp]$ of
some prescribed slope $s$ to isoclinic $\sigma$-modules over $\Gamma[\fp]$ of
slope $s$, is fully faithful.
\item[(b)]
The base change functor, from isoclinic $\sigma$-modules over $\Gcon[\fp]$ of
slope $s$ to isoclinic $\sigma$-modules over $\Gancon$ of
slope $s$, is an equivalence of categories.
\end{enumerate}
In particular, any isoclinic $\sigma$-module over $\Galgcon[\fp]$ is isomorphic
to a direct sum of standard $\sigma$-modules of the same slope (and hence
is semistable by Proposition~\ref{P:which stable}).
\end{theorem}
\begin{proof}
To see that the functors are fully faithful, let $M$ and $N$
be isoclinic $\sigma$-modules over $\Gcon[\fp]$ of the same slope
$s$. Then $M^\dual \otimes N$ is unit-root, and $\Hom(M,N) = H^0(M^\dual
\otimes N)$, so the full faithfulness assertion follows from
Proposition~\ref{P:isoclinic h0}.

To see that the functor in (b) is essentially surjective,
apply Proposition~\ref{P:good model} to produce an $F$-stable
$\GLcon[\fp]$-lattice $N$ in $M \otimes \GLancon$ for some finite separable
extension $L$ of $K$, which we may take to be Galois. Note that $N$ is unique by
full faithfulness of the base change functor (from $\GLcon[\fp]$ to $\GLancon$); 
hence the action of $G = \Gal(L/K)$ on $M \otimes
\GLancon$ induces an action on $N$. By ordinary Galois descent,
there is a unique $\Gcon[\fp]$-lattice $M_b$ of $N$ such that $N = 
M_b \otimes
\GLcon[\fp]$; because of the uniqueness, $M_b$ is $F$-stable. This yields 
the desired
result.
(Compare \cite[Theorem~6.10]{me-local}.)
\end{proof}

\begin{remark}
The functor in (a) is not essentially surjective in general. For instance,
if $K = k((t))$ with $k$ perfect, $M$ is an isoclinic $\sigma$-module over $\Gcon$,
and $b_n$ is the highest break of the representation of $\Gal(K^{\sep}/K)$
on the images modulo $\pi^n$ of the eigenvectors of $M \otimes \Galg$, then
$b_n/n$ is bounded. By contrast, if $M$ is an isoclinic $\sigma$-module over $\Gamma$,
$b_n/n$ need not be bounded.
\end{remark}

Incidentally,
Theorem~\ref{T:isoclinic lattice} allows us to give a more succinct
characterization of isoclinic $\sigma$-modules, which one could take
as an alternate definition. 
\begin{prop} \label{P:isoclinic integral}
Let $c,d$ be integers with $d>0$.
Then a 
$\sigma$-module $M$ over $\Gancon$ is isoclinic of slope $s = c/d$ if and 
only if $M$ admits a $\Gcon$-lattice $N$ such that $\pi^{-c} F^d$
maps some (any) basis of $N$ to another basis of $N$.
\end{prop}
\begin{proof}
If such a lattice exists, then we may apply Proposition~\ref{P:nonnegative}
to $([d]_* M)(-c)$ to deduce that its generic HN-slopes are
all zero. By Proposition~\ref{P:slope arith}, $M$ is isoclinic of slope $c/d$.

Conversely, suppose $M$ is isoclinic of slope $s$; then $([d]_* M)(-c)$
is isoclinic of slope 0. By Theorem~\ref{T:isoclinic lattice}, 
$([d]_* M)(-c)$ admits a unique $F$-stable $\Gcon[\fp]$-lattice $N'$.
By Proposition~\ref{P:lattice}, $N'$ in turn admits an $F$-stable
$\Gcon$-lattice $N$; by Proposition~\ref{P:nonnegative}, the Frobenius
on $N$ carries any basis to another basis.
\end{proof}

\subsection{Descent of the HN-filtration}

At last, we are ready to establish the slope filtration theorem
\cite[Theorem~6.10]{me-local}.

\begin{theorem} \label{T:descend special}
Let $M$ be a $\sigma$-module over $\Gancon$. Then there exists a unique
filtration $0 = M_0 \subset M_1 \subset \cdots \subset M_l = M$ of $M$
by saturated $\sigma$-submodules with the following properties.
\begin{enumerate}
\item[(a)] For $i=1, \dots, l$, the quotient $M_i/M_{i-1}$ is isoclinic
of some slope $s_i$.
\item[(b)] $s_1 < \cdots < s_l$.
\end{enumerate}
Moreover, this filtration coincides with the Harder-Narasimhan filtration of $M$.
\end{theorem}
\begin{proof}
Since isoclinic $\sigma$-modules are semistable by Theorem~\ref{T:isoclinic
lattice}, any filtration as in (a) and (b) 
is a Harder-Narasimhan filtration. In particular, the filtration is unique
if it exists.

To prove existence, it suffices to show that the HN-filtration of 
$M' = M \otimes \Galgancon$ descends to $\Gancon$. Let $M'_1$ be the first
step of that filtration; by Lemma~\ref{L:wedge descent}, it is enough to check 
that $M'_1$ descends to $\Gancon$ in the case $\rank(M'_1) = 1$.

By Proposition~\ref{P:good model}, there exists a finite separable extension
$L$ of $K$ and a $\sigma$-module $N$ over $\GLcon[\fp]$ such that
$M \otimes \GLancon \cong N \otimes \GLancon$, and $N$ has the same
generic and special HN-polygons. Of course there is no harm in assuming $L$ is Galois
with $\Gal(L/K) = G$.
By Theorem~\ref{T:coincide}, $M'_1$ descends to $\GLancon$; let $M_1$
be the descended submodule of $M \otimes \GLancon$.

Let $\be_1, \dots, \be_n$ be a basis of $M$, let $\bv$ be a generator
of $M_1$, and write $\bv = a_1 \be_1 + \cdots + a_n \be_n$ with $a_i \in \GLancon$.
Then for each $i,j$ with $a_j \neq 0$, 
$a_i/a_j \in \Frac \GLancon$ is invariant under $G$.
By Corollary~\ref{C:invariant}, $a_i/a_j \in \Frac \GKancon$ for each $i,j$ with
$a_j \neq 0$; clearing denominators, we obtain a nonzero element of $M_1 \cap M$.
Hence $M_1$ descends to $\Gancon$, as desired.
(Compare \cite[Theorem~6.10]{me-local}.)
\end{proof}
\begin{cor}
For any extension $K'$ of $K$ (complete with respect to a valuation $v_{K'}$
extending $v_K$, such that $\Gamma^{K'}$ has enough units), and any
$\sigma$-module $M$ over $\GKancon$, the HN-filtration of $M \otimes 
\Gamma^{K'}_{\an,\con}$ coincides with the result of tensoring the
HN-filtration of $M$ with $\Gamma^{K'}_{\an,\con}$. In other words,
the formation of the HN-filtration commutes with base change.
\end{cor}
\begin{proof}
The characterization of the HN-filtration given by Theorem~\ref{T:descend
special} is stable under base change.
\end{proof}
\begin{cor}
A $\sigma$-module over $\Gancon$ is semistable if and only if it is isoclinic.
\end{cor}
\begin{proof}
If $M$ is an isoclinic $\sigma$-module over $\Gancon$, then $M$ is 
semistable by Proposition~\ref{P:which stable}. Conversely, if $M$
is not isoclinic, then by Theorem~\ref{T:descend special},
$M$ admits a nonzero $\sigma$-submodule $M_1$ with $\mu(M_1) < \mu(M)$,
so $M$ is not semistable.
\end{proof}

\section{Complements}
\label{sec:complements}

\subsection{Differentials and the slope filtration}

The slope filtration turns out to behave nicely with respect to differentials;
this is what allows the application to Crew's conjecture.
\begin{defn}
Let $S/R$ be an extension of topological rings. A \emph{module of
continuous differentials} is a topological $S$-module $\Omega^1_{S/R}$
equipped with a continuous $R$-linear derivation $d: S \to \Omega^1_{S/R}$,
having the following universal property: for any topological $S$-module
$M$ equipped with a continuous $R$-linear derivation $D: S \to M$,
there exists a unique morphism $\phi: \Omega^1_{S/R} \to M$ of topological
$S$-modules such that $D = \phi \circ d$. Since the definition is via
a universal property, the module of continuous differentials is unique
up to unique isomorphism if it exists at all.
\end{defn}

\begin{convention}
For the remainder of this section, assume that $\Gcon$, viewed
as a topological $\calO$-algebra via the levelwise topology,
has a module
of continuous differentials $\Omega^1_{\Gcon/\calO}$
which is finite free over $\Gcon$. In
this case, for any finite separable extension $L$ of $K$,
$\Omega^1_{\Gcon/\calO} \otimes_{\Gcon} \GLancon$ is a module
of continuous differentials of $\GLancon$ over $\calO[\fp]$.
\end{convention}

\begin{example}
If $K = k((t))$ as in Section~\ref{subsec:robba}, we have a module
of continuous differentials for $\Gcon$ over $\calO$ given by
$\Omega^1_{\Gcon/\calO} = \Gcon\,dt$ with the formal derivation
$d$ sending $\sum c_i t^i$ to $(\sum ic_i t^{i-1})dt$.
\end{example}

\begin{remark} \label{R:omega-sigma}
Note that $\Omega^1_{\Gcon/\calO}$ may be viewed as a $\sigma$-module,
via the action of $d\sigma$. Since $d\sigma(\omega) \equiv 0
\pmod{\pi}$ for any $\omega \in \Omega^1_{\Gcon/\calO}$,
by Proposition~\ref{P:nonnegative} the generic slopes of 
$\Omega^1_{\Gcon/\calO}$ as a $\sigma$-module are all positive.
By Proposition~\ref{P:comparison},
the special slopes of $\Omega^1_{\Gancon/\calO}$
as a $\sigma$-module are also nonnegative.
\end{remark}

\begin{defn} \label{D:nabla-module}
For $S = \Gcon, \Gcon[\fp], \Gancon$,
define a \emph{$\nabla$-module} over $S$ to be a 
finite free $S$-module equipped with an integrable connection
$\nabla: M \to M \otimes \Omega^1_{S/\calO}$.
(Integrability here means that the composition of $\nabla$
with the induced map $M \otimes \Omega^1_{S/\calO}
\to M \otimes \wedge^2_S \Omega^2_{S/\calO}$ is the zero map.)
Define a \emph{$(\sigma,\nabla)$-module} over $S$
to be a finite free $S$-module $M$ equipped with the structures of
both a $\sigma$-module and a $\nabla$-module, which commute as in
the following diagram:
\[
\xymatrix{
M \ar^-{\nabla}[r] \ar^{F}[d] & M \otimes \Omega^1_{S/\calO} \ar^{F \otimes 
d\sigma}[d] \\ M \ar^-{\nabla}[r] & M \otimes \Omega^1_{S/\calO}.
}
\]
\end{defn}

\begin{prop} \label{P:del stable1}
Let $M$ be a $(\sigma,\nabla)$-module over $\Gancon$. Then each
step of the HN-filtration of $M$ is a $(\sigma,\nabla)$-submodule
of $M$.
\end{prop}
\begin{proof}
Let $M_1$ be the first step of the HN-filtration, which is isoclinic
by Theorem~\ref{T:descend special}; it suffices to
check that $M_1$ is a $(\sigma,\nabla)$-submodule of $M$.
To simplify notation, write $N$ for $\Omega^1_{\Gancon/\calO}$.
Then the map $M_1 \to (M/M_1) \otimes N$ induced by $\nabla$
is a morphism of $\sigma$-modules, and the slopes of
$(M/M_1) \otimes N$ are all strictly greater than the slope of $M_1$
by Remark~\ref{R:omega-sigma}.
By Proposition~\ref{P:hom ext3}, the map $M_1 \to (M/M_1) \otimes N$
must be zero; that is, $M_1$ is a $\nabla$-submodule of $M$.
This proves the desired result.
\end{proof}

\begin{prop} \label{P:del stable2}
Let $M$ be an isoclinic $\sigma$-module over $\Gcon[\fp]$ such that
$M \otimes \Gancon$ admits the structure of a $(\sigma,\nabla)$-module.
Then $M$ admits a corresponding structure of a $(\sigma, \nabla)$-module.
\end{prop}
\begin{proof}
Write $N$ for $\Omega^1_{\Gcon[\fp]/\calO}$, so that $\nabla$
induces an additive map $M \otimes \Gancon \to (M \otimes N) \otimes \Gancon$.
Pick a basis $\be_1, \dots, \be_n$ of $M$, and define the map
$f: M \to M \otimes N$ by
\[
f\left(\sum c_i \be_i\right) = \sum_i \be_i \otimes dc_i.
\]
Then $\nabla - f$ is $\Gancon$-linear, so we may view it as an element
$\bv$ of $M^\dual \otimes M \otimes N \otimes \Gancon$. That element
satisfies $F\bv - \bv = \bw$ for some $\bw \in M^\dual \otimes M \otimes N$
by the commutation relation between $F$ and $\nabla$.
However, $M^\dual \otimes M$ is unit-root, so the 
(generic and special) slopes of $M^\dual \otimes M \otimes N$ are all
positive by Remark~\ref{R:omega-sigma}. By Theorem~\ref{T:Dieudonne-Manin} and
Proposition~\ref{P:h1}(b1), it follows that $\bv \in M^\dual \otimes M
\otimes N$, so $\nabla$ acts on $M$, as desired.
(Compare \cite[Theorem~6.12]{me-local} and \cite[Lemme~V.14]{berger-cst}.)
\end{proof}
\begin{remark}
We suspect there is a more conceptual way of saying this in terms
of splitting a certain exact sequence of $\sigma$-modules.
\end{remark}

\subsection{The $p$-adic local monodromy theorem}

We recall briefly how the slope filtration theorem,
plus a theorem of Tsuzuki, yields the $p$-adic local monodromy theorem
(formerly Crew's conjecture).

\begin{convention}
Throughout this section,
retain notation as in Section~\ref{subsec:robba}, i.e.,
$\Gancon$ is the Robba ring. Note that in this case, the integrability
condition in Definition~\ref{D:nabla-module} is vacuous, since
$\Omega^1_{\Gancon/\calO}$ is free of rank 1.
\end{convention}

\begin{defn}
We say a $\nabla$-module $M$ 
over $\Gancon$ is \emph{constant} if it is spanned by the kernel of
$\nabla$; equivalently, $M$ is isomorphic to a direct sum of trivial
$\nabla$-modules. (The ``trivial'' $\nabla$-module here means $\Gancon$
itself with the connection given by $d$.) 
We say $M$ is \emph{quasi-constant}
if $M \otimes \GLancon$ is constant for some finite separable
extension $L$ of $K$. We say a $(\sigma,\nabla)$-module is
(quasi-)constant if its underlying $\nabla$-module is (quasi-)constant.
\end{defn}

The key external result we need here is the following result essentially
due to Tsuzuki \cite{tsuzuki-unitroot}. A simplified proof of Tsuzuki's
result has been given by Christol \cite{christol}; however, see
\cite{andre-divizio} for the corrections of some errors in \cite{christol}.
\begin{prop} \label{P:tsuzuki}
Let $M$ be a unit-root $(\sigma, \nabla)$-module over $\Gcon[\fp]$.
Then $M \otimes \Gancon$ is quasi-constant; in fact, for some finite separable
extension $L$ of $K$, $M \otimes \GLcon[\fp]$ admits a basis in the kernel of
$\nabla$.
\end{prop}
\begin{proof}
Suppose first that
$M = M_0 \otimes \Gcon[\fp]$ for a unit-root $(\sigma, \nabla)$-module $M_0$
over $\Gcon$. Then our desired statement is Tsuzuki's theorem
\cite[Theorem~4.2.6]{tsuzuki-unitroot}. To reduce to this case,
apply Proposition~\ref{P:lattice} to obtain a $\sigma$-stable 
$\Gcon$-lattice $M_0$; by applying Frobenius repeatedly, we see that
$M_0$ must also be $\nabla$-stable. Thus Tsuzuki's theorem applies to give
the claimed result.
\end{proof}

\begin{defn}
We say a $\nabla$-module $M$ 
over $\Gancon$ is said to be \emph{unipotent} if it admits
an exhaustive filtration by $\nabla$-submodules whose successive quotients
are constant. We say $M$ is \emph{quasi-unipotent}
if $M \otimes \GLancon$ is unipotent for some finite separable
extension $L$ of $K$.
We say a $(\sigma,\nabla)$-module is
(quasi-)unipotent if its underlying $\nabla$-module is (quasi-)unipotent.
\end{defn}

\begin{theorem}[$p$-adic local monodromy theorem] \label{T:pLMT}
With notations as in Section~\ref{subsec:robba} (i.e.,
$\Gancon$ is the Robba ring), 
let $M$ be a $(\sigma, \nabla)$-module over $\Gancon$. Then $M$
is quasi-unipotent.
\end{theorem}
\begin{proof}
By Proposition~\ref{P:del stable1}, each step of the HN-filtration of $M$
is $\nabla$-stable. It is thus enough to check that each successive quotient
is quasi-unipotent; in other words, we may assume that $M$ itself is 
isoclinic.

Note that the definition of unipotence does not depend on the Frobenius
lift, so there is no harm in either applying the functor $[b]_*$ or in
twisting. We may thus reduce to the case where $M$ is isoclinic of slope
zero (i.e., is unit-root). Applying Proposition~\ref{P:tsuzuki}
then yields the desired result.
\end{proof}
\begin{example}
In Theorem~\ref{T:pLMT}, it can certainly happen that $M$ fails to be
quasi-constant. For instance, if $u^\sigma = qu$, and
$M$ has rank two with a basis
$\be_1, \be_2$ such that
\begin{gather*}
F\be_1 = \be_1, \qquad F\be_2 = q\be_2 \\
\nabla \be_1 = 0, \qquad \nabla \be_2 = \be_1 \otimes \frac{du}{u},
\end{gather*}
then $M$ is a $(\sigma, \nabla)$-module which is unipotent but not
quasi-constant.
\end{example}

\begin{remark}
The fact that quasi-unipotence follows from the existence of a slope
filtration as in Theorem~\ref{T:descend special}
was originally pointed out by Tsuzuki \cite[Theorem~5.2.1]{tsuzuki-slope}.
Indeed, that observation was the principal motivation for
the construction of slope filtrations of $\sigma$-modules in
\cite{me-local}.
\end{remark}

\begin{remark} \label{R:compare pLMT}
We remind the reader that Theorem~\ref{T:pLMT} has also been proved 
(independently) by Andr\'e \cite{andre} and by Mebkhout \cite{mebkhout},
using the index theory for $p$-adic differential equations developed in
a series of papers by Christol and Mebkhout
\cite{cm1}, \cite{cm2}, \cite{cm3}, \cite{cm4}. This represents a completely
orthogonal approach to ours, as it primarily involves the structure of the
connection rather than the Frobenius. The different approaches seem to have
different strengths. For example, on one hand, the Christol-Mebkhout approach
seems to say more about $p$-adic differential equations on annuli over
$p$-adic fields which are not discretely valued. On the other hand, our 
approach has a certain flexibility that the Christol-Mebkhout approach lacks;
for instance, it carries over directly to the $q$-difference situation
considered by Andr\'e and di Vizio in \cite{andre}, whereas the analogue of
the Christol-Mebkhout theory seems much more difficult to develop. It also
carries over to the setting of ``fake annuli'' arising in the problem
of semistable reduction for overconvergent $F$-isocrystals: in this setting,
one replaces $k((t))$ by the completion of $k[t_1, t_1^{-1}, \dots, t_n, 
t_n^{-1}]$ for a valuation which totally orders monomials (i.e., the
valuations of $t_1, \dots, t_n$ are linearly independent over $\QQ$).
See \cite{me-fake} for further details.
\end{remark}

\subsection{Generic versus special revisited}
\label{subsec:gen spec2}

The adjectives ``generic'' and ``special'' were introduced in
Chapter~\ref{sec:generic special} to describe the two paradigms
for attaching slopes to $\sigma$-modules over $\Gcon[\fp]$.
Here is a bit of clarification as to why this was done.
Throughout this section, retain notation as in 
Section~\ref{subsec:robba}.

Let $X \to \Spec k \llbracket t \rrbracket$ 
be a smooth proper morphism, for $k$ a field
of characteristic $p>0$. Then the crystalline cohomology
of $X$, equipped with the action of the absolute Frobenius,
gives a $(\sigma,\nabla)$-module $M$ over the ring 
$\calO \llbracket u \rrbracket$; the crystalline cohomology
of the generic fibre corresponds to $M \otimes \Gamma$, whereas
the crystalline cohomology of the special fibre corresponds to
$M \otimes (\calO \llbracket u \rrbracket/u \calO \llbracket u \rrbracket)$.
However, the latter is isomorphic to $M \otimes
(\Gancon \cap \calO[\fp] \llbracket u \rrbracket)$ by ``Dwork's trick''
\cite[Lemma~6.3]{dejong}. Thus the generic and special HN-polygons
correspond precisely to the Newton polygons of the generic and special fibres;
the fact that the special HN-polygon lies above the generic HN-polygon
(i.e., Proposition~\ref{P:comparison}) in this case follows from Grothendieck's
specialization theorem.

This gives a theoretical explanation for why Proposition~\ref{P:comparison}
holds, but a more computationally explicit example may also be useful.
(Thanks to Frans Oort for suggesting this presentation.)
Suppose that $\sigma$ is chosen so that $u^\sigma = u^p$.
Let $M$ be the rank 2 $\sigma$-module over $\Gcon$ defined by
\[
F\bv_1 = \bv_2, \qquad F\bv_2 = p \bv_1 + u \bv_2.
\]
Then $\bv_1$ is a cyclic vector and $F^2 \bv_1 - u F\bv_1 - p\bv_1 = 0$,
so by Lemma~\ref{L:use cyclic vector},
the generic HN-polygon of $M$ has the same slopes as the Newton
polygon of the polynomial $x^2 - ux - p$, namely 0 and 1.
On the other hand, Dwork's trick implies that the special HN-polygon of $M$
has slopes 1/2 and 1/2.

\subsection{Splitting exact sequences (again)}
\label{subsec:split again}

For reference,
we collect here some more results about computing $H^1$ of $\sigma$-modules.
\begin{prop} \label{P:good model ext}
For any $\sigma$-modules $M_1, M_2$ over $\Gcon[\fp]$, the map
$\Ext(M_1, M_2) \to \Ext(M_1 \otimes \Gancon, M_2 \otimes \Gancon)$
is surjective.
\end{prop}
\begin{proof}
Let 
\[
0 \to M_2 \otimes \Gancon \to M \to M_1 \otimes \Gancon \to 0
\]
be a short exact sequence of $\sigma$-modules.
Choose a basis of $M_2$, then lift to $M$ a basis of $M_1$;
the result is a basis of $M$. Let $A$ be the matrix via which
$F$ acts on this basis; after rescaling the basis of $M_2$ suitably,
we can put ourselves into the situation of 
Lemma~\ref{L:frob turnover}. We can now perform the iteration of
Lemma~\ref{L:frob turnover} in such a way as to respect the
short exact sequence (i.e., take $u_{ijlm} = 0$ whenever
the pair $(i,j)$ falls in the lower left block);
as in the proof of Proposition~\ref{P:good model}, we end up
with a model $M_b$ of $M$ over $\Gcon[\fp]$, which by construction
sits in an exact sequence $0 \to M_2 \to M_b \to M_1 \to 0$.
This yields the desired surjectivity.
\end{proof}

We next give some generalizations of parts of
Proposition~\ref{P:h1}.
\begin{prop} \label{P:h1 injective}
Let $M$ be a $\sigma$-module over $\Gcon[\fp]$ whose
generic HN-slopes are all nonnegative. Then the map
$H^1(M) \to H^1(M \otimes \Gancon)$ is injective.
\end{prop}
\begin{proof}
By Proposition~\ref{P:lattice}, we can choose an $F$-stable
$\Gcon$-lattice $M_0$ of $M$.
Let $\be_1, \dots, \be_n$ be a basis of $M_0$,
and define the matrix $A$ over $\Gcon$ by $F\be_j = \sum_i A_{ij} \be_i$.
Choose $r>0$ such that $A$ has entries in $\Gamma_r$ and $w_r(A) \geq 0$.

Suppose $\bv \in M$ and $\bw \in M \otimes \Gancon$ satisfy
$\bv = \bw - F\bw$, and write $\bv = \sum_i x_i \be_i$ and 
$\bw = \sum_i y_i \be_i$ with $x_i \in \Gcon[\fp]$ and 
$y_i \in \Gancon$; then $x_i = y_i - \sum_j A_{ij} y_j^\sigma$.

If $\bw \notin M$, we can choose $m<0$ such that
$v_m(x_i) = \infty$ and
$0 < \min_i \min_{l \leq m} \{v_{l,r}(y_i)\} < \infty$. Then
\begin{align*}
\min_{l \leq m} \{ v_{l,r}(y_i)\} &>
q^{-1} \min_{l \leq m} \{ v_{l,r}(y_i)\} \\
&\geq q^{-1} \min_{l \leq m} \min_j \{ v_{l,r}(A_{ij} y_j^\sigma) \} \\
&\geq q^{-1} \min_{l \leq m} \min_j \{ v_{l,r}(y_j^\sigma) \} \\
&= q^{-1} \min_{l \leq m} \min_j \{ r v_l(y_j^\sigma) + l \} \\
&\geq \min_{l \leq m} \min_j \{ r v_l(y_j) + l \} \\
&= \min_{l \leq m} \min_j \{ v_{l,r}(y_j)\};
\end{align*}
taking the minimum over all $i$ yields a contradiction. Hence
$\bw \in M$, yielding the injectivity of the map
$H^1(M) \to H^1(M \otimes \Gancon)$.
\end{proof}

\begin{prop} \label{P:h0h1 vanishes}
Let $M$ be a $\sigma$-module over $\Gamma[\fp]$
whose HN-slopes are all positive. Then $F-1$ is a bijection
on $M$, i.e., $H^0(M) = H^1(M) = 0$.
\end{prop}
\begin{proof}
By Proposition~\ref{P:lattice}, we can choose an $F$-stable
$\Gamma$-lattice $M_0$ of $M$ such that $F(M_0) \subseteq \pi M_0$.
Then $F-1$ is a bijection on $M_0/\pi M_0$, hence also on $M_0$
and $M$.
\end{proof}

\begin{lemma} \label{L:Galois}
For $L/K$ a finite Galois extension, and $M$ a 
$\sigma$-module over $\Gancon$, the map
$H^1(M) \to H^1(M \otimes \GLancon)$ is injective.
\end{lemma}
\begin{proof}
Put $G = \Gal(L/K)$.
Given $\bw \in M$, suppose that there exists $\bv \in M
\otimes \GLancon$ such that $\bw = \bv - F\bv$. Then we also have
$\bw = \bv' - F\bv'$ for
\[
\bv' = \frac{1}{\#G} \sum_{g \in G} \bv^g,
\]
which is $G$-invariant and hence belongs to $M$.
\end{proof}

\begin{theorem} \label{T:split slope zero}
Let $M$ be a $\sigma$-module over $\Gancon$ whose
HN-slopes are all positive. Then the exact sequence
\begin{equation} \label{eq:split slope}
0 \to M \to N \to \Gancon \to 0
\end{equation}
splits if and only if $N$
has smallest HN-slope zero.
\end{theorem}
\begin{proof}
If the sequence splits, then $P(N) = P(M) + P(\Gancon)$ by
Proposition~\ref{P:split exact}, so $N$ has smallest HN-slope zero.
We verify the converse first in the case where $K$ is algebraically closed,
so that $\Gancon = \Galgancon$.

We proceed by induction on $\rank(M)$.
In case $M$ is isoclinic, then the inequality
$P(N) \geq P(M) + P(\Gancon)$ from Proposition~\ref{P:split exact}
and the hypothesis that $N$ has smallest HN-slope zero
together
force $P(N) = P(M) + P(\Gancon)$; by Proposition~\ref{P:split exact}
again, \eqref{eq:split slope} splits.
In case $M$ is not isoclinic,
let $M_1$ be the first step in the HN-filtration of $M$.
By Proposition~\ref{P:split exact}, we have
\[
P(N) \geq P(M_1) + P(N/M_1) \geq P(M_1) + P(M/M_1) + P(\Gancon);
\] 
since $P(N)$ has smallest slope 0, so does $P(M_1) + P(N/M_1)$. Since
$P(M_1)$ has positive slope, $P(N/M_1)$ must have smallest slope zero.
Hence by the induction hypothesis, the exact sequence $0 \to M/M_1 \to N/M_1
\to \Gancon \to 0$ splits, say as $N/M_1 \cong
M/M_1 \oplus M'$ with $M' \cong \Gancon$.
Let $N'$ be the inverse image of $M'$ under the surjection
$N \to N/M_1$; it follows now that \eqref{eq:split slope}
splits if and only if $0 \to M_1 \to N' \to M' \to 0$ splits.
By Proposition~\ref{P:split exact} again,
$P(N) \geq P(N') + P(M/M_1) \geq P(M_1) + P(M/M_1) + P(\Gancon)$,
and $P(M/M_1)$ has all slopes positive, so $P(N')$ has smallest
slope zero. Again by the induction hypothesis,
$0 \to M_1 \to N' \to M' \to 0$ splits, yielding the splitting of
\eqref{eq:split slope}.

To summarize, we have proved that if $N$ has smallest HN-slope
zero, then \eqref{eq:split slope}
splits after tensoring with $\Galgancon$; it remains to descend this splitting back to $\Gancon$. To do this,
apply Proposition~\ref{P:good model} to produce a $\sigma$-module
$M_0$ over $\GLcon[\fp]$, for some finite Galois extension $L$ of 
$K$, with $M_0 \otimes \GLancon \cong M \otimes \GLancon$.
Then apply Proposition~\ref{P:good model ext} to descend the given
exact sequence, after tensoring up to $\GLancon$, to an exact sequence $0 \to M_0 \to N_0 \to \GLcon[\fp] \to 0$.
We have shown that the exact sequence
$0 \to M_0 \otimes \Galgancon \to N_0 \otimes \Galgancon \to \Galgancon \to 0$
splits; by Proposition~\ref{P:h1 injective}, the exact sequence
$0 \to M_0 \otimes \Galgcon[\fp] \to N_0 \otimes 
\Galgcon[\fp] \to \Galgcon[\fp] \to 0$ splits.

Choose $\bw \in M_0$ whose image in $H^1(M_0)$ corresponds to the
exact sequence $0 \to M_0 \to N_0 \to \GLcon[\fp] \to 0$;
we have now shown that 
the class of $\bw$ in $H^1(M_0 \otimes \Galgcon[\fp])$ vanishes.
By Proposition~\ref{P:h0h1 vanishes}, there exists a unique
$\bv \in M_0 \otimes \Gamma^L[\fp]$ with $\bw = \bv - F\bv$.
Since the class of $\bw$ in $H^1(M_0 \otimes \Galgcon[\fp])$ vanishes,
we have $\bv \in M_0 \otimes \Galgcon[\fp]$; since $M_0$
is free over $\GLcon[\fp]$ and $\Galgcon[\fp] \cap \Gamma^L[\fp] = 
\GLcon[\fp]$, we have $\bv \in M_0$. Hence the
sequence \eqref{eq:split slope} splits after tensoring with
$\GLancon$. By Lemma~\ref{L:Galois}, \eqref{eq:split slope} also
splits, as desired.
\end{proof}

\subsection{Full faithfulness}
\label{subsec:fullfaith}

Here is de Jong's original application of
the reverse filtration \cite[Proposition~8.2]{dejong}.
\begin{prop} \label{P:split reverse}
Suppose that $K$ admits a valuation $p$-basis.
Let $M$ be a $\sigma$-module over $\Gcon[\fp]$ admitting an
injective $F$-equivariant $\Gcon[\fp]$-linear morphism
$\phi: M \to \Gamma[\fp](m)$ for some integer $m$.
Then $\phi^{-1}(\Gcon[\fp])$ 
is a  $\sigma$-submodule of $M$ of rank $1$,
and its extension to $\Galgcon[\fp]$ 
is equal to the first step of the reverse
filtration of $M$. In particular, $M$ 
has highest generic slope $m$ with multiplicity $1$.
\end{prop}
\begin{proof}
We first suppose $K$ is algebraically closed (this case
being \cite[Corollary~5.7]{dejong}).
Let $M_1$ be the first step of the reverse filtration of $M$.
Then $M_1 \otimes \Galg[\fp]$ is isomorphic to a direct sum of standard
$\sigma$-modules of some slope $s_1 = c/d$; by Lemma~\ref{L:negative eigen}
(applied to $([d]_* \Galgcon[\fp])(-c) \otimes M_1$), $M_1$ itself is 
isomorphic
to a direct sum of standard $\sigma$-modules of slope $s_1$.

The map $\phi$ induces a nonzero $F$-equivariant map $M_1 \otimes
\Galg[\fp] \to \Galg[\fp](m)$, and hence a nonzero element of 
$H^0(M_1^\dual \otimes \Galg[\fp](m))$. By
Corollary~\ref{C:hom ext2}, we must have $m = s_1$, and so
$M_1 \cong \Galgcon[\fp](m)^{\oplus n}$ for some $n$.
By Proposition~\ref{P:h0}, we have
$H^0(M_1^\dual(m)) = H^0(M_1^\dual \otimes \Galg[\fp](m))$,
and so $\phi$ actually induces an injective $F$-equivariant map
$M_1 \to \Galgcon[\fp](m)$.

To summarize, $M$ has highest generic slope $m$, and
the first step of the reverse filtration is contained in
$\phi^{-1}(\Galgcon[\fp])$. Since the latter is a $\sigma$-submodule of $M$
of rank no more than 1, we have the desired result.

We now suppose $K$ is general. Put $M' = M \otimes_{\Gcon[\fp]} \Galgcon[\fp]$;
by Proposition~\ref{P:injective}(c),
the composite map
\[
\psi: M' \stackrel{\phi \otimes 1}{\to} \Gamma[\fp] \otimes_{\Gcon[\fp]}
\Galgcon[\fp] \stackrel{\mu}{\to} \Galg[\fp]
\]
is also injective.
By the above, $\psi^{-1}(\Galgcon[\fp])$ is a $\sigma$-submodule of $M'$
of rank 1, and coincides with the first step of the reverse filtration of 
$M'$. Let $\be_1, \dots, \be_n$ be a basis of $M$, put $\bv = 
\psi^{-1}(1)$, and write $\bv = \sum x_i \be_i$ with $x_i \in \Galgcon[\fp]$.
Then $1 = \psi(\bv) = \sum x_i \phi(\be_i)$; by
Proposition~\ref{P:injective}(b), we have $x_i \in \Gcon[\fp]$ for
$i=1, \dots, n$. Hence $\bv \in \phi^{-1}(\Gcon[\fp])$, so the latter
is a $\sigma$-submodule of $M$ of rank 1. This yields the desired result.
\end{proof}

\begin{remark}
Proposition~\ref{P:split reverse} can be used to reduce instances of
showing $H^0(M) = H^0(M \otimes \Gamma[\fp])$, for $M$ an
$F$-module over $\Gcon[\fp]$, to showing that a certain class in
$H^1(N)$, where $N$ is a related $F$-module over $\Gcon[\fp]$ with
positive HN-slopes, vanishes in $H^1(N \otimes \Gancon)$.
Thanks to Proposition~\ref{P:h1}(b), this
can in turn be checked over $\Gancon$,
where either Dwork's trick, in the case of \cite{dejong},
or the $p$-adic local monodromy theorem, in the case of \cite{me-full},
can be brought to bear. Note that by
Theorem~\ref{T:split slope zero}, it is enough to check vanishing
of a class in $H^1(N \otimes \Gancon)$ after replacing $K$ by
a finite separable extension.
\end{remark}

\end{document}